\documentclass[10pt,a4paper,twoside]{article}

\usepackage[utf8]{inputenc}
\usepackage[T1]{fontenc}
\usepackage[english]{babel}
\usepackage{enumerate}
\usepackage{amsmath,amsfonts,amssymb,amsthm,mathtools}
\usepackage{graphicx,multirow}
\usepackage{fancyref} 
\usepackage[hidelinks]{hyperref}
\usepackage{float}
\usepackage{longtable}
\usepackage{colortbl}
\usepackage[table,svgnames]{xcolor}
\usepackage{comment}
\usepackage{dsfont} 
\usepackage{stmaryrd}
\usepackage{calc}
\usepackage{relsize} 
\usepackage{txfonts} 
\usepackage{tikz,pgf}
\usetikzlibrary{arrows,chains,matrix,positioning,scopes}

\usepackage{natbib} 

\usepackage[left=2.5cm,right=2.5cm,top=2.5cm,bottom=2.5cm]{geometry} 

\usepackage{enumitem}

\newcommand\given[1][]{\:#1\vert\:} 
\newcommand\numberthis{\addtocounter{equation}{1}\tag{\theequation}}

\DeclareMathOperator*{\argmin}{argmin}
\DeclareMathOperator*{\argmax}{argmax}

\newtheorem{lemme}{Lemma}
\newtheorem{theorem}{Theorem}
\newtheorem*{theorem*}{Theorem}
\newtheorem{corollary}{Corollary}
\newtheorem{remark}{Remark}

\newcommand{\proba}{\ensuremath{\mathbb{P}}}
\newcommand{\esp}{\ensuremath{\mathbb{E}}}
\newcommand{\probat}{\ensuremath{\mathbb{P}_{\theta}}}
\newcommand{\espt}{\ensuremath{\esp_{\theta}}}
\newcommand{\probaetoile}{\ensuremath{\mathbb{P}^*_{\theta^*}}}
\newcommand{\espetoile}{\ensuremath{\esp_{\theta^*}}}
\newcommand{\var}{\mathrm{Var}}

\newcommand{\Xtij}{\ensuremath{X^t_{ij}}}
\newcommand{\ind}{\ensuremath{\mathds{1}}}
\newcommand{\Q}{\ensuremath{\llbracket 1, Q \rrbracket}}
\newcommand{\T}{\ensuremath{\llbracket 1, T \rrbracket}}
\newcommand{\n}{\ensuremath{\llbracket 1, n \rrbracket}}

\title{Consistency of the maximum likelihood and variational estimators in a dynamic stochastic block model}
\author{Léa Longepierre and Catherine Matias\\
Sorbonne Université, Université Paris
Diderot, Centre National de la Recherche Scientifique, \\
Laboratoire de  Probabilités, Statistique et Modélisation,  \\
4 place Jussieu, 75252 PARIS Cedex 05, FRANCE. \\
\texttt{\{lea.longepierre,catherine.matias\}@sorbonne-universite.fr}
}
\date{}

\begin{document}
\DeclareGraphicsExtensions{.pdf, .jpg, .jpeg, .png, .gif,.eps}
\thispagestyle{empty}
\maketitle

\begin{abstract}
  We consider a dynamic version of the stochastic block model, in which the nodes are partitioned into latent classes
  and the connection between two nodes is drawn from a Bernoulli distribution depending on the classes of these two
  nodes. The temporal evolution is modeled through a hidden Markov chain on the nodes memberships. We
  prove the consistency (as the number of nodes and time steps increase) of the maximum likelihood and variational
  estimators of the model parameters, and obtain upper bounds on the rates of convergence of these estimators.
  We also explore the particular case where the number of time steps is fixed and connectivity parameters are allowed
  to vary.
\end{abstract}


\section{Introduction}
Random graphs are a suitable tool to model and describe interactions in many kinds of datasets such as
biological, ecological, social or transport networks. Here we are interested in time-evolving networks, which is a
powerful tool for modeling real-world phenomena, where the role or behaviour of the nodes in the network and the
relationships between them are allowed to change over time. Indeed, it is important to take into account the
evolutionary behaviour of the graphs, instead of just studying separate snapshots as static graphs. 
We focus on graphs evolving in discrete time and refer to  \cite{holme2015modern} for an introduction to dynamic networks.

A myriad of dynamic graph models has been introduced in the past few years, see for instance \cite{Zhang2017}. We focus here on those which are  based on the (static) stochastic block model~\citep[SBM,][]{HOLLAND1983109} in which the nodes are
partitioned into classes. In the SBM, class memberships of the nodes are represented by latent variables and the connection
between two nodes is drawn from a distribution depending on the classes of these two nodes (a Bernoulli distribution in
the case of binary graphs). A first dynamic version of the SBM with discrete time is proposed in \cite{Yang2011}. There, the
nodes are partitioned into $Q$ classes and the graphs are binary or weighted. The nodes are allowed to change membership
over time, and these changes are governed by independent Markov chains with values in the $Q$ classes, while the
connection probabilities are constant over time.  \cite{xu2014dynamic} introduce a state-space model on the logit of the
connection probabilities for dynamic (binary) networks with connection probabilities and group memberships varying over
time. Unfortunately, their model presents parameter identifiability issues \citep{matias2017statistical}. 
\cite{Xu2015StochasticBT} proposes a stochastic block transition model  in which the presence or absence of an
edge between two nodes at a particular time affects the presence or absence of such an edge at a future time. There, the
nodes can change classes over time, new nodes can enter the network, and the connection probabilities are allowed to
vary over time. The model in \cite{matias2017statistical} and in \cite{Becker2018} is quite similar to that of
\cite{Yang2011} except that it allows the connection probabilities to vary and the latter is moreover nonparametric. \cite{Bartolucci:2018:DRD:3221748.3222094}
extend the model of \cite{Yang2011} to deal with different forms of reciprocity in directed graphs, by directly modeling
dyadic relations and with the assumption that the dyads are conditionally independent given the latent
variables.
\cite{paul2016} and \cite{han2015consistent} study multi-graph SBM, arising in settings including dynamic
networks and multi-layer networks where each layer corresponds to a type of edge.
In these two models, the nodes memberships stay constant over the layers.
\cite{pensky2016dynamic,pensky2017spectral}  study a dynamic SBM  for undirected and binary edges where both connection probabilities and group memberships vary over time, assuming that the connection probabilities between groups are a smooth function of time.
\cite{xing2010} and \cite{ho2011evolving} introduce dynamic versions of the mixed-membership
stochastic block model, allowing each actor to carry out different roles when interacting with different peers.
\cite{zreik:hal-01122393} introduce the dynamic random subgraph model, given a known decomposition of the graph into
subgraphs, in which the latent class membership depends on the subgraph membership and the edges are categorical
variables, their types being sampled from a distribution depending on the latent classes of the two nodes. There, a
state-space model is used to  characterize the temporal evolution of the latent classes
proportions.

As far as estimation is concerned, different methods of inference are proposed to estimate groups and model parameters.
The maximum likelihood estimator (MLE)  is not tractable in the SBM, thus neither in its dynamic versions. 
Variational methods are rather popular to approximate that MLE
\citep{xing2010,ho2011evolving,han2015consistent,paul2016,zreik:hal-01122393,matias2017statistical,Bartolucci:2018:DRD:3221748.3222094}.
\cite{Yang2011} rely on Gibbs sampling and simulated annealing.
\cite{pensky2017spectral} propose an estimator of the connection probabilities matrix at each time step by a discrete
kernel-type method and obtain a clustering of the nodes thanks to spectral clustering on this estimated matrix. They
also give an estimator for the number of clusters. Spectral clustering algorithms are also used by
\cite{han2015consistent} on the mean graph over time
and by \cite{Liu927} who use eigenvector smoothing to get some similarity across time periods (and allow the number  of
classes to be unknown and possibly varying over time).

Some theoretical results on the convergence of the procedures have been proven, mainly for static graphs. In the static SBM, \cite{celisse2012consistency} prove the consistency of the MLE and variational estimates as the number of nodes increases, and \cite{bickel2013} establish their asymptotic normality. \cite{mariadassou2015} have a different approach and give sufficient conditions for the groups posterior distribution to converge to a Dirac mass located at the actual groups configuration, for every parameter in a neighborhood of the true one. \cite{rohe2011} give asymptotic results on the normalized graph Laplacian and its eigenvectors for the spectral clustering algorithm, allowing the number of clusters to grow with the number of nodes. They also provide bounds on the number of misclustered nodes, requiring an assumption on the degree distribution. \cite{lei2015} prove consistency for the recovery of communities in the spectral clustering on the adjacency matrix, with milder conditions on the degrees, and also extend this result to degree corrected stochastic block models. 
\cite{klopp2017} derive oracle inequalities for the connection probabilities estimator and obtain minimax estimation
rates, including the sparse case where the density of edges converges to zero as the number of nodes increase thus
extending previous results of \cite{gao2015}. \cite{Gaucher_2019} propose a bound on the risk of the
maximum likelihood estimator of network connection probabilities, and show that it is minimax optimal in the sparse
graphon model. 

In the dynamic setting, fewer theoretical results have been established. 
\cite{pensky2016dynamic} derives a penalized least squares estimator of the connection probabilities adaptive to
the number of blocks and which does not require knowledge of the number of classes $Q$.
She shows that it satisfies an oracle inequality. Under the additional assumption that at most $n_0$ nodes change groups
between two time steps, this 
estimator attains minimax lower bounds for the risk. She also introduces a dynamic graphon model and shows that the
estimators (that do not require knowledge of a degree of smoothness of the graphon function) are minimax optimal within a logarithmic factor of the number of time steps.
Based on the same dynamic SBM with at most $n_0$ nodes changing groups between two time steps, \cite{pensky2017spectral}
give an upper bound for the (non asymptotic) error of their estimators of the connection probabilities matrix and group memberships (and also an
estimator for the number of clusters). \cite{han2015consistent} show consistency (as the number of time steps increases
but the number of nodes is fixed) of two estimators of the class memberships for dynamic SBM (and more generally
multi-graph SBM) in which the nodes memberships are constant over time but the connection probabilities are allowed to
vary and the considered graphs are binary and symmetric. They show that the spectral clustering (on the mean graph over
time) estimator of the class memberships is consistent under some stationarity and ergodicity conditions on the connection probabilities. They
also prove that the MLE of the class memberships is consistent (i.e. that the fraction of misclustered nodes converges to $0$) in the general case (without any structure on the connection probabilities),
provided certain sufficient conditions are satisfied. 
In their multi-layer model, \cite{paul2016} give minimax rates of misclassification  under certain conditions on the
growth of the types of relations, number of nodes and number of classes, extending the result of \cite{han2015consistent}.

Here, we consider a dynamic version of the binary SBM as in \cite{Yang2011}, where each node is allowed to change group
membership at each time step according to a Markov chain, independently of other nodes. 
We prove the consistency of the connectivity parameter MLE and, under some additional conditions, of the transition
matrix MLE, when the number of nodes and of time steps are increasing. We also give upper bounds on the rates of
convergence of these estimators. While these upper bounds are known to be non optimal in the static case where
asymptotic normality is obtained with classical parametric rates of convergence \citep{bickel2013}, these are the first to be
established in a dynamic setting for the MLE. 
As already mentioned,  the log-likelihood is intractable (except for very small values of the number of nodes $n$ and
the number of time steps $T$), as it requires to sum over 
$Q^{nT}$ terms. Thus, while its consistency remains an important result, the estimator cannot be computed.
A possible alternative is to rely on a variational estimator to approximate the MLE~\cite[see for instance][]{matias2017statistical}. We also establish the consistency of the variational
estimator of the connectivity parameter and under some additional assumptions, that of the variational estimator of the transition matrix and obtain the same upper bounds on the rates of convergence as for the MLE.
In the particular case where the number of time steps $T$ is fixed, we also consider the model of
\cite{matias2017statistical}, in which the connection probabilities are allowed to vary over time and generalise these
results with only the number of nodes increasing. When $T=1$, we not only recover the results of
\cite{celisse2012consistency} but extend these by giving rates of convergence. Unlike the model studied in \cite{han2015consistent} and \cite{paul2016}, the node memberships in our model evolve over time. 
Our context is different from \cite{pensky2016dynamic} that focuses on least squares estimate.

This article is organized as follows. Section~\ref{sec:model} introduces our model and notation. More precisely,
Section~\ref{sub:model} describes the dynamic stochastic block model as introduced in \cite{Yang2011},
Section~\ref{sub:assumptions} gives the assumptions we make on the model parameters, Section~\ref{sub:finite} describes
the dynamic stochastic block model as in \cite{matias2017statistical} for the finite time case and
Section~\ref{sub:likelihood} states the expression of the likelihood of this model to define the
MLE. Section~\ref{sec:consistency_MLE} establishes the consistency and upper bounds of the rates of convergence for the
MLE of the connection probabilities in Section~\ref{sub:consistency_pi} and of the transition matrix in
Section~\ref{sub:transition_matrix}. Section~\ref{sec:variational} is dedicated to variational estimators:
Section~\ref{sub:consistency_pi_variational} and \ref{sub:consistency_gamma_variational} establish  the consistency of
the variational estimators of the connection probabilities and transition matrix, respectively, along with upper bounds
of the associated rates of convergence. All the proofs of the main results are postponed to
Section~\ref{sec:proofs_main}, except those for the fixed $T$ case that are in Appendix~\ref{appendix:proofs_fixed_T}, while the more technical proofs are deferred to Appendix~\ref{appendix:proofs_tech}.

\section{Model and notation} \label{sec:model}
\subsection{Dynamic stochastic block model} \label{sub:model}
We consider a set of $n$ vertices, forming a sequence of binary undirected graphs with no self-loops at each time $t=1,\ldots,T$. The case of a set of directed graphs, with or without self-loops, may be handled similarly.  These vertices are
assumed to be split into $Q$ latent classes, and we denote by $Z^t_i$ the label of the $i$-th vertex at time $t$.
Letting $Z_i=(Z_i^1,\dots,Z_i^T)$, we assume that the $\{Z_i\}_{1\le i \le n}$ are independent and
identically distributed (iid) and each
$Z_i$ is a homogeneous and stationary Markov chain with transition probabilities
\[
\proba(Z_i^{t+1}=l\given Z_i^{t}=q)=\gamma_{ql} , \quad \forall 1\le q,l\le Q
\]
where  $\Gamma=(\gamma_{ql})_{1\leq q,l \leq Q}$ is a stochastic matrix, i.e. with nonnegative coefficients and with each row summing to 1.
We let 
$\alpha=(\alpha_1,\ldots,\alpha_Q)$  the stationary distribution of the Markov chain. For any $i\in \n$, the probability distribution of $Z_i$ is then
\begin{align*}
\proba_{\theta}(Z_i)= \alpha_{Z_i^1} \prod_{t=1}^{T-1} \gamma_{Z_i^t Z_i^{t+1}}.
\end{align*}
We will also denote $Z^t =(Z^t_1,\dots, Z^t_n)$ and $Z^{1:T}=(Z^1,\dots,Z^T)=(Z^t_i)_{1\le t\le T, 1\le i \le n}$.

Consider $X^t=\{\Xtij\}_{1\leq i,j\leq n}$ the symmetric binary adjacency matrix of the graph at time $t$
such that for every nodes $1\le i,j\le n$, we have $X^t_{ii}=0$ and $\Xtij=X^t_{ji}$.
Each $X^t$ follows a stochastic block model so that, conditional on the latent groups $\{Z^t_i\}_{1\le i \le n}$, the
$\{\Xtij\}_{1\leq i,j\leq n}$ are independent Bernoulli random variables  
\[
\Xtij \given Z^t_i=q, Z^t_j=l  \sim \mathcal{B}(\pi_{ql})
\]
where 
$(\pi_{ql})_{1\leq q,l\leq Q} \in [0,1]^{Q^2}$ are the connectivity parameters. More precisely, conditional on the whole sequence of latent groups $\{Z^t_i\}_{1\le t\le T, 1\le i \le n}$, the  graphs
$X^{1:T}=X^1,\dots, X^T$ are assumed to be independent, each $X^t$ having a distribution depending only on
$\{Z^t_i\}_{1\le i \le n}$. 
The model  is thus parameterized by $\theta=(\Gamma,\pi)$, with $\Gamma=(\gamma_{ql})_{1\leq q,l\leq Q}$ and
$\pi=(\pi_{ql})_{1\leq q,l\leq Q}$. Note that $\pi$ is a symmetric matrix in the undirected setup. We denote by $\probat$ (resp.  $\espt$) the probability distribution (resp. expectation) of all the random variables $\{Z_i^t, \Xtij\}_{t\ge 1 ; i,j\ge 1}$, under the parameter value $\theta$. In the following, we assume that we observe $\{\Xtij\}_{1\leq i,j, \leq n,\; 1\leq t \leq T}$ and we denote by $\theta^*=(\Gamma^*,\pi^*)=((\gamma^{*}_{ql})_{1\leq q,l\leq Q},
(\pi^{*}_{ql})_{1\leq q,l\leq Q})$ the true parameter value, with corresponding probability distribution $\proba_{\theta^*}$ and expectation $\espetoile$, and by $\alpha^*=(\alpha_q^*)_{1\leq q\leq Q}$ the (true) stationary distribution 
corresponding to the transition matrix $\Gamma^*$. 
We also let $\ind_A$ denote the indicator function of the set $A$ and $A^c$ the complementary set of $A$ in the ambient set. For any integer $M\ge 1$, the set $\llbracket 1, M \rrbracket$ is the set of integers between $1$ and $M$. For any finite set $A$, let $|A|$ denote its cardinality. 
For any configuration $z^{1:T}$, we  denote  $N_{q}(z^t)$ (resp. $N_{ql}(z^{1:T})$)  the number of nodes assigned to class $q$ by the configuration $z^t$ (resp. the number of transitions from class $q$ to class $l$ in configuration $z^{1:T}$), that is 
  \begin{align}
N_{q}(z^t)=|\{i \in \n ;z^t_i=q\}| \quad \text{ and } \quad 
    N_{ql}(z^{1:T})=\sum_{t=1}^{T-1} \sum_{i=1}^n \ind_{z_i^{t}=q,z_i^{t+1}=l}  . \label{eq:def_nql}
\end{align}
We also define for any two parameters $\theta=(\Gamma,\pi)$ and $\theta'=(\Gamma',\pi')$ the following distances
\begin{align*}
\| \pi - \pi'\|_{\infty} = \max_{ 1\leq q,l \leq Q}  |\pi_{ql}-\pi'_{ql}| \quad \text{ and } \quad 
\| \Gamma - \Gamma'\|_{\infty}= \max_{ 1\leq q,l \leq Q}  |\gamma_{ql}-\gamma'_{ql}|.
\end{align*}

\subsection{Assumptions} \label{sub:assumptions}
The assumptions we make on the model parameters are the following. 
\begin{enumerate}
\item \label{itm:hyp_id} For every $1\le q\neq q'\le Q$, there exists some $l\in \Q$ such that $\pi_{ql} \neq \pi_{q'l}$. 
\item \label{itm:hyp_gamma}
There exists some $0<\delta< 1/Q$ such that for any $ (q,l) \in \Q^2$, we have $ \gamma_{ql} \in [\delta,1-\delta] $. 
\item \label{itm:hyp_pi} There exists some $\zeta>0$ such that for any $(q,l) \in \Q^2$, we have $\pi_{ql} \in [\zeta,1-\zeta] $.
\end{enumerate}
      
Assumption~\ref{itm:hyp_id} is necessary for identifiability of the model. Indeed, if it does not hold, we cannot distinguish between classes $q$ and $q'$.
Assumption~\ref{itm:hyp_gamma} ensures that each Markov chain $Z_i$ is irreducible, aperiodic and recurrent. This
assumption could be weakened at the cost of technicalities. In particular, it implies that the stationary distribution
$\alpha$ exists. Moreover, Assumption~\ref{itm:hyp_gamma} also implies that for any $q\in \Q$, we have $\alpha_q \in [\delta,1-\delta]$. Note that this can be seen as an equivalent of Assumption 2 in \cite{celisse2012consistency} (on the probability distribution of the class memberships) in the dynamic case. \cite{celisse2012consistency} however also have an additional assumption that is an empirical version of this assumption (which states that the observed class proportions are bounded away from $0$) that is true with high probability. We do not make such an assumption and use the fact that the probability of this event converges to $1$.
Assumption~\ref{itm:hyp_pi} is technical and could also be weakened with additional technicalities. For example, \cite{celisse2012consistency} also consider the case $\pi_{ql} \in \{0,1\}$ (i.e. $\pi_{ql} \in \{0,1\} \cup [\zeta,1-\zeta]$) whereas we do not. 
The whole parameter set defined by these constraints is denoted by $\Theta$. In the following, we assume that $\theta^* \in \Theta$.

In what follows, we work up to label permutation on the groups.
Indeed, as in any latent group model, the parameters can only be recovered up to label switching on the
latent groups. We then define the following notation for any permutation $\sigma \in \mathfrak{S}_Q$ with $\mathfrak{S}_Q$ the set of permutations on $\Q$
\[
\theta_{\sigma}=(\Gamma_{\sigma},\pi_{\sigma})
=\left((\gamma_{\sigma(q)\sigma(l)})_{1\leq q,l \leq Q},(\pi_{\sigma(q)\sigma(l)})_{1\leq q,l \leq Q}\right).
\]
\subsection{Finite time case} \label{sub:finite}
If the number of time steps $T$ is fixed, it is possible to let the connection probabilities vary over time. We then consider this case, the connection parameter now being $\pi^{1:T}=(\pi^1,\ldots,\pi^T)$ with $\pi^t=(\pi^t_{ql})_{1 \leq q,l \leq Q}$ for every $t \in \T$ and $\pi_{ql}^t=\probat(\Xtij=1 \given Z^t_i=q, Z^t_j=l)$ for any $(t,q,l) \in \T \times \Q^2$. Note that this is the more general model of \cite{matias2017statistical}, in which the model parameter is $\theta=(\Gamma,\pi^{1:T})$.
Moreover, we introduce the following Assumptions~\ref{itm:hyp_id_fixed_T} and~\ref{itm:hyp_pi_fixed_T} that are alternate versions of Assumptions~\ref{itm:hyp_id} and~\ref{itm:hyp_pi} respectively for the finite time case.
\begin{enumerate}[label=\arabic*']
	\item\hspace{-7pt}. \label{itm:hyp_id_fixed_T} For every $t \in \T$, for every $1\le q\neq q'\le Q$, there exists some $l\in \Q$ such that $\pi^t_{ql} \neq \pi^t_{q'l}$.
	\setcounter{enumi}{2}
	\item\hspace{-7pt}. \label{itm:hyp_pi_fixed_T} There exists some $\zeta>0$ such that for every $t \in \T$, for any $(q,l) \in \Q^2$, we have $\pi^t_{ql} \in [\zeta,1-\zeta] $.
\end{enumerate}
Assumption~\ref{itm:hyp_id_fixed_T} (resp. Assumption~\ref{itm:hyp_pi_fixed_T}) expresses that for every $t \in \T$, $\pi^t$ satisfies Assumption~\ref{itm:hyp_id} (resp. Assumption~\ref{itm:hyp_pi}).
We also introduce the following additional assumption, which ensures (together with Assumption~\ref{itm:hyp_id_fixed_T}) that the model is identifiable (up to a label permutation). See \cite{matias2017statistical}.
\begin{enumerate} \setcounter{enumi}{3}
	\item \label{itm:hyp_pi_variation} For every $q \in \Q$, for every $t_1, t_2 \in \T$, \, $\pi^{t_1}_{qq}=\pi^{t_2}_{qq} \coloneqq \pi_{qq}$ and $\{\pi_{qq}; q\in \Q \}$ are $Q$ distinct values.
\end{enumerate}
Assumption~\ref{itm:hyp_pi_variation} states that the diagonal of $\pi$ does not change over time, and that its values are distinct. We denote by $\Theta^T$ the set of parameters satisfying Assumptions~\ref{itm:hyp_id_fixed_T}, \ref{itm:hyp_gamma}, \ref{itm:hyp_pi_fixed_T} and~\ref{itm:hyp_pi_variation}. As before, we assume in the following that $\theta^* \in \Theta^T$ in the fixed $T$ case. We also define as before for any $\pi^{1:T}$ and $\pi'^{1:T}$ the distance
\begin{align*}
\| \pi^{1:T} - \pi'^{1:T} \|_{\infty}= \max_{ (q,l,t) \in \Q^2 \times \T}  |\pi_{ql}^t-\pi'^t_{ql}|.
\end{align*}

\subsection{Likelihood} \label{sub:likelihood}
The conditional log-likelihood and the log-likelihood write 
\begin{align*}
\ell_{c} (\theta ; Z^{1:T})&=
\log\probat(X^{1:T} \given Z^{1:T})=\sum_{t=1}^T \log\probat(X^{t} \given Z^{t})
=\sum_{t=1}^T \sum_{1\le i<j \le n} \Xtij \log \pi_{Z^{t}_i Z^{t}_j} + (1-\Xtij) \log (1-\pi_{Z^{t}_i Z^{t}_j}) \\
\text{and} \quad  \ell(\theta)&=
\log\probat(X^{1:T})=\log \left( \sum_{z^{1:T}\in \Q^{nT}} e^{\ell_c(\theta; z^{1:T})}\probat(Z^{1:T}=z^{1:T})\right) , \numberthis \label{eq:log-likelihood}
\end{align*}
respectively. 
We then denote the maximum likelihood estimator (MLE) by
\[
  \hat{\theta}=(\hat{\Gamma},\hat{\pi})=\argmax_{\theta \in \Theta} \ell(\theta).
\]
In the next section, we study separately the consistency of the connectivity parameter estimator $\hat  \pi$ and that of the transition matrix estimator $\hat \Gamma$.

\section{Consistency of the maximum likelihood estimate} \label{sec:consistency_MLE}
\subsection{Connectivity parameter} \label{sub:consistency_pi}
We first prove the consistency of the maximum likelihood estimator of the connectivity parameter $\pi=(\pi_{ql})_{1 \leq q,l \leq Q}$ when the number of nodes and time steps increase. We denote the normalized log-likelihood by 
\[
M_{n,T}(\Gamma,\pi) = \frac{2}{n(n-1)T} \ell(\theta)= \frac{2}{n(n-1)T}\log\probat(X^{1:T})
\]
and introduce the quantities, for any $A=(a_{ql})_{1 \leq q,l \leq Q} \in \mathcal{A}$ the set of $Q \times Q$ stochastic matrices,
\begin{align} \label{eq:def_M}
  \mathbb{M}(\pi,A) &= \sum_{1\le q,l \le Q} \alpha^*_q \alpha^*_l \sum_{1\le q',l'\le Q} a_{qq'} a_{ll'} [ \pi_{ql}^{*} \log\pi_{q'l'} +
  (1-\pi_{ql}^{*}) \log(1-\pi_{q'l'}) ] \nonumber \\
\text{and} \quad \mathbb{M}(\pi) 
&=\sup_{A \in \mathcal{A}} \mathbb{M}(\pi,A) = \mathbb{M}(\pi,\bar{A}_{\pi}),  
\end{align}
where $\bar{A}_{\pi}=\argmax_{A \in \mathcal{A}}
\mathbb{M}(\pi,A)$. 
It is worth noticing that $\mathbb{M}(\pi)$, which will be the limiting value for $M_{n,T}(\Gamma,\pi)$ when 
$n$ and $T$ increase (see below), does not depend on $\Gamma$.

\begin{theorem}\label{prop:vitesse_cv_M}
For any sequence $\{r_{n,T}\}_{n,T \geq 1}$ increasing to infinity, if $\log(T)=o(n)$, we have for all $\epsilon >0$ 
 \begin{align*}
 \proba_{\theta^*}\left(\sup_{(\Gamma,\pi)\in \Theta}\left| M_{n,T}(\Gamma,\pi) - \mathbb{M}(\pi) \right| >
  \frac{\epsilon r_{n,T}}{\sqrt{n}}\right) \xrightarrow[n,T\to +\infty]{} 0 . 
 \end{align*}
\end{theorem}
We then conclude on the consistency of the maximum likelihood estimator of the connection probabilities with the following corollary. Note that we also obtain an upper bound of the rate of convergence of this estimator.

\begin{corollary}\label{cor:vitesse_cv_pi_bis}
  	For any sequence $\{r_{n,T}\}_{n,T \geq 1}$ increasing to infinity such that $r_{n,T}=o(n^{1/4})$ and if $\log(T)=o(n)$, we have for every $\epsilon>0$
  	\begin{align*}
  	\proba_{\theta^*}\left(\min_{\sigma \in \mathfrak{S}_Q} \| \pi^* - \hat{\pi}_{\sigma} \|_{\infty} > \frac{\epsilon r_{n,T}}{n^{1/4}} \right) \xrightarrow[n,T \rightarrow \infty]{}0.
  	\end{align*}
  \end{corollary}

We want to get equivalent consistency results if the number of time steps $T$ is fixed and only the number of nodes $n$ increases. In that case, denoting by $\hat{\theta}=(\hat{\Gamma},\hat{\pi}^{1:T})$ the MLE of $\theta$, we have the following Corollary that is the equivalent of Corollary~\ref{cor:vitesse_cv_pi_bis}.

  \begin{corollary}\label{cor:vitesse_cv_pi_fixedT}
  	If the number of time steps $T$ is fixed, we have for every $\epsilon>0$ and for any sequence $\{r_n\}_{n \geq 1}$ increasing to infinity such that $r_n=o(n^{1/4})$
  	\begin{align*}
  	\proba_{\theta^*}\left(\min_{\sigma \in \mathfrak{S}_Q }\| \pi^{*1:T} - \hat{\pi}_{\sigma}^{1:T} \|_{\infty} > \frac{\epsilon r_n}{n^{1/4}}\right) \xrightarrow[n \rightarrow \infty]{}0,
  	\end{align*}
	denoting $\hat{\pi}^{1:T}_{\sigma}=(\hat{\pi}_{\sigma}^{t})_{t \in \T}$.
  \end{corollary}
  
This result states that $\min_{\sigma \in \mathfrak{S}_Q } \| \pi^{*1:T} - \hat{\pi}_{\sigma}^{1:T} \|_{\infty}$ converges to $0$ in $\proba_{\theta^*}$-probability as $n$ increases, i.e. the MLE of the connection probabilities is consistent up to label switching, and gives an upper bound of the rate of convergence of the MLE of the connection probabilities. 
The particular case when $T=1$ is then a stronger result than that of \cite{celisse2012consistency} where no rate of convergence is given.

   \begin{remark} \label{rq:cor1_2}
   	Note that in Corollaries~\ref{cor:vitesse_cv_pi_bis} and~\ref{cor:vitesse_cv_pi_fixedT}, the results still hold for any sequences $r_{n,T}$ and $r_n$ increasing to infinity, respectively. However, we are interested in sequences increasing slowly to infinity, giving the strongest results, namely the smallest lower bounds. Indeed, whenever these assumptions are not satisfied, the lower bounds appearing in the inequalities are larger, and the results may even become trivial.
   \end{remark}
   
\subsection{Latent transition matrix}
\label{sub:transition_matrix}
We now prove that the MLE for the transition matrix $\Gamma$ is consistent when the number of nodes and time steps increase.

\begin{lemme} \label{estgamma} 
 Any critical point $\breve{\theta}=(\breve{\Gamma}, \breve \pi)$  of the likelihood function $\ell(\cdot)$ is such that $\breve \Gamma$ satisfies the fixed  point equation 
	\begin{align} \label{eq:criticalpoint}
\forall (q,l)\in \Q^2, \quad \breve{\gamma}_{ql}= \frac{\sum_{t=1}^{T-1} \sum_{i=1}^n
          \proba_{\breve{\theta}}\left(Z^t_i=q, Z^{t+1}_i=l \given X^{1:T}\right)} {\sum_{t=1}^{T-1} \sum_{i=1}^n
          \proba_{\breve {\theta}}\left(Z^t_i=q \given X^{1:T} \right)} .
	\end{align}
      \end{lemme}
There are two different possible cases for the MLE $\hat{\theta}$
\begin{itemize}
	\item Either $\hat{\theta}$ is a critical point of the likelihood function. Then $\hat{\Gamma}$ satisfies equation~\eqref{eq:criticalpoint}.
	\item Or $\hat{\theta}$ is not a critical point (this can happen if it belongs to the boundary of $\Theta$) and
          we assume that there exists $\breve{\Gamma}$ such that $(\breve{\Gamma},\hat{\pi}) \in \Theta$ and
          $(\breve{\Gamma},\hat{\pi})$ satisfies equation~\eqref{eq:criticalpoint}
          (at least for $n$ and $T$ large enough). We then choose as our estimator $(\breve{\Gamma},\hat{\pi})$. By an abuse of notation, we will denote this estimator $\hat{\theta}=(\hat{\Gamma},\hat{\pi})$ and call it MLE in the following.
\end{itemize}

In what follows, for any fixed configuration $z^{1:T}$, any $\theta \in \Theta$ and any $\epsilon >0$, we consider the event 
\[
\mathcal{E}(z^{1:T}, \theta, \epsilon) \coloneqq\left\{ \frac{\proba_{\theta}(Z^{1:T}\neq z^{1:T} \given X^{1:T}
      )}{\proba_{\theta}(Z^{1:T}=z^{1:T} \given X^{1:T} )} > \epsilon \right\} .
  \]
The following result establishes that asymptotically, any estimator that correctly estimates the transition probability
matrix $\pi$ also recovers the group memberships. This result is similar to Theorem 1 in~\cite{mariadassou2015}.

\begin{theorem}
  \label{prop:ratio} 
  For any estimator $\breve{\theta} \in \Theta$ (at least for $n$ and $T$ large enough), if $\log(T)=o(n)$, there exist some positive constants $C, C_1, C_2, C_3, C_4$  such that for any  $\epsilon >0$, for any positive sequence $\{y_{n,T}\}_{n,T \geq 1}$ such that $\log(1/y_{n,T})=o(n)$, any  $\eta \in
  (0,\delta)$  and for $n$ and $T$ large enough, we have
\begin{multline*}
\proba_{\theta^*}\left( \mathcal{E}(Z^{1:T}, \breve{\theta}, \epsilon y_{n,T}) \right)
\leq  
QT\exp(-2\eta^2n) + 
\proba_{\theta^*} \left( \| \breve{\pi} - \pi^* \|_{\infty} > v_{n,T} \right)\\
+C nT \left\{\exp \Bigg[ - (\delta - \eta )^2 C_1 n + C_2 \log(nT) - C_4 \log(\epsilon y_{n,T}) \Bigg] + \exp\Bigg[ - C_3
\frac{(\log(nT))^2}{n v_{n,T}^2}+ 3 n \log(nT)  \Bigg] \right\}, 
\end{multline*}
whenever $\{v_{n,T}\}_{n,T \geq 1}$ is a sequence decreasing to $0$ such that $v_{n,T}=o(\sqrt{\log(nT)}/n)$. 
\end{theorem}

\begin{theorem}\label{prop:Gamma_consistency_bis}
	If $\log(T)=o(n)$, for any $\epsilon>0$ and $\{r_{n,T}\}_{n,T\geq 1}$ any sequence increasing to infinity such that $r_{n,T}=o\left(\sqrt{nT/\log n}\right)$, we have for any $\sigma\in \mathfrak{S}_Q$ 
	\begin{align*}
	&\proba_{\theta^*}\left( \| \hat{\Gamma}_{\sigma}-\Gamma^* \|_{\infty} > \epsilon r_{n,T}\frac{\sqrt{\log n}}{\sqrt{nT}} \right) \leq Q^2(3Q+1)
	\proba_{\theta^*} \left( \| \hat{\pi}_{\sigma} - \pi^* \|_{\infty} > v_{n,T} \right) + o(1)
	\end{align*}
	with $\{v_{n,T}\}_{n,T \geq 1}$ a sequence decreasing to $0$ such that $v_{n,T}=o(\sqrt{\log(nT)}/n)$.
\end{theorem}

\begin{corollary} \label{cor:Gamma_consistency}
	Assume that $\log(T)=o(n)$ and $\min_{\sigma\in \mathfrak{S}_Q}\| \hat{\pi}_{\sigma}- \pi^* \|_{\infty}=o_{\proba_{\theta^*}}(v_{n,T})$ with $\{v_{n,T}\}_{n,T \geq 1}$ a sequence decreasing to $0$ such that $v_{n,T}=o(\sqrt{\log(nT)}/n)$. Then for any $\epsilon>0$ and $\{r_{n,T}\}_{n,T\geq 1}$ any sequence increasing to infinity such that $r_{n,T}=o\left(\sqrt{nT/\log n}\right)$, we have the convergence
	\begin{align*}
	\proba_{\theta^*}\left( \min_{\sigma \in \mathfrak{S}_Q} \| \hat{\Gamma}_{\sigma}-\Gamma^* \|_{\infty} > \epsilon r_{n,T}\frac{\sqrt{\log n}}{\sqrt{nT}} \right) \xrightarrow[n,T \to \infty]{}0.
	\end{align*}
\end{corollary}

\begin{remark}
	Note that the upper bound obtained  in Corollary~\ref{cor:vitesse_cv_pi_bis} on the rate of convergence in
        probability of $\hat{\pi}$ does not ensure that  $\min_{\sigma\in \mathfrak{S}_Q}\| \hat{\pi}_{\sigma}- \pi^*
        \|_{\infty}=o_{\proba_{\theta^*}}(v_{n,T})$ holds. While the latter has never been established (to our
        knowledge), it is a reasonable assumption.  
\end{remark}

We want an equivalent result than that of Corollary~\ref{cor:Gamma_consistency} when the number of time steps $T$ is fixed, and the connection probabilities are varying over time (the connection parameter being $\pi=\pi^{1:T}=(\pi^1,\ldots,\pi^T)$ with $\pi^t=(\pi^t_{ql})_{q,l}$).
For that, we are going to need an equivalent of Theorem~\ref{prop:ratio} in that case.

\begin{theorem} \label{prop:ratio_alt}
	For any fixed $T\geq 2$, for any estimator $\breve{\theta} \in \Theta^T$ (at least for $n$ large enough), there exist some positive constants $C, C_1, C_2, C_3, C_4$  such that for any  $\epsilon >0$, for any positive sequence $\{y_{n}\}_{n\geq 1}$ such that $\log(1/y_{n})=o(n)$, any $\eta \in
	(0,\delta)$ and for $n$ large enough, we have
	\begin{multline*}
	\proba_{\theta^*}\left( \mathcal{E}(Z^{1:T}, \breve{\theta}, \epsilon y_{n}) \right)
	\leq  
	QT\exp(-2\eta^2n) + 
	\proba_{\theta^*} \left( \| \breve{\pi}^{1:T} - \pi^{*1:T} \|_{\infty} > v_{n} \right)\\
	+C nT \left\{\exp \Bigg[ - (\delta - \eta )^2 C_1 n + C_2 \log(nT) - C_4 \log(\epsilon y_{n}) \Bigg] + \exp\Bigg[ - C_3
	\frac{(\log(nT))^2}{n v_{n}^2}+ 5 n \log(nT)  \Bigg] \right\} , 
	\end{multline*}
	whenever $\{v_{n}\}_{n \geq 1}$ is a sequence decreasing to $0$ such that $v_{n}=o(\sqrt{\log(n)}/n)$. 
\end{theorem}

The following corollary gives the expected result.
\begin{corollary} \label{cor:Gamma_consistency_fixedT}
	Let the number of time steps $T \geq 2$ be fixed. Assume that $\min_{\sigma\in \mathfrak{S}_Q}\| \hat{\pi}^{1:T}_{\sigma}- \pi^{*1:T} \|_{\infty}=o_{\proba_{\theta^*}}(v_{n})$ with $\{v_{n}\}_{n \geq 1}$ a sequence decreasing to $0$ such that $v_{n}=o(\sqrt{\log(n)}/n)$. Then for any $\epsilon>0$ and $\{r_{n}\}_{n\geq 1}$ any sequence increasing to infinity such that $r_{n}=o\left(\sqrt{n/\log n}\right)$, we have the convergence
	\begin{align*}
	\proba_{\theta^*}\left( \min_{\sigma \in \mathfrak{S}_Q} \| \hat{\Gamma}_{\sigma}-\Gamma^* \|_{\infty} > \epsilon r_{n} \frac{\sqrt{\log n}}{\sqrt{n}} \right) \xrightarrow[n \to \infty]{}0.
	\end{align*}
\end{corollary}

	The proof of Corollary~\ref{cor:Gamma_consistency_fixedT} is the same as that of Corollary~\ref{cor:Gamma_consistency}, but relying on
        Theorem~\ref{prop:ratio_alt} instead of Theorem~\ref{prop:ratio} and is therefore omitted.
\begin{remark} 
	As in Remark~\ref{rq:cor1_2} for Corollaries~\ref{cor:vitesse_cv_pi_bis} and ~\ref{cor:vitesse_cv_pi_fixedT}, the results of Corollaries~\ref{cor:Gamma_consistency} and~\ref{cor:Gamma_consistency_fixedT}  still hold for sequences $r_{n,T}$ and $r_n$ increasing to infinity at any rate.
\end{remark}

\section{Variational estimators} \label{sec:variational}
In practice, we cannot compute the MLE except for very small values of $n$ and $T$, because it involves a summation over all the $Q^{nT}$ possible latent configurations. We cannot either use the Expectation-Maximization (EM) algorithm to approximate it because it involves the computation of the conditional distribution of the latent variables given the observations which is not tractable. A common solution is to use the Variational Expectation-Maximization (VEM) algorithm that optimizes a lower bound of the log-likelihood (see for example \cite{Daudin2008}).
Let us denote $Z^t_{iq}=\ind_{Z^t_i=q}$ for every $t,i$ and $q$. Using the same approach as in \cite{matias2017statistical} for the VEM algorithm in the dynamic SBM, we consider a variational approximation of the conditional distribution of the latent variable $Z^{1:T}$ given the observed variable $X^{1:T}$ in the class of probability distributions parameterized by $\chi=(\tau,\eta)=\left(\{\tau^t_{iq}\}_{t,i,q},\{\eta^t_{iql}\}_{t,i,q,l}\right)$ of the form
\[
\mathbb{Q}_{\chi}(Z^{1:T})=\prod_{i=1}^n \mathbb{Q}_{\chi}(Z^{1}_i) \prod_{t=2}^T \mathbb{Q}_{\chi}(Z^{t}_i \given Z^{t-1}_i)=\prod_{i=1}^n \left\{ \left[\prod_{q=1}^Q (\tau^1_{iq})^{Z^{1}_{iq}} \right] \prod_{t=1}^{T-1} \prod_{1\leq q,l\leq Q} \left(\frac{\eta^t_{iql}}{\tau^t_{iq}}\right)^{Z_{iq}^{t} Z^{t+1}_{il}} \right\},
\]
i.e. with $\mathbb{Q}_{\chi}$ such that $\esp_{\mathbb{Q}_{\chi}}\left[Z^t_{iq} Z^{t+1}_{il}\right]=\eta^t_{iql}$ and $\esp_{\mathbb{Q}_{\chi}}\left[ Z^t_{iq}\right]=\tau^t_{iq}$. Notice that $\mathbb{Q}_{\chi}(Z^{t+1}_i=l \given Z^{t}_i=q)=\eta_{iql}^t / \tau_{iq}^t=\eta_{iql}^t / \sum_{q'=1}^Q \eta^t_{iqq'}$.
The quantity to optimize in the VEM algorithm is then
\[
\mathcal{J}(\chi,\theta)=\ell(\theta) - KL(\mathbb{Q}_{\chi}, \probat(\cdot | X^{1:T}))=\esp_{\mathbb{Q}_{\chi}} \left[\log \proba_{\theta}(X^{1:T},Z^{1:T})\right] + \mathcal{H}(\mathbb{Q}_{\chi})
\] 
with $KL(\cdot,\cdot)$ denoting the Kullback-Leibler divergence and $\mathcal{H}(\cdot)$ denoting the entropy.
Define
\[
\hat{\chi}(\theta)=(\hat{\tau}(\theta),\hat{\eta}(\theta))=\argmax_{\chi \in [0,1]^{T^2n^2Q^3}}\mathcal{J}(\chi,\theta),
\]
and the variational estimator of $\theta$
\[
\tilde{\theta}=(\tilde{\Gamma},\tilde{\pi})=\argmax_{\theta \in \Theta}\mathcal{J}(\hat{\chi}(\theta),\theta).
\]
Moreover, we denote $\tilde{\chi}=(\tilde{\tau},\tilde{\eta})=\hat{\chi}(\tilde{\theta})=(\hat{\tau}(\tilde\theta),\hat{\eta}(\tilde\theta))$. 
In practice, the VEM algorithm is an iterative algorithm that maximizes the function $\mathcal{J}$ alternatively with
respect to $\chi$ and $\theta$ in order to find $\tilde{\theta}$.

\subsection{Connectivity parameter} \label{sub:consistency_pi_variational}

\begin{theorem} \label{prop:vitesse_cv_M_variationnel}
	For any sequence $\{r_{n,T}\}_{n,T \geq 1}$ increasing to infinity, if $\log(T)=o(n)$, we have for all $\epsilon >0$
	\begin{align*}
	\proba_{\theta^*}\left(\sup_{\theta \in \Theta}\left| \frac{2}{n(n-1)T} \mathcal{J}(\hat{\chi}(\theta),\theta) - \mathbb{M}(\pi) \right| >
	\frac{\epsilon r_{n,T}}{\sqrt{n}}\right) \mathop{\longrightarrow}_{n,T\to +\infty} 0 . 
	\end{align*}
\end{theorem}

We conclude on the consistency of the connection probabilities variational estimators as $n$ and $T$ increase thanks to the following corollary. 
\begin{corollary}\label{cor:vitesse_pitilde_bis}
	For any sequence $\{r_{n,T}\}_{n,T \geq 1}$ increasing to infinity such that $r_{n,T} = o(n^{1/4})$, we have for any $\epsilon>0$
	\begin{align*}
	\frac{1}{2} \proba_{\theta^*}\left(\min_{\sigma \in \mathfrak{S}_Q} \| \tilde{\pi}_{\sigma}-\pi^*\|_{\infty} > \frac{\epsilon r_{n,T}}{n^{1/4}} \right) \xrightarrow[n,T \rightarrow \infty]{} 0.
	\end{align*}
\end{corollary}
We have the equivalent following corollary for a fixed number of time steps.
\begin{corollary}\label{cor:vitesse_cv_pi_fixedT_variationnel}
	If the number of time steps $T$ is fixed, we have for every $\epsilon>0$ and for any sequence $\{r_n\}_{n \geq 1}$ increasing to infinity such that $r_n = o(n^{1/4})$
	\begin{align*}
	\frac{1}{2}\proba_{\theta^*}\left(\min_{\sigma \in \mathfrak{S}_Q}\| \tilde{\pi}_{\sigma}^{1:T}-\pi^{*1:T}  \|_{\infty} > \frac{\epsilon r_n}{n^{1/4}}\right) \xrightarrow[n \rightarrow \infty]{}0.
	\end{align*}
\end{corollary}

\begin{remark} 
	As for Corollaries~\ref{cor:vitesse_cv_pi_bis} to~\ref{cor:Gamma_consistency_fixedT}, the results of Corollaries~\ref{cor:vitesse_pitilde_bis} and~\ref{cor:vitesse_cv_pi_fixedT_variationnel} still hold for any sequences $r_{n,T}$ and $r_n$ increasing to infinity.
\end{remark}

\subsection{Latent transition matrix} \label{sub:consistency_gamma_variational}
We now prove that $\tilde{\Gamma}$ is consistent when the number of nodes and time steps increase. 
\begin{lemme}\label{lem:fixed_point_variational}
	Any critical point 
	$(\breve{\chi},\breve{\theta})$ of the function $\mathcal{J}(\cdot,\cdot)$ 
	is such that $\breve{\Gamma}$ satisfies the fixed-point equation
	\begin{align} \label{eq:fixed_point_variationnel}
	\forall (q,l) \in \Q^2, \quad \breve{\gamma}_{ql} = \frac{\sum_{i=1}^n \sum_{t=1}^{T-1} \breve{\eta}^t_{iql}}{\sum_{i=1}^n \sum_{t=1}^{T-1} \breve{\tau}^{t}_{iq}}.
	\end{align}
\end{lemme}
We assume that $(\tilde{\chi},\tilde{\theta})$ is a critical point of $\mathcal{J}(\cdot,\cdot)$. Then we have the fixed-point equation 
\begin{align}
\forall (q,l) \in \Q^2, \quad \tilde{\gamma}_{ql} = \frac{\sum_{i=1}^n \sum_{t=1}^{T-1} \hat{\eta}^t_{iql}(\tilde{\theta})}{\sum_{i=1}^n \sum_{t=1}^{T-1} \hat{\tau}^{t}_{iq}(\tilde{\theta})}.
\end{align}
The following proposition gives the consistency and a rate of convergence of this estimator, under an assumption on the rate of convergence of $\tilde{\pi}$.

\begin{theorem}\label{prop:consistency_gamma_tilde}
	If $\log(T)=o(n)$, for any $\epsilon>0$ and $\{r_{n,T}\}_{n,T\geq 1}$ any sequence increasing to infinity such that $r_{n,T}=o\left(\sqrt{nT/\log n}\right)$ and for any $\sigma \in \mathfrak{S}_Q$
	\begin{align*}
	&\proba_{\theta^*}\left( \| \tilde{\Gamma}_{\sigma}-\Gamma^* \|_{\infty} > \epsilon r_{n,T} \frac{\sqrt{\log n}}{\sqrt{nT}} \right) \leq 2Q^2(3Q+1)
	\proba_{\theta^*} \left( \| \tilde{\pi}_{\sigma}- \pi^* \|_{\infty} > v_{n,T} \right) + o(1)
	\end{align*}
	with $\{v_{n,T}\}_{n,T \geq 1}$ a sequence decreasing to $0$ such that $v_{n,T}=o(\sqrt{\log(nT)}/n)$.
\end{theorem}
\begin{corollary} \label{cor:Gamma_consistency_variationnel}
	Assume that $\log(T)=o(n)$ and $\min_{\sigma\in \mathfrak{S}_Q}\| \tilde{\pi}_{\sigma}- \pi^* \|_{\infty}=o_{\proba_{\theta^*}}(v_{n,T})$ with $\{v_{n,T}\}_{n,T \geq 1}$ a sequence decreasing to $0$ such that $v_{n,T}=o(\sqrt{\log(nT)}/n)$. Then for any $\epsilon>0$ and $\{r_{n,T}\}_{n,T\geq 1}$ any sequence increasing to infinity such that  $r_{n,T}=o\left(\sqrt{nT/\log n}\right)$, we have the convergence
	\begin{align*}
	\proba_{\theta^*}\left( \min_{\sigma \in \mathfrak{S}_Q} \| \tilde{\Gamma}_{\sigma}-\Gamma^* \|_{\infty} > \epsilon r_{n,T} \frac{\sqrt{\log n}}{\sqrt{nT}} \right) \xrightarrow[n,T \to \infty]{}0.
	\end{align*}
      \end{corollary}
The proof of Corollary~\ref{cor:Gamma_consistency_variationnel} is the same as that of Corollary~\ref{cor:Gamma_consistency}, using Theorem~\ref{prop:consistency_gamma_tilde} instead of Theorem~\ref{prop:Gamma_consistency_bis} and is therefore omitted.

When the number of time steps $T$ is fixed and the connection probabilities can vary over time, we have the following Corollary that is the equivalent of Corollary~\ref{cor:Gamma_consistency_variationnel}. 

\begin{corollary} \label{cor:Gamma_consistency_variationnel_fixedT}
	Let the number of time steps $T \geq 2$ be fixed. Assume that $\min_{\sigma\in \mathfrak{S}_Q}\| \tilde{\pi}^{1:T}_{\sigma}- \pi^{*1:T} \|_{\infty}=o_{\proba_{\theta^*}}(v_{n})$ with $\{v_{n}\}_{n \geq 1}$ a sequence decreasing to $0$ such that $v_{n}=o(\sqrt{\log(n)}/n)$. Then for any $\epsilon>0$ and $\{r_{n}\}_{n\geq 1}$ any sequence increasing to infinity such that  $r_{n}=o\left(\sqrt{n/\log n}\right)$, we have the convergence
	\begin{align*}
	\proba_{\theta^*}\left( \min_{\sigma \in \mathfrak{S}_Q} \| \tilde{\Gamma}_{\sigma}-\Gamma^* \|_{\infty} > \epsilon r_{n} \frac{\sqrt{\log n}}{\sqrt{n}} \right) \xrightarrow[n \to \infty]{}0.
	\end{align*}
\end{corollary}

	The proof of Corollary~\ref{cor:Gamma_consistency_variationnel_fixedT} is the same as that of
        Corollary~\ref{cor:Gamma_consistency_variationnel}, but relying on Theorem~\ref{prop:ratio_alt} instead of
        Theorem~\ref{prop:ratio} and is therefore omitted.

\begin{remark} 
	As for Corollaries~\ref{cor:vitesse_cv_pi_bis} to~\ref{cor:vitesse_cv_pi_fixedT_variationnel}, the results of Corollaries~\ref{cor:Gamma_consistency_variationnel} and~\ref{cor:Gamma_consistency_variationnel_fixedT} still hold for any sequences $r_{n,T}$ and $r_n$ increasing to infinity.
\end{remark}

\section{Proofs of main results}
\label{sec:proofs_main}

\subsection{Proof of Theorem~\ref{prop:vitesse_cv_M}}
The proof follows the lines of the proof of Theorem 3.6 in \cite{celisse2012consistency}. Nonetheless, our result
  is sharper as we establish an upper bound of the rate of convergence (in probability) of the normalised likelihood.
We fix some $\theta \in \Theta$ and introduce the quantities
\begin{align}
\hat{z}^{1:T} & =\argmax_{z^{1:T}\in \Q^{nT}} \log\probat(X^{1:T} \given Z^{1:T}=z^{1:T})  , \label{eq:zhat} \\
 \tilde{Z}^{1:T} &=\argmax_{z^{1:T}\in \Q^{nT}} \espetoile \left[\log\probat(X^{1:T} \given Z^{1:T}=z^{1:T}) \given[\Big]
                   Z^{1:T} \right]  . \label{eq:ztilde}
\end{align}
Note that $\tilde{Z}^{1:T}$ is a random variable that depends on $Z^{1:T}$ and that
\begin{align}
 \hat{z}^{1:T}
     =\argmax_{z^{1:T}\in \Q^{nT}} \sum_{t=1}^T \log\probat(X^{t} \given Z^{t}=z^{t}) =  \left(\argmax_{z\in \Q^n} \log\probat(X^{1} \given
       Z^{1}=z),\ldots,\argmax_{z\in \Q^n} \log\probat(X^{T} \given Z^{T}=z)\right). \label{eq:zthat}
\end{align}
Similarly, for any $t \in \T$,
 we have $\tilde{Z}^{t} =\argmax_{z\in \Q^n} \espetoile \left[\log\probat(X^{t} \given Z^{t}=z) \given Z^t  \right]$.  

We bound the difference between $M_{n,T}(\Gamma,\pi)$ and $\mathbb{M}(\pi)$ by introducing three intermediate terms so that we can write, for any sequence $\{r_{n,T}\}_{n,T\geq 1}$ and any $\epsilon>0$
	\begin{align*}
	&\proba_{\theta^*}\left( \sup_{\theta \in \Theta} \left| M_{n,T}(\Gamma,\pi) - \mathbb{M}(\pi) \right| > \frac{\epsilon r_{n,T}}{\sqrt{n}} \right)\\ \leq& \proba_{\theta^*}\left( \sup_{\theta \in \Theta} \left| \frac{2}{n(n-1)T}\log\probat(X^{1:T}) - \frac{2}{n(n-1)T}\log\probat(X^{1:T} \given {Z}^{1:T} = \hat{z}^{1:T}) \right| > \frac{\epsilon r_{n,T}}{3 \sqrt{n}} \right) \\
	&+ \proba_{\theta^*}\left( \sup_{\theta \in \Theta} \left| \frac{2}{n(n-1)T}\log\probat(X^{1:T} \given {Z}^{1:T} = \hat{z}^{1:T}) -
	\frac{2}{n(n-1)T} \espetoile \left[\log\probat(X^{1:T} \given {Z}^{1:T}=\tilde{Z}^{1:T}) \given[\Big] Z^{1:T} \right] \right| > \frac{\epsilon r_{n,T}}{3 \sqrt{n}}  \right) \\
	&+ \proba_{\theta^*}\left(\sup_{\theta \in \Theta} \left| \frac{2}{n(n-1)T} \espetoile \left[ \log\probat(X^{1:T} \given {Z}^{1:T}=\tilde{Z}^{1:T}) \given[\Big] Z^{1:T} \right] -
	\mathbb{M}(\pi) \right| > \frac{\epsilon r_{n,T}}{3 \sqrt{n}}  \right) . \numberthis  \label{decomposition}
		\end{align*}
In the following, we prove separately the convergence (in $\proba_{\theta^*}$-probability) to zero of the three terms of this sum (while controlling for the rate of these convergences).
Before starting, let us remark that we have 
\begin{align}
  \log\probat(X^{1:T} \given Z^{1:T}=z^{1:T})
      &= \sum_{t=1}^T \sum_{1\le i<j\le n} \Xtij \log \pi_{z^{t}_i z^{t}_j} + (1-\Xtij) \log (1-\pi_{z^{t}_i z^{t}_j})   \label{eq:logvrais_cond} \\
\text{and} \quad \espetoile \left[ \log\probat(X^{1:T} \given Z^{1:T}=z^{1:T}) \given[\Big] Z^{1:T} \right] 
&= \sum_{t=1}^T \sum_{1\le i<j\le n} \pi^{*}_{Z^{t}_i Z^{t}_j} \log \pi_{z^{t}_i z^{t}_j} + (1-\pi^{*}_{Z^{t}_i Z^{t}_j}) \log
(1-\pi_{z^{t}_i z^{t}_j}). \label{eq:logvrais_cond_E*} 
\end{align}
In particular, for every $t\in \T$, we have
\begin{align*}
\hat{z}^{t} & =\argmax_{z=(z_1,\dots,z_n)\in \Q^n} \sum_{1\le i<j\le n} \Xtij \log \pi_{z_i z_j} + (1-\Xtij) \log (1-\pi_{z_i z_j})  ,\\
 \tilde{Z}^{t} &=\argmax_{z=(z_1,\dots,z_n)\in \Q^n} \sum_{1\le i<j\le n } \pi^{*}_{Z^{t}_i Z^{t}_j} \log \pi_{z_i z_j} + (1-\pi^{*}_{Z^{t}_i Z^{t}_j}) \log
(1-\pi_{z_i z_j}).
\end{align*}

\paragraph*{First term of the right-hand side of~\eqref{decomposition}.} 
We let
\begin{align}
	T_1\coloneqq & \left| \frac{2}{n(n-1)T}\log\probat(X^{1:T}) - \frac{2}{n(n-1)T}\log\probat(X^{1:T} \given       {Z}^{1:T}=\hat{z}^{1:T}) \right| \nonumber \\
	\leq &\frac{2}{n(n-1)T} \sum_{t=1}^T \left| \log\probat(X^{t} \given X^{1:t-1}) -  \log\probat(X^{t} \given {Z}^{t}=\hat{z}^{t}) \right|. \label{t1} 
\end{align}

\begin{lemme}
  \label{lemmet1}
For every $t  \in \llbracket 1,T \rrbracket $, we have
\begin{equation*}
\left| \log\probat(X^{t} | X^{1:t-1}) -  \log\probat(X^{t} | {Z}^t=\hat{z}^t) \right| \leq \left| \log
 \probat({Z}^t=\hat{z}^t | X^{1:t-1}) \right| .
\end{equation*}
\end{lemme}
Going back to~\eqref{t1} and applying Lemma~\ref{lemmet1}, we get 
	\begin{equation*}
T_1  \leq \frac{2}{n(n-1)T} \sum_{t=1}^T \left|\log\probat({Z}^{t}=\hat{z}^{t} \given X^{1:t-1}) \right| = - \frac{2}{n(n-1)T} \sum_{t=1}^T \log\probat({Z}^{t}=\hat{z}^{t} \given X^{1:t-1}).
    \end{equation*}
Now, using classical dependency rules in directed acyclic graphs~\cite[see for e.g.][]{Lauritzen} combined with Assumption~\ref{itm:hyp_gamma}, we get
    \begin{align*}
T_1	&\leq - \frac{2}{n(n-1)T} \sum_{t=1}^T \log \sum_{z^{t-1}\in \Q^n} \probat({Z}^{t}=\hat{z}^{t} \given Z^{t-1}=z^{t-1}) \probat(Z^{t-1}=z^{t-1} \given X^{1:t-1}) \\
	&\leq - \frac{2}{n(n-1)T} \sum_{t=1}^T \log \sum_{z^{t-1}\in \Q^n} \delta^n \probat(Z^{t-1}=z^{t-1} \given X^{1:t-1}) 
	\leq - \frac{2}{n(n-1)T} \sum_{t=1}^T n \log \delta =  \frac{2}{n-1} \log (1/\delta) .
	\end{align*}
This implies that $\proba_{\theta^*}(\sup_{\theta \in \Theta} T_1>\epsilon r_{n,T} /(3\sqrt{n})) =0$ as soon as $\epsilon r_{n,T}/\sqrt{n} \geq 6 \log(1/\delta)/(n-1)$. Then for any sequence $\{r_{n,T}\}_{n,T \geq 1}$ increasing to infinity, for any $\epsilon>0$, we have that $\proba_{\theta^*}(\sup_{\theta \in \Theta} T_1>\epsilon r_{n,T} /(3\sqrt{n})) \to 0$ as $n$ and $T$ increase.

\paragraph*{Second term of the right-hand side of~\eqref{decomposition}.} 
Let us denote 
\begin{align*}
T_2(Z^{1:T}) \coloneqq& \left| \frac{2}{n(n-1)T}\log\probat(X^{1:T} \given {Z}^{1:T}=\hat{z}^{1:T}) -
\frac{2}{n(n-1)T} \espetoile \left[ \log\probat(X^{1:T} | {Z}^{1:T}=\tilde{Z}^{1:T}) \given[\Big] Z^{1:T} \right]   \right| .
\end{align*}
For the sake of clarity, we study this term on the event $\{Z^{1:T}=z^{*1:T}\}$ where $z^{*1:T}\in \Q^{nT}$ is a fixed configuration. This event induces the
definition of $\tilde Z^{1:T}$ following Equation~\eqref{eq:ztilde}  as
\[
\tilde{Z}^{1:T} =\argmax_{z^{1:T}\in \Q^{nT}} \espetoile \left[\log\probat(X^{1:T} \given Z^{1:T}=z^{1:T}) \given[\Big]
  Z^{1:T}=z^{*1:T} \right], 
\]
or equivalently for every $t\in \T$, 
\[
\tilde{Z}^{t} =\argmax_{z=(z_1,\dots,z_n)\in \Q^n} \sum_{1\le i<j\le n} \pi^{*}_{z^{*t}_i z^{*t}_j} \log \pi_{z_i z_j} + (1-\pi^{*}_{z^{*t}_i z^{*t}_j}) \log (1-\pi_{z_i z_j}).
\]
By definition of $\hat{z}^{1:T}$ and $\tilde{Z}^{1:T}$ respectively, we have the two inequalities
\[
\log\probat(X^{1:T} \given {Z}^{1:T}=\hat{z}^{1:T}) \geq \log\probat(X^{1:T} \given {Z}^{1:T}=\tilde{Z}^{1:T})
\]
and 
\[
\espetoile \left[ \log\probat(X^{1:T} \given {Z}^{1:T}=\tilde{Z}^{1:T}) \given[\Big] Z^{1:T}=z^{*1:T} \right] \geq \espetoile \left[ \log\probat(X^{1:T} \given {Z}^{1:T}=\hat{z}^{1:T}) \given[\Big] Z^{1:T}=z^{*1:T} \right],
\] 
implying the lower and upper bounds 
\begin{multline*}
\log\probat(X^{1:T} \given {Z}^{1:T}=\tilde{Z}^{1:T}) -
\espetoile \left[ \log\probat(X^{1:T} \given {Z}^{1:T}=\tilde{Z}^{1:T}) \given[\Big] Z^{1:T}=z^{*1:T} \right] 	\\ 
\leq \log\probat(X^{1:T} \given {Z}^{1:T}=\hat{z}^{1:T}) -
\espetoile \left[ \log\probat(X^{1:T} \given {Z}^{1:T}=\tilde{Z}^{1:T}) \given[\Big] Z^{1:T}=z^{*1:T} \right] \\ 
\leq	\log\probat(X^{1:T} \given {Z}^{1:T}=\hat{z}^{1:T}) -
\espetoile \left[ \log\probat(X^{1:T} \given {Z}^{1:T}=\hat{z}^{1:T}) \given[\Big] Z^{1:T}=z^{*1:T} \right].
\end{multline*}
Taking the absolute value gives us an upper bound for $T_2(z^{*1:T})$
\[
T_2 (z^{*1:T})\leq \max_{z^{1:T} \in \{\hat{z}^{1:T}, \tilde{Z}^{1:T}\}} \frac{2}{n(n-1)T}\left|  \log\probat(X^{1:T} \given Z^{1:T}=z^{1:T}) - 
\espetoile \left[ \log\probat(X^{1:T} \given Z^{1:T}=z^{1:T}) \given[\Big] Z^{1:T}=z^{*1:T} \right] \right|.
\]
Using Equations~\eqref{eq:logvrais_cond} and~\eqref{eq:logvrais_cond_E*}, we then obtain the following upper bound for $T_2(z^{*1:T})$ 
\[
T_2(z^{*1:T}) \leq \max_{z^{1:T} \in \{\hat{z}^{1:T}, \tilde{Z}^{1:T}\}}  \left|\frac{2}{n(n-1)T} \sum_{t=1}^T
  \sum_{1\le i<j\le n } (\Xtij-\pi^{*}_{z^{*t}_i z^{*t}_j}) \log\left(\frac{\pi_{z^t_i z^t_j}}{1-\pi_{z^t_i z^t_j}}\right) \right|.
\]
%
We use the following concentration result to conclude.

\begin{lemme} \label{lem:app_talagrand}
  Let $\epsilon, \beta >0$ and   $\{x_{n,T}\}_{n,T\geq 1}$  a sequence of positive real numbers.
We let $\probaetoile(\cdot)$ denote the probability conditional on $\{Z^{1:T}=z^{*1:T}\}$ under parameter
$\theta^*$, i.e. $\probaetoile(\cdot)=\proba_{\theta^*}(\cdot \given Z^{1:T}=z^{*1:T})$.
Denoting $\Lambda=2\log[(1-\zeta)/\zeta]>0$ we have for any $\theta \in \Theta$
\begin{align} \label{eq:majprobsup}
& \probaetoile \left(\sup_{z^{1:T}\in \Q^{nT}} \sup_{\pi \in [\zeta, 1-\zeta]^{Q^2}} \frac{2}{n(n-1)T} \left| \sum_{t=1}^T  \sum_{1\leq i<j\leq n} (\Xtij-\pi^{*}_{z^{*t}_i z^{*t}_j}) \log\left(\frac{\pi_{z^t_i z^t_j}}{1-\pi_{z^t_i z^t_j}}\right) \right| > \epsilon \right) \nonumber \\
 \leq & \probaetoile \left[ \frac{(1+\beta) \Lambda}{\sqrt{n(n-1)T/2}}
+ \frac{\Lambda \sqrt{ x_{n,T}/2}}{\sqrt{n(n-1)T/2}} + (1/\beta+1/3)  \frac{(\Lambda/2) x_{n,T}}{n(n-1)T/2} > \epsilon \right] + 2 e^{-x_{n,T}} \nonumber \\
\leq&  \ind_{2\Omega/(n(n-1)T) > \epsilon} + 2 e^{-x_{n,T}}
\end{align}
with $\Omega=(1+\beta) \Lambda\sqrt{n(n-1)T/2}  + \Lambda\sqrt{n(n-1) T x_{n,T}/4} + (1/\beta+1/3) (\Lambda/2) x_{n,T}$.
\end{lemme}
	Let us choose $x_{n,T}= \log(n)$ in the above lemma. For any $\epsilon>0$, for any sequence $\{r_{n,T}\}_{n,T \geq 1}$ increasing to infinity, we have for $n$ and $T$ large enough 
	\[
	\frac{\epsilon r_{n,T}}{3 \sqrt{n}} \geq \frac{2\Omega}{n (n-1)T}.
	\]
	Then for $n$ and $T$ large enough, the first term in the right-hand side of inequality \eqref{eq:majprobsup} is equal to $0$ and we have 
	\begin{align*}
	& \probaetoile\left(\sup_{\theta \in \Theta} T_2(z^{*1:T}) > \frac{\epsilon r_{n,T}}{3 \sqrt{n}} \right) 
	 \leq  \frac{2}{n}\\
	\text{ and }\quad 
	&\proba_{\theta^*} \left(\sup_{\theta \in \Theta} T_2 (Z^{1:T}) > \frac{\epsilon r_{n,T}}{3 \sqrt{n}} \right) \le \sum_{z^{*1:T}}
	\probaetoile\left(\sup_{\theta \in \Theta} T_2(z^{*1:T}) > \frac{\epsilon r_{n,T}}{3 \sqrt{n}} \right)   \proba_{\theta^*} (Z^{1:T}= z^{*1:T}) \leq \frac{2}{n}.
	\end{align*}
	
\paragraph*{Third term of the right-hand side of~\eqref{decomposition}.} 
	Let us denote 
	\begin{align*}
	T_3(Z^{1:T}) \coloneqq& \left| \frac{2}{n(n-1)T} \espetoile \left[ \log\probat(X^{1:T} \given {Z}^{1:T}=\tilde{Z}^{1:T}) \given[\Big] Z^{1:T} \right] - \mathbb{M(\pi)} \right| \\
	=& \left| \frac{2}{n(n-1)T} \sum_{t=1}^T \espetoile \left[ \log\probat(X^t \given {Z}^t=\tilde{Z}^t) \given[\Big] Z^{t} \right] -
	\mathbb{M}(\pi,\bar{A}_{\pi}) \right| .
	\end{align*}
For any fixed configuration $z^t\in \Q^{n}$, analogous to  Equation~\eqref{eq:logvrais_cond_E*}, we write
	\begin{align*}
\espetoile \left[ \log\probat(X^t \given Z^t=z^t) \given[\Big] Z^{t} \right]
	=& \sum_{1\le i<j\le n} \pi^{*}_{Z^{t}_i Z^{t}_j} \log \pi_{z^{t}_i z^{t}_j} + (1-\pi^{*}_{Z^{t}_i Z^{t}_j}) \log (1-\pi_{z^{t}_i z^{t}_j})\\
	=& \frac{1}{2} \sum_{1\le i\neq j\le n } \pi^{*}_{Z^{t}_i Z^{t}_j} \log \pi_{z^{t}_i z^{t}_j} + (1-\pi^{*}_{Z^{t}_i Z^{t}_j}) \log (1-\pi_{z^{t}_i z^{t}_j})\\
	=& \frac{1}{2} \sum_{1\le q,l,q',l'\le Q} \sum_{1\le i\neq j\le n } \left( \pi^{*}_{ql} \log \pi_{q'l'}+ (1-\pi^{*}_{ql}) \log (1-\pi_{q'l'}) \right) \ind_{\{Z^{t}_i=q, Z^{t}_j=l,z^{t}_i=q',z^{t}_j=l'\}}\\
	=& \frac{1}{2} \sum_{1\le q,l,q',l'\le Q} C_{qq'}(Z^t,z^t)C_{ll'}(Z^t,z^t) \left( \pi^{*}_{ql} \log \pi_{q'l'}+
           (1-\pi^{*}_{ql}) \log (1-\pi_{q'l'}) \right) ,
	\end{align*}
where $C_{qq'}(Z^t,z^t)=|\{i \in \n;Z^{t}_i=q, z_i^t=q'\}|$ is the (random variable) number of nodes classified in group $q$  in the current  (random) configuration $Z^{t}$, while they belong to group $q'$ in (deterministic) configuration $z^t$.
Recall that   $N_{q}(z^t)$ is the number of
nodes assigned to class $q$ by the configuration $z^t$ and let us denote $a^t_{qq'}=a_{qq'}(Z^t,z^t) =C_{qq'}(Z^t,z^t)/N_{q}(Z^{t})$ the
(random) proportion of vertices from class $q$ in $Z^{t}$ attributed to class $q'$ by $z^t$.
We write
\begin{align*}
\frac{2}{n(n-1)} \espetoile \left[ \log\probat(X^t \given Z^t=z^t) \given[\Big] Z^{t}\right] &= \sum_{1\le q,l,q',l'\le Q} \frac{N_{q}(Z^{t}) N_{l}(Z^{t})}{n(n-1)} a^t_{qq'} a^t_{ll'} \left( \pi^{*}_{ql} \log \pi_{q'l'}+ (1-\pi^{*}_{ql}) \log (1-\pi_{q'l'}) \right) \\
&\coloneqq \Phi^t(A^t,\pi) , 
\end{align*}
with $A^t=(a^t_{qq'})_{1\leq q,q'\leq Q}$. 

Now extending these notations to the case where $z^t=\tilde Z^t$, we let  $\tilde{A}^t=(\tilde{a}^t_{qq'})_{1\leq q,q'\leq Q}$ where
$\tilde{a}^t_{qq'}=a_{qq'}(Z^t, \tilde{Z}^t)$. We remark that the definition of $\tilde{Z}^t$ implies that $\tilde{A}^t=\argmax_{A^t \in \mathcal{A}^t(Z^{1:T})} \Phi^t(A^t,\pi)$
with $\mathcal{A}^t(Z^{1:T})$ the (random) subset of stochastic matrices defined for every $t\in \T$ by 
\[
\mathcal{A}^t(Z^{1:T}) = \Big\{A=(n_{ql}/N_q(Z^t))_{1 \leq q,l  \leq Q}; n_{ql} \in \llbracket 0,N_q(Z^t) \rrbracket, \sum_{l=1}^Q n_{ql}=N_q(Z^t) \Big\}.
\] 
Let us also
denote $\bar{A}_{\pi}^t=\argmax_{A \in \mathcal{A}^t(Z^{1:T})} \mathbb{M}(\pi,A) $. Then
\begin{align} 
\sup_{\theta\in \Theta} T_3(Z^{1:T}) &\leq \sup_{\pi\in [\zeta, 1-\zeta]^{Q^2}} \frac{1}{T} \sum_{t=1}^T \left| \Phi^t(\tilde{A}^t,\pi) -
\mathbb{M}(\pi,\bar{A}_{\pi}) \right| \nonumber\\
& \leq \sup_{\pi\in [\zeta, 1-\zeta]^{Q^2}} \frac{1}{T} \sum_{t=1}^T \left| \Phi^t(\tilde{A}^t,\pi) -
\mathbb{M}(\pi,\bar{A}_{\pi}^t) \right| + \frac{1}{T} \sum_{t=1}^T \sup_{\pi\in [\zeta, 1-\zeta]^{Q^2}} \left|\mathbb{M}(\pi,\bar{A}_{\pi}^t)  -\mathbb{M}(\pi,\bar{A}_{\pi}) \right|. \label{eq:T3_decomp}
\end{align}
We start by stating a concentration lemma on the random variable $N_q(Z^t)$ for any $q\in \Q$ and any $t\in \T$.  
\begin{lemme}\label{prop:whp}
  For any $\theta\in \Theta$ and any $\eta \in (0, \delta)$, let 
\begin{equation*}
\Omega_\eta(\theta) \coloneqq \left\{z^{1:T} \in \Q^{nT} ; \forall t \in \T, \forall q \in \Q,
  \frac{N_q(z^{t})}{n} \geq \alpha_q-\eta \right\} .
\end{equation*}
Then $\probat \left( Z^{1:T} \in \Omega_\eta (\theta) \right) \geq 1 -  QT\exp(-2 \eta^2 n)$. 
\end{lemme}
Building on the previous concentration lemma, the following one gives the convergence in $\proba_{\theta^*}$-probability of the second term in the right-hand side of~\eqref{eq:T3_decomp}.
\begin{lemme} \label{lem:cv_matriceA}
	For any $\epsilon>0$, any $\eta \in (0,\delta)$ and $\{r_{n,T}\}_{n,T \geq 1}$ any positive sequence, 
	\begin{align} \label{eq:lem:cv_matriceA}
	\proba_{\theta^*} \left(\frac{1}{T} \sum_{t=1}^T \sup_{\pi\in [\zeta,1-\zeta]^{Q^2}} \left|\mathbb{M}(\pi,\bar{A}_{\pi}^t)  -
	\mathbb{M}(\pi,\bar{A}_{\pi}) \right| > \frac{\epsilon  r_{n,T} }{6 \sqrt{n}}\right) &\leq QT \exp\left(-2 \eta^2 n \right) +\ind_{n \leq  6c\sqrt{n} / [\epsilon r_{n,T} (\delta-\eta)]}
	\end{align}
	with $c=6 (1-\delta)^2 (1-\zeta) \log(1/\zeta) Q^4$.
\end{lemme}
Then taking any $\eta \in (0,\delta)$, for any $\epsilon>0$, for any sequence $\{r_{n,T}\}_{n,T \geq 1}$ increasing to infinity, we have the following inequality for $n$ and $T$ large enough
\begin{equation} \label{eq:vn_third_term}
r_{n,T} > \frac{6c\sqrt{n}}{\epsilon (\delta-\eta)n},
\end{equation}
implying that the probability in Lemma~\ref{lem:cv_matriceA} converges to $0$ as $n$ and $T$ increase for any $\epsilon>0$, as long as $\log T=o(n)$.
Now, for the first term in the right-hand side of \eqref{eq:T3_decomp}, note that we have for every $\pi$ and every $t$
\[ 
\left\{
  \begin{array}{ll}
\Phi^t(\tilde{A}^t,\pi) \geq \Phi^t(\bar{A}_{\pi}^t,\pi) &\text{ because } \tilde{A}^t=\argmax_{A \in \mathcal{A}^t} \Phi^t(A,\pi) \\
\mathbb{M}(\pi,\bar{A}_{\pi}^t) \geq \mathbb{M}(\pi,\tilde{A}^t) &\text{ because }  \bar{A}_{\pi}^t=\argmax_{A \in \mathcal{A}^t} \mathbb{M}(\pi,A).
  \end{array}
\right. 
\]
Then, either $\mathbb{M}(\pi,\bar{A}_{\pi}^t) \leq \Phi^t(\tilde{A}^t,\pi) $  and 
\[
0 \leq
\Phi^t(\tilde{A}^t,\pi)-\mathbb{M}(\pi,\bar{A}_{\pi}^t) \leq \Phi^t(\tilde{A}^t,\pi)-\mathbb{M}(\pi,\tilde{A}^t)
\]
or  $\mathbb{M}(\pi,\bar{A}_{\pi}^t) \geq \Phi^t(\tilde{A}^t,\pi)$ and
\[
0 \leq \mathbb{M}(\pi,\bar{A}_{\pi}^t) -
\Phi^t(\tilde{A}^t,\pi) \leq \mathbb{M}(\pi,\bar{A}_{\pi}^t) - \Phi^t(\bar{A}_{\pi}^t,\pi). 
\]
In both cases, we get that  $\left|
  \Phi^t(\tilde{A}^t,\pi) - \mathbb{M}(\pi,\bar{A}_{\pi}^t)\right| \leq \sup_{A \in \mathcal{A}} \left|
  \Phi^t(A,\pi) - \mathbb{M}(\pi,A) \right|$ 
  for every $t$ and $\pi$, thus obtaining the upper bound
 \begin{align*}
 \sup_{\pi\in [\zeta, 1-\zeta]^{Q^2} } \frac{1}{T} \sum_{t=1}^T \left| \Phi^t(\tilde{A}^t,\pi) -
 \mathbb{M}(\pi,\bar{A}_{\pi}^t) \right|
 \leq& \frac{1}{T} \sum_{t=1}^T \sup_{\pi \in [\zeta, 1-\zeta]^{Q^2}} \sup_{A^t \in \mathcal{A}} \left|
 \Phi^t(A^t,\pi) - \mathbb{M}(\pi,A^t) \right|.
 \end{align*}
Letting
 \[
   \Delta(\zeta)= \sup_{\pi \in [\zeta, 1-\zeta]} \sup_{\pi^* \in [\zeta, 1-\zeta]} |\pi^{*} \log \pi+ (1-\pi^{*}) \log
   (1-\pi)| \in (0,+\infty)
 \]
and recalling that $0\leq a_{ql} \leq1$ (for every $q,l \in \Q$) for every $A=(a_{ql})_{1\leq q,l \leq Q} \in \mathcal{A}$, we have 
\begin{align*}
&\sup_{\pi\in [\zeta, 1-\zeta]^{Q^2}} \sup_{A^t \in \mathcal{A}} \left| \Phi^t(A^t,\pi) - \mathbb{M}(\pi,A^t) \right| \\
&\leq \sup_{\pi\in [\zeta, 1-\zeta]^{Q^2}} \sup_{A^t \in \mathcal{A}} \sum_{1\le q,l,q',l'\le Q} \left| \left(\frac{N_{q}(Z^{t}) N_{l}(Z^{t})}{n(n-1)} - \alpha^*_q \alpha^*_l \right) a^t_{qq'} a^t_{ll'} \left( \pi^{*}_{ql} \log \pi_{q'l'}+ (1-\pi^{*}_{ql}) \log (1-\pi_{q'l'}) \right) \right|\\
&\leq \Delta(\zeta) Q^2\sum_{1\le q,l\le Q} \left|  \frac{ N_{q}(Z^{t}) N_{l}(Z^{t})}{n(n-1)} - \alpha^*_q \alpha^*_l \right| .
\end{align*}
Finally, we bound the first term of the right-hand-side of~\eqref{eq:T3_decomp} as follows
\begin{align} \label{eq:maj_third_term_prop}
\sup_{\pi\in [\zeta, 1-\zeta]^{{Q^2}}} \frac{1}{T} \sum_{t=1}^T \left| \Phi^t(\tilde{A}^t,\pi) -
\mathbb{M}(\pi,\bar{A}_{\pi}^t) \right| \leq \Delta(\zeta) Q^2\sum_{1\le q,l\le Q}  \frac{1}{T} \sum_{t=1}^T \left|  \frac{N_{q}(Z^{t}) N_{l}(Z^{t})}{n(n-1)} - \alpha^*_q \alpha^*_l \right|.
\end{align}
Applying Markov's Inequality, we obtain
\begin{align*}
\proba_{\theta^*}\left( \sup_{\pi \in [\zeta,1-\zeta]^{Q^2} } \frac{1}{T} \sum_{t=1}^T \left| \Phi^t(\tilde{A}^t,\pi) -
\mathbb{M}(\pi,\bar{A}_{\pi}^t) \right| > \frac{\epsilon r_{n,T}}{6 \sqrt{n}} \right) \leq& \sum_{q,l}  \proba_{\theta^*}\left(\frac{1}{T} \sum_{t=1}^T \left| \frac{N_{q}(Z^{t}) N_{l}(Z^{t})}{n(n-1)} - \alpha^*_q \alpha^*_l \right| > \frac{\epsilon r_{n,T}}{6\Delta(\zeta) Q^4 \sqrt{n}}\right)\nonumber \\
\leq& \frac{6\Delta(\zeta) Q^4 \sqrt{n}}{\epsilon r_{n,T}} \sum_{q,l} \frac{1}{T} \sum_{t=1}^T \espetoile\left[ \left| \frac{N_{q}(Z^{t}) N_{l}(Z^{t})}{n(n-1)} - \alpha^*_q \alpha^*_l \right| \right] \\
\leq& \frac{6\Delta(\zeta) Q^4 \sqrt{n}}{\epsilon r_{n,T}} \sum_{q,l} \espetoile\left[ \left| \frac{N_{q}(Z^{1}) N_{l}(Z^{1})}{n(n-1)} - \alpha^*_q \alpha^*_l \right| \right].
\end{align*}	
The following lemma gives an upper bound of the expectation appearing in the previous inequality, for any $q,l \in \Q$.
\begin{lemme} \label{lem:thirdterm}
For any $q,l\in \Q$ and any $t \in \T$, we have the following inequality
	\begin{align*} 
	\espetoile \left[ \left| \frac{N_{q}(Z^{t}) N_{l}(Z^{t})}{n(n-1)} - \alpha^*_q
	\alpha^*_l \right| \right] \leq \frac{2 \sqrt{n}}{n-1}.
	\end{align*}	
\end{lemme}
This leads to
\begin{align*}
	\proba_{\theta^*}\left( \sup_{\pi \in [\zeta,1-\zeta]^{Q^2} } \frac{1}{T} \sum_{t=1}^T \left| \Phi^t(\tilde{A}^t,\pi) -
	\mathbb{M}(\pi,\bar{A}_{\pi}^t) \right| > \frac{\epsilon r_{n,T}}{6 \sqrt{n}} \right)
	\leq& \frac{12\Delta(\zeta) Q^6 n}{\epsilon r_{n,T} (n-1)}.
\end{align*} 
Then for any $\epsilon>0$, for any sequence $\{r_{n,T}\}_{n,T \geq 1}$ increasing to infinity, we have the convergence
\[
\proba_{\theta^*}\left( \sup_{\pi \in [\zeta,1-\zeta]^{Q^2} } \frac{1}{T} \sum_{t=1}^T \left| \Phi^t(\tilde{A}^t,\pi) -
\mathbb{M}(\pi,\bar{A}_{\pi}^t) \right| > \epsilon r_{n,T} /(6\sqrt{n}) \right) \xrightarrow[n,T \rightarrow \infty]{} 0.
\]
We proved the convergence to $0$ of the three terms in the right-hand side of~\eqref{decomposition} for any sequence $\{r_{n,T}\}_{n,T \geq 1}$ increasing to infinity and as long as $\log T=o(n)$. This gives the expected result and concludes the proof.
\qed
\subsection{Proof of Corollary~\ref{cor:vitesse_cv_pi_bis}}
To prove this corollary, we establish the following lemma that allows us to obtain a rate of convergence of $\hat{\pi}$ to $\pi^*$ from a rate of convergence of $M_{n,T}$ to $\mathbb{M}$. Note that this lemma is a bit more general than what we need and gives an equivalent result when the number of time steps $T$ is fixed, which will be useful for Corollary~\ref{cor:vitesse_cv_pi_fixedT}.
\begin{lemme}\label{lem:cv_M_to_cv_pi_bis} 
	Let $\{F_{n,T}\}_{n,T \geq 1}$ be any random functions on the set $\Theta$ (resp. $\Theta^T$) 
	and $\mathbb{M}$ (resp. $\mathbb{M}^T$) 
	defined as before. Assume that there exists a sequence $\{v_{n,T}\}_{n,T \geq 1}$ (resp. $\{v_{n}\}_{n\geq 1}$) a sequence decreasing to $0$ such that for every $\epsilon>0$, we have the following convergence as $n,T \rightarrow \infty$ (resp. $n \rightarrow \infty$) 
	\[
	\proba_{\theta^*} \left( \sup_{(\Gamma,\pi)\in \Theta}\left| F_{n,T}(\Gamma,\pi) - \mathbb{M}(\pi) \right| > \epsilon v_{n,T} \right) \xrightarrow[n,T \rightarrow \infty]{} 0
	\]
	\[
	\left( \textrm{resp. } \proba_{\theta^*} \left( \sup_{(\Gamma,\pi)\in \Theta^T}\left| F_{n,T}(\Gamma,\pi^{1:T}) - \mathbb{M}^T(\pi^{1:T}) \right| > \epsilon v_n \right) \xrightarrow[n \rightarrow \infty]{} 0 \right).
	\]
	If for any $n$ and $T$, $\hat{\theta}=(\hat{\Gamma},\hat{\pi})$ (resp. $\hat{\theta}=(\hat{\Gamma},\hat{\pi}^{1:T})$) 
	is defined as the maximizer of $F_{n,T}$ on the set $\Theta$, (resp. $\Theta^T$)
	we have the following convergence
	\[
	\proba_{\theta^*}\left(\min_{\sigma \in \mathfrak{S}_Q} \|\hat{\pi}_{\sigma}-\pi^*\|_{\infty}> \epsilon \sqrt{v_{n,T}} \right) \xrightarrow[n,T \rightarrow \infty]{} 0
	\]
	\[
	\left( \textrm{resp. } \proba_{\theta^*}\left(\min_{\sigma^{1},\ldots,\sigma^{T} \in \mathfrak{S}_Q}\|\hat{\pi}_{\sigma^{1:T}}^{1:T}-\pi^{*1:T} \|_{\infty}> \epsilon \sqrt{v_n} \right) \xrightarrow[n \rightarrow \infty]{} 0 \right)
	\]
	with $\hat{\pi}_{\sigma^{1:T}}^{1:T}=(\hat{\pi}_{\sigma^t}^{t})_{t \in \T}$.
	\end{lemme}
	The result of Corollary~\ref{cor:vitesse_cv_pi_bis} is then a direct consequence of Theorem~\ref{prop:vitesse_cv_M} (choosing the sequence $\{r_{n,T}^2\}_{n,t \geq 1}$) and  Lemma~\ref{lem:cv_M_to_cv_pi_bis} applied with $F_{n,T}=M_{n,T}$.
	\qed
\subsection{Proof of Theorem~\ref{prop:ratio}}
The proof follows the lines of the proof of Theorem 3.8 in \cite{celisse2012consistency}. Nonetheless, our result
is sharper as we will establish an upper bound of the rate of convergence (in probability) of the quantity at stake.
For any $\epsilon>0$, any sequence $\{y_{n,T}\}_{n,T \geq 1}$ and $\eta \in (0, \delta)$, we write
\begin{align}\label{eq:ratio_cond}
  \proba_{\theta^*} (\mathcal{E}(Z^{1:T},\breve{\theta},\epsilon y_{n,T}) ) =& \sum_{z^{*1:T}\in \Q^{nT}} \proba_{\theta^*}
  (\mathcal{E}(z^{*1:T},\breve{\theta},\epsilon y_{n,T})  ; Z^{1:T}=z^{*1:T})
\leq \proba_{\theta^*} (Z^{1:T}\in \Omega^c_\eta(\theta^*)) 
\nonumber \\ &+ \sum_{z^{*1:T}\in \Omega_\eta(\theta^*)}  \proba_{\theta^*}\left( \frac{\proba_{\breve{\theta}}\left(Z^{1:T} \neq z^{*1:T} \given
	X^{1:T}\right)}{\proba_{\breve{\theta}}\left(Z^{1:T} = z^{*1:T} \given X^{1:T}\right)}  > \epsilon y_{n,T} \given[\Big]
Z^{1:T}=z^{*1:T}\right) \proba_{\theta^*} \left(Z^{1:T}=z^{*1:T}\right)
\end{align}
with $\Omega_\eta(\theta^*)$ as defined in Lemma~\ref{prop:whp}.
We will establish that there exist some positive constants $C, C_1, C_2, C_3, C_4$  such that for any fixed configuration
$z^{*1:T}\in \Omega_{\eta}(\theta^*)$, any $\epsilon>0$, any positive sequence $\{y_{n,T}\}_{n,T \geq 1}$ such that $\log(1/y_{n,T})=o(n)$ and $n$ and $T$ large enough, we have
\begin{multline}\label{eq:prop:ratio}
\proba_{\theta^*}\left[ \frac{\proba_{\breve{\theta}}(Z^{1:T}\neq z^{*1:T} \given X^{1:T}
  )}{\proba_{\breve{\theta}}(Z^{1:T}=z^{*1:T} \given X^{1:T} )} > \epsilon y_{n,T} \given[ \Big]   Z^{1:T}=z^{*1:T}\right] \leq
\proba_{\theta^*} \left( \| \breve{\pi}- \pi^* \|_{\infty} > v_{n,T} \given Z^{1:T}=z^{*1:T}\right)\\
+C nT \left\{\exp \Bigg[ - (\delta - \eta )^2 C_1 n + C_2 \log(nT) + C_4 \log(1/(\epsilon y_{n,T})) \Bigg] + \exp\Bigg[ - C_3
    \frac{(\log(nT))^2}{n v_{n,T}^2}+ 3 n \log(nT)    \Bigg] \right\} . 
\end{multline}
Combined with~\eqref{eq:ratio_cond} and applying Lemma~\ref{prop:whp}, this gives the desired result. 
So now we focus on establishing~\eqref{eq:prop:ratio}.

In what follows, we consider a fixed configuration  $z^{*1:T}\in \Omega_\eta(\theta^*)$ and 
introduce the Hamming distance between $z^{*1:T}$ and any other configuration $z^{1:T}$ defined as 
\[
  \|z^{1:T}-z^{*1:T}\|_0 =\sum_{t=1}^T\sum_{i=1}^n \ind_{z^t_i\neq z^{*t}_i} .
\]
We let $\probaetoile(\cdot)$ denote the probability conditional on $\{Z^{1:T}=z^{*1:T}\}$ under parameter
$\theta=\theta^*$, i.e. $\probaetoile(\cdot)=\proba_{\theta^*}(\cdot \given Z^{1:T}=z^{*1:T})$.
In the following, we will often use the fact that the variables $\{X_{ij}^t\}$ are independent under $\probaetoile$ (with mean value $\pi^*_{z_i^{*t}z_j^{*t}}$) so that we can rely on Hoeffding's Inequality. We introduce a sequence $\{v_{n,T}\}_{n,T \geq 1}$ 
decreasing to 0 and $\Omega_{n,T}$ the event defined as
\[
\Omega_{n,T}=\{ \| \breve{\pi}-\pi^* \|_{\infty} \leq v_{n,T} \}.
\]   
We  bound the probability of interest in~\eqref{eq:prop:ratio} by splitting it on the two complementary events $\Omega_{n,T}$ and $\Omega_{n,T}^c$. For any $\epsilon>0$ and any positive sequence $\{y_{n,T}\}_{n,T \geq 1}$
	\begin{align} \label{eq:decomp_omega}
	\probaetoile \left[ \frac{\proba_{\breve{\theta}}(Z^{1:T}\neq z^{*1:T}\given
          X^{1:T})}{\proba_{\breve{\theta}}(Z^{1:T}=z^{*1:T}\given X^{1:T})} > \epsilon y_{n,T} \right] \leq
\probaetoile \left( \Omega_{n,T}^c  \right) +
          \probaetoile \left[ \left\{ \frac{\proba_{\breve{\theta}}(Z^{1:T}\neq z^{*1:T}\given X^{1:T})}{\proba_{\breve{\theta}}(Z^{1:T}=z^{*1:T} \given X^{1:T})} > \epsilon y_{n,T} \right\} \cap \Omega_{n,T}  \right] .
	\end{align}
Thus, the proof of~\eqref{eq:prop:ratio} boils down to establishing the desired upper bound  on the second term appearing in
the right-hand side of~\eqref{eq:decomp_omega}. We have  
 \begin{multline*}
 \probaetoile \left[ \left\{ \frac{\proba_{\breve{\theta}}(Z^{1:T}\neq z^{*1:T}\given
                 X^{1:T})}{\proba_{\breve{\theta}}(Z^{1:T}=z^{*1:T}\given X^{1:T})} > \epsilon y_{n,T} \right\} \cap \Omega_{n,T}       \right] \\
 \leq \sum_{r=1}^{nT}
 \sum_{z^{1:T} ; \| z^{1:T} -z^{*1:T} \|_{0}=r} \probaetoile \left[ \left\{ 
\frac{\proba_{\breve{\theta}}(Z^{1:T}=z^{1:T}\given X^{1:T})}{\proba_{\breve{\theta}}(Z^{1:T}=z^{*1:T}\given X^{1:T})} 
> \epsilon y_{n,T}/(Q^r (nT)^{r+1}) \right\} \cap
   \Omega_{n,T}\right] ,
   \end{multline*}
by using the bound $ (Q-1)^r \binom{nT}{r} \leq Q^r (nT)^{r}$ on the number of terms in the sum over $\{z^{1:T} ; \|
z^{1:T} -z^{*1:T} \|_{0}=r\}$ (for each value of $r$). Then, 
\begin{align}
&\probaetoile \left[ \left\{ \frac{\proba_{\breve{\theta}}(Z^{1:T}\neq z^{*1:T}\given X^{1:T})}{\proba_{\breve{\theta}}(Z^{1:T}=z^{*1:T}\given X^{1:T})} > \epsilon y_{n,T} \right\} \cap \Omega_{n,T} \right] \nonumber \\
\leq& \sum_{r=1}^{nT} \sum_{z^{1:T} ; \| z^{1:T} -z^{*1:T} \|_{0}=r} \probaetoile \left[ \left\{ \log\frac{\proba_{\breve{\theta}}(Z^{1:T}=z^{1:T} \given X^{1:T})}{\proba_{\breve{\theta}}(Z^{1:T}=z^{*1:T}\given X^{1:T})} > \log(\epsilon y_{n,T})- r\log Q - (r+1)\log(nT) \right\} \cap \Omega_{n,T}\right] \nonumber\\
\leq& \sum_{r=1}^{nT} \sum_{z^{1:T} ; \| z^{1:T} -z^{*1:T} \|_{0}=r} \probaetoile \left[ \left\{ \log\frac{\proba_{\breve{\theta}}(Z^{1:T}=z^{1:T}\given     X^{1:T})}{\proba_{\breve{\theta}}(Z^{1:T}=z^{*1:T}\given X^{1:T})} > - \log(1/(\epsilon y_{n,T})) - 3r\log (nT) \right\} \cap \Omega_{n,T}\right],\label{eq:borne_r_bis}
\end{align}
 as long as $nT \geq Q$. 
For any configuration $z^{1:T}$ such that $\|z^{1:T}-z^{*1:T}\|_0=r$, we denote by $r(1),\ldots,r(T)$ the number of
differences between the two configurations at each time step $t \in \T$, i.e. $ r(t)=\| z^{t} -z^{*t} \|_{0}$ such that $r=\sum_{t} r(t)$. 
Moreover, for any parameter $\pi$, we define $D_{n,T}(z^{1:T}, \pi)$ the subset of indexes $(i,j,t)\in \n^2\times\T$ such that $i<j$ for which
the parameter $\pi$  differs between the configuration $z^{*1:T}$ and $z^{1:T}$, namely  
\begin{equation*}
D_{n,T}(z^{1:T}, \pi) \coloneqq \left\{(i,j,t) \in I_{n,T}; \pi_{z^{t}_i z^{t}_j} \neq \pi_{z^{*t}_i z^{*t}_j} \right\},
	\end{equation*}
	with $I_{n,T}=\{(i,j,t) \in \n^2 \times \T ; i<j \}$ the set of indexes over which we sum to compute the conditional log-likelihood.
In what follows, we abbreviate to $D^*$ (resp. $\breve D$), the set $D_{n,T}(z^{1:T}, \pi^*) $
(resp. $D_{n,T}(z^{1:T},\breve \pi)$).      
Next lemma gives a decomposition of the main term at stake in~\eqref{eq:borne_r_bis}. 
\begin{lemme}\label{lem:ratio_3terms}
 We have the decomposition
\begin{equation*}
	\log \frac{\proba_{\breve{\theta}}(Z^{1:T}=z^{1:T}\given X^{1:T})}{\proba_{\breve{\theta}}(Z^{1:T}=z^{*1:T}\given      X^{1:T})}= U_1 + U_2 - U_3, 
	\end{equation*}
where  
   \begin{align}
U_1 \coloneqq& \sum_{(i,j,t)\in D^*} \left( X_{ij}^t \log\frac{\pi^{*}_{z^t_i z^t_j}}{\pi^{*}_{z^{*t}_i z^{*t}_j}} +
  (1-X_{ij}^t) \log\frac{1-\pi^{*}_{z^t_i z^t_j}}{1-\pi^{*}_{z^{*t}_i z^{*t}_j}} \right) + \sum_{i=1}^n
    \log \frac{\breve{\alpha}_{z^1_i}}{\breve{\alpha}_{z^{*1}_i}} + \sum_{t=1}^{T-1} \sum_{i=1}^n \log
  \frac{\breve{\gamma}_{z^{t}_i z^{t+1}_i}}{\breve{\gamma}_{z^{*t}_i z^{*t+1}_i}} \label{eq:U1}\\
U_2 \coloneqq& \sum_{(i,j,t)\in D^*\cup \breve{D}}\log\left[1+\frac{(\breve{\pi}_{z^t_i z^t_j}-\pi^{*}_{z^t_i  z^t_j})(X_{ij}^t-\pi^{*}_{z^t_i z^t_j})}{\pi^{*}_{z^{t}_i z^{t}_j}(1-\pi^{*}_{z^t_i z^t_j})} \right] \label{eq:U2}\\
U_3 \coloneqq& \sum_{(i,j,t)\in D^*\cup \breve{D}} \log\left[1+\frac{(\breve{\pi}_{z^{*t}_i
    z^{*t}_j}-\pi^{*}_{z^{*t}_i z^{*t}_j})(X_{ij}^t-\pi^{*}_{z^{*t}_i z^{*t}_j})}{\pi^{*}_{z^{*t}_i
 z^{*t}_j}(1-\pi^{*}_{z^{*t}_i z^{*t}_j})} \right] . \label{eq:U3}
\end{align}     
\end{lemme}

Combining~\eqref{eq:borne_r_bis} and Lemma~\ref{lem:ratio_3terms}, we obtain
  \begin{multline}\label{eq:decomp_r}
  \probaetoile \left[ \left\{ \sum_{z^{1:T} \neq z^{*1:T}} \frac{\proba_{\breve{\theta}}(Z^{1:T}=z^{1:T}\given X^{1:T})}{\proba_{\breve{\theta}}(Z^{1:T}=z^{*1:T}\given X^{1:T})} > \epsilon y_{n,T} \right\} \cap \Omega_{n,T} \right]\\
\leq \sum_{r=1}^{nT} \sum_{z^{1:T} ; \| z^{1:T} -z^{*1:T} \|_{0}=r}
\probaetoile \left[ \left\{ U_1+U_2-U_3 > - \log(1/(\epsilon y_{n,T})) - 3r\log (nT) \right\} \cap \Omega_{n,T}\right].
  \end{multline}
We then decompose
	\begin{align} \label{eq:decomp_intermediate_bis}
	&\probaetoile \left[ \left\{ U_1+U_2-U_3 > - \log(1/(\epsilon y_{n,T})) - 3r\log (nT) \right\} \cap \Omega_{n,T}\right] 
	\nonumber \\ \leq& \probaetoile \left[ \left\{ U_1+U_2-U_3 > - \log(1/(\epsilon y_{n,T})) - 3r\log (nT) \right\} \cap \Omega_{n,T} \cap \left\{|U_3|\leq r \log(nT)\right\} \right] \nonumber \\ &+ \probaetoile \left[ \Omega_{n,T} \cap \left\{|U_3| >r \log(nT)\right\} \right]\nonumber \\
	\leq& \probaetoile \left[ \left\{ U_1+U_2 > - \log(1/(\epsilon y_{n,T})) - 4r\log (nT) \right\} \cap \Omega_{n,T}  \right] + \probaetoile \left[ \Omega_{n,T} \cap \left\{|U_3| >r \log(nT)\right\} \right] \nonumber\\
	\leq &\probaetoile \left[U_1 > - \log(1/(\epsilon y_{n,T})) - 5r\log (nT) \right] + \probaetoile \left[ \Omega_{n,T} \cap \left\{|U_2| >r \log(nT)\right\} \right] \nonumber \\
	&+ \probaetoile \left[ \Omega_{n,T} \cap \left\{|U_3| >r \log(nT)\right\} \right].
	\end{align}
We handle these three terms separately in the following.
From now on, we consider a configuration $z^{1:T}$ such that $\| z^{1:T} -z^{*1:T} \|_{0}=r=\sum_{t}r(t)$. 

\paragraph*{First term in the right-hand side of~\eqref{eq:decomp_intermediate_bis}.}
Recall that $U_1$ is given by~\eqref{eq:U1}. 
We can further decompose this term
\begin{align*}
U_1 =& \sum_{(i,j,t)\in D^*} \left( (X_{ij}^t- \pi^{*}_{z^{*t}_i z^{*t}_j}) \log\frac{\pi^{*}_{z^t_i z^t_j}}{\pi^{*}_{z^{*t}_i z^{*t}_j}} \frac{1-\pi^{*}_{z^{*t}_i z^{*t}_j}}{1-\pi^{*}_{z^t_i z^t_j}}\right) \\
&+ \sum_{(i,j,t)\in D^*} \left( \pi^{*}_{z^{*t}_i z^{*t}_j} \log\frac{\pi^{*}_{z^t_i z^t_j}}{\pi^{*}_{z^{*t}_i z^{*t}_j}} + (1-\pi^{*}_{z^{*t}_i z^{*t}_j}) \log\frac{1-\pi^{*}_{z^t_i z^t_j}}{1-\pi^{*}_{z^{*t}_i z^{*t}_j}} \right)+ \sum_{i=1}^n \log \frac{\breve{\alpha}_{z^1_i}}{\breve{\alpha}_{z^{*1}_i}}
+ \sum_{t=1}^{T-1} \sum_{i=1}^n \log \frac{\breve{\gamma}_{z^{t}_i z^{t+1}_i}}{\breve{\gamma}_{z^{*t}_i z^{*t+1}_i}}.
\end{align*}
For $n$ and $T$ large enough such that $\breve{\Gamma} \in [\delta,1-\delta]^{Q^2}$ (implying for the corresponding stationary distribution $\breve{\alpha} \in [\delta,1-\delta]^{Q}$), we have 
\begin{align*}
\sum_{i=1}^n \log \frac{\breve{\alpha}_{z^1_i}}{\breve{\alpha}_{z^{*1}_i}} + \sum_{t=1}^{T-1} \sum_{i=1}^n \log
  \frac{\breve{\gamma}_{z^{t}_i z^{t+1}_i}}{\breve{\gamma}_{z^{*t}_i z^{*t+1}_i}}
& = \sum_{i=1}^n \ind_{\{z^1_i\neq z^{*1}_i\}}\log \frac{\breve{\alpha}_{z^1_i}}{\breve{\alpha}_{z^{*1}_i}} + \sum_{t=1}^{T-1} \sum_{i=1}^n \ind_{\{(z^{t}_i,  z^{t+1}_i)\neq (z^{*t}_i ,z^{*t+1}_i)\}}\log
 \frac{\breve{\gamma}_{z^{t}_i z^{t+1}_i}}{\breve{\gamma}_{z^{*t}_i z^{*t+1}_i}}\\
  &\leq r(1) \log \frac{1-\delta}{\delta} + \sum_{t=1}^{T-1} [r(t)+r(t+1)] \log \frac{1-\delta}{\delta} \leq 2r \log \frac{1-\delta}{\delta}.
\end{align*}
To handle the term $U_1$, we need to lower bound the cardinality of the set
$D^*$. This is the purpose of Lemma~\ref{minor} which is a generalization of
Proposition B.4 in~\cite{celisse2012consistency}.
This can be done for all the configurations $z^{1:T}$ and all the configurations $z^{*1:T}$ that belong to some $\Omega_\eta(\theta)$. 

\begin{lemme}
  \label{minor}
   For any $\eta \in (0, \delta)$, any parameter $\theta\in \Theta$, any configuration
   $z^{1:T}$ and any $z^{*1:T}\in \Omega_\eta(\theta)$  such that $\|z^{1:T}-z^{*1:T}\|_0=r$, we have 
\[
\left| D_{n,T}(z^{1:T},\pi)\right| \geq \frac{(\delta-\eta)^2}{4}n r.
\]
\end{lemme}
Combining Lemma~\ref{minor} with the previous bound, we get that 
\begin{align}\label{eq:maj_apha_gamma}
(|D^*|)^{-1} \left( \sum_{i=1}^n \log \frac{\breve{\alpha}_{z^1_i}}{\breve{\alpha}_{z^{*1}_i}} + \sum_{t=1}^{T-1}    \sum_{i=1}^n \log \frac{\breve{\gamma}_{z^{t}_i z^{t+1}_i}}{\breve{\gamma}_{z^{*t}_i z^{*t+1}_i}} \right) \leq    \frac{2r}{|D^*|} \log \frac{1-\delta}{\delta} \leq \frac{8}{n (\delta-\eta)^2} \log \frac{1-\delta}{\delta}          \xrightarrow[n \rightarrow + \infty]{} 0.
\end{align}
We also have
\[
(|D^*|)^{-1}\sum_{(i,j,t)\in D^*} \left( \pi^{*}_{z^{*t}_i z^{*t}_j} \log\frac{\pi^{*}_{z^t_i z^t_j}}{\pi^{*}_{z^{*t}_i z^{*t}_j}} + (1-\pi^{*}_{z^{*t}_i z^{*t}_j}) \log\frac{1-\pi^{*}_{z^t_i z^t_j}}{1-\pi^{*}_{z^{*t}_i z^{*t}_j}} \right) \leq \max_{q,l, q',l'; \pi^*_{ql} \neq \pi^*_{q'l'}} -k(\pi^*_{ql}, \pi^*_{q'l'}) 
\]
with $k(x, y)= x \log(x/y) + (1-x) \log[(1-x)/(1-y)]$ for $(x,y) \in (0,1)^2$. 
The function $k$ is positive  for every
$(x,y)$ such that $x\neq y$, hence, introducing the notation $K^*=\min_{q,l, q',l'; \pi^*_{ql} \neq \pi^*_{q'l'}} k(\pi^*_{ql}, \pi^*_{q'l'})/2$,
\[
\max_{q,l, q',l'; \pi^*_{ql} \neq \pi^*_{q'l'}} -k(\pi^*_{ql}, \pi^*_{q'l'}) \coloneqq -2K^* <0.
\]
So, by \eqref{eq:maj_apha_gamma}, we have for $n$ large enough 
\[
(|D^*|)^{-1} \left\{\sum_{(i,j,t)\in D^*} \left( \pi^{*}_{z^{*t}_i z^{*t}_j} \log\frac{\pi^{*}_{z^t_i
 z^t_j}}{\pi^{*}_{z^{*t}_i z^{*t}_j}} + (1-\pi^{*}_{z^{*t}_i z^{*t}_j}) \log\frac{1-\pi^{*}_{z^t_i
 z^t_j}}{1-\pi^{*}_{z^{*t}_i z^{*t}_j}} \right) + \sum_{i=1}^n \log \frac{\breve{\alpha}_{z^1_i}}{\breve{\alpha}_{z^{*1}_i}} +
\sum_{t=1}^{T-1} \sum_{i=1}^n \log \frac{\breve{\gamma}_{z^{t}_i  z^{t+1}_i}}{\breve{\gamma}_{z^{*t}_i z^{*t+1}_i}} \right\} \leq -K^*. 
\]
This leads to
\[
\probaetoile(U_1 > u) \leq \probaetoile\left[ \sum_{(i,j,t)\in D^*} \left( (X_{ij}^t- \pi^{*}_{z^{*t}_i z^{*t}_j})
    \log\frac{\pi^{*}_{z^t_i z^t_j}}{\pi^{*}_{z^{*t}_i z^{*t}_j}} \frac{1-\pi^{*}_{z^{*t}_i z^{*t}_j}}{1-\pi^{*}_{z^t_i
        z^t_j}}\right)- |D^*| K^* > u \right] 
\]
for any $u>0$ and large enough $n$.
Moreover, thanks to Hoeffding's Inequality and Assumption~\ref{itm:hyp_pi},
\begin{align*}
\probaetoile(U_1 > u) \leq &\probaetoile \left( \sum_{(i,j,t)\in D^*} \left( (X_{ij}^t- \pi^{*}_{z^{*t}_i z^{*t}_j}) \log\frac{\pi^{*}_{z^t_i z^t_j}}{\pi^{*}_{z^{*t}_i z^{*t}_j}} \frac{1-\pi^{*}_{z^{*t}_i z^{*t}_j}}{1-\pi^{*}_{z^t_i z^t_j}}\right) > u + |D^*| K^* \right)\\
\leq & \exp \left[ - \frac{u^2 + |D^*|^2 K^{*2} + 2u |D^*| K^{*}}{|D^*| C_{\zeta}} \right]\\ 
\leq & \exp \left[ - \frac{ |D^*|^2 K^{*2} + 2u |D^*| K^{*}}{|D^*| C_{\zeta}} \right] = \exp \left[ - \frac{ 2u K^{*}}{ C_{\zeta}} \right] \exp \left[ - \frac{ |D^*| K^{*2}}{C_{\zeta}} \right],
\end{align*}
where $C_{\zeta}$ is a constant depending on $\zeta$. Finally using Lemma~\ref{minor}, we have
\begin{align*}
\probaetoile \left(U_1 > - \log(1/(\epsilon y_{n,T}))- 5r\log (nT) \right) \leq 	& \exp \left[ \left[ \log(1/(\epsilon y_{n,T}))+ 5 r \log(nT)\right] \frac{ 2 K^{*}}{ C_{\zeta}} \right] \exp \left[ - \frac{ |D^*| K^{*2}}{C_{\zeta}} \right]\\
\leq & \exp \left[ \left[ \log(1/(\epsilon y_{n,T}))+ 5 r \log(nT)\right] \frac{ 2 K^{*}}{ C_{\zeta}} \right] \exp \left[ - nr \frac{ {(\delta-\eta)}^2 K^{*2}}{4 C_{\zeta}} \right].
\end{align*}

\paragraph*{Second term in the right-hand side of \eqref{eq:decomp_intermediate_bis}.}
We have
\[
U_2 \coloneqq \sum_{(i,j,t)\in D^*\cup \breve{D}}\log\left[1+\frac{(\breve{\pi}_{z^t_i z^t_j}-\pi^{*}_{z^t_i
      z^t_j})(X_{ij}^t-\pi^{*}_{z^t_i z^t_j})}{\pi^{*}_{z^{t}_i z^{t}_j}(1-\pi^{*}_{z^t_i z^t_j})}  \right] 
\leq 
\sum_{(i,j,t)\in D^*\cup \breve{D}}\frac{(\breve{\pi}_{z^t_i z^t_j}-\pi^{*}_{z^t_i z^t_j})(X_{ij}^t-\pi^{*}_{z^t_i z^t_j})}{\pi^{*}_{z^{t}_i z^{t}_j}(1-\pi^{*}_{z^t_i z^t_j})}.
\]
For any $q,l,q',l' \in \Q$, we introduce the sets
\begin{align*}
  F_{qlq'l'}&=  F_{qlq'l'} (z^{1:T},z^{*1:T})\coloneqq \{(i,j,t)\in I_{n,T} ;
              z_i^t=q,z_j^t=l,z_i^{*t}=q',z_j^{*t}=l'\}\\
 F_{ql}&=  F_{ql} (z^{1:T})\coloneqq \cup_{1\leq q',l' \leq Q} F_{qlq'l'} =\{(i,j,t)\in I_{n,T} ; z_i^t=q,z_j^t=l\}\\
   G_{qlq'l'} &=   G_{qlq'l'} (z^{1:T}, z^{*1:T},\pi^*,\breve \pi)\coloneqq (D^*\cup \breve{D}) \cap F_{qlq'l'}\\
   &=\{(i,j,t)\in I_{n,T};  z_i^t=q,z_j^t=l,z_i^{*t}=q',z_j^{*t}=l' \text{ and } (\pi^*_{z_i^t z_j^t} \neq
 \pi^*_{z_i^{*t} z_j^{*t}} \text{ or } \breve \pi_{z_i^t z_j^t} \neq \breve \pi_{z_i^{*t} z_j^{*t}} )\} \\
 G_{ql} &=   G_{ql} (z^{1:T}, z^{*1:T},\pi^*,\breve \pi)\coloneqq (D^*\cup \breve{D}) \cap  F_{ql} \\
   &=\{(i,j,t)\in I_{n,T} ;  z_i^t=q,z_j^t=l \text{ and } (\pi^*_{z_i^t z_j^t} \neq
 \pi^*_{z_i^{*t} z_j^{*t}} \text{ or } \breve \pi_{z_i^t z_j^t} \neq \breve \pi_{z_i^{*t} z_j^{*t}} )\} .
\end{align*}
Then we bound
\begin{align} \label{eq:prop3_secondterm}
|U_2| \leq& \sum_{1 \leq q,l \leq Q} \frac{\left| \breve{\pi}_{ql}-\pi^{*}_{ql}\right| }{\pi^{*}_{ql}(1-\pi^{*}_{ql})} \left| \sum_{(i,j,t)\in D^*\cup \breve{D}} (X_{ij}^t-\pi^{*}_{ql}) \ind_{z_i^t=q,z_j^t=l} \right|
  \leq \sum_{1 \leq q,l \leq Q} \frac{\left| \breve{\pi}_{ql}-\pi^{*}_{ql}\right| }{\pi^{*}_{ql}(1-\pi^{*}_{ql})} \left|  \sum_{(i,j,t)\in    G_{ql}} (X_{ij}^t-\pi^{*}_{ql}) 
  \right| \nonumber\\
\leq& \sum_{1 \leq q,l \leq Q} \frac{\left| \breve{\pi}_{ql}-\pi^{*}_{ql}\right| }{\pi^{*}_{ql}(1-\pi^{*}_{ql})} \left|
      \sum_{(i,j,t)\in    G_{ql} } (X_{ij}^t-\pi^{*}_{z_i^{*t}z_j^{*t}} ) \right| + \sum_{1 \leq q,l \leq Q} \frac{\left|
      \breve{\pi}_{ql}-\pi^{*}_{ql}\right| }{\pi^{*}_{ql}(1-\pi^{*}_{ql})} \left|   \sum_{(i,j,t)\in    G_{ql} }
      (\pi^{*}_{z_i^{*t}z_j^{*t}}-\pi^{*}_{ql})  \right|  \nonumber\\
\leq& \sum_{1 \leq q,l \leq Q} \frac{\left| \breve{\pi}_{ql}-\pi^{*}_{ql}\right| }{\pi^{*}_{ql}(1-\pi^{*}_{ql})} \left|
      \sum_{(i,j,t)\in    G_{ql} } (X_{ij}^t-\pi^{*}_{z_i^{*t}z_j^{*t}} ) \right| + \sum_{1 \leq q,l \leq Q} \frac{\left|
      \breve{\pi}_{ql}-\pi^{*}_{ql}\right| }{\pi^{*}_{ql}(1-\pi^{*}_{ql})} \left| \sum_{q',l'}
      (\pi^{*}_{q'l'}-\pi^{*}_{ql}) | G_{qlq'l'}| \right| .
\end{align}
For every $u>0$, we thus have
\begin{align}
  \label{eq:decomp_U2}
\probaetoile(\Omega_{n,T}\cap \{|U_2|>u\}) \leq& \probaetoile\left(\left\{  \sum_{1 \leq q,l \leq Q} 
   \frac{\left| \breve{\pi}_{ql}-\pi^{*}_{ql} \right| }{\pi^{*}_{ql}(1-\pi^{*}_{ql})}
   \left|  \sum_{(i,j,t)\in    G_{ql} }  (X_{ij}^t-\pi^{*}_{z_i^{*t}z_j^{*t}}) \right| > u/2 \right\} \cap \Omega_{n,T} \right) \nonumber \\
&+ \probaetoile \left( \left\{ \sum_{1 \leq q,l \leq Q}  \frac{\left|\breve{\pi}_{ql}-\pi^{*}_{ql}\right| }{\pi^{*}_{ql}(1-\pi^{*}_{ql})} \left| \sum_{1\leq q',l' \leq Q}  (\pi^{*}_{q'l'}-\pi^{*}_{ql}) |   G_{qlq'l'}| \right|  > u/2 \right\} \cap \Omega_{n,T} \right) .
\end{align}
We start by dealing with the first term of the right-hand side of~\eqref{eq:decomp_U2}. Notice that on the event
$\Omega_{n,T}$, we have $\left| (\breve{\pi}_{ql}-\pi^{*}_{ql})/(\pi^{*}_{ql}(1-\pi^{*}_{ql}))\right|\leq
v_{n,T}/\zeta^2$ for every $q,l \in \Q$.
The next lemma establishes that any set $D_{n,T}(z^{1:T}, \pi) $ is included in a larger set, whose cardinality is
bounded. In particular, the random set $\breve D$ is included in a larger deterministic subset. 

 \begin{lemme}
 \label{major}
Let $z^{1:T}$ and $z^{*1:T}$ denote two configurations such that
$   \|z^{1:T}-z^{*1:T}\|_0=r$. Then  for any parameter $\pi=(\pi_{ql})_{1\leq q,l \leq Q}$, we have
 \begin{equation*}
     D_{n,T}(z^{1:T}, \pi)  \subset D_{n,T}(z^{1:T}) \coloneqq 
  \left\{ (i,j,t) \in \n^2\times \T; (z_i^t ,z_j^t) \neq (z_i^{*t}, z_j^{*t}) \right\} 
   \text{ and }  \left| D_{n,T}(z^{1:T})    \right| \leq 2nr.
 \end{equation*}
\end{lemme}
As the set $G_{ql}$ is random (because $\breve D$ is random), we write 
\begin{align*}
& \probaetoile\left(\left\{  \sum_{1 \leq q,l \leq Q} \frac{\left|\breve{\pi}_{ql}-\pi^{*}_{ql}\right| }{\pi^{*}_{ql}(1-\pi^{*}_{ql})} \left|  \sum_{(i,j,t)\in G_{ql}}  (X_{ij}^t-\pi^{*}_{z_i^{*t}z_j^{*t}}) \right| > u/2 \right\} \cap \Omega_{n,T} \right) \\
\leq&  \probaetoile\left(\sum_{1 \leq q,l \leq Q} \left|  \sum_{(i,j,t)\in G_{ql}}  (X_{ij}^t-\pi^{*}_{z_i^{*t}z_j^{*t}}) \right| > \frac{u\zeta^{2}}{2v_{n,T}} \right)
\leq \sum_{D\subset D_{n,T}(z^{1:T}) } \probaetoile\left(  \sum_{1 \leq q,l \leq Q} \left|  \sum_{(i,j,t)\in
  F_{ql}\cap D}  (X_{ij}^t-\pi^{*}_{z_i^{*t}z_j^{*t}}) \right| > \frac{u\zeta^{2}}{2v_{n,T}} \right)  ,
\end{align*}
where now $D $ is a deterministic set. 
By a union bound and Hoeffding's inequality, we have for any $D \subset D_{n,T}(z^{1:T})$
\begin{align*}
\probaetoile\left(\sum_{1 \leq q,l \leq Q} \left|  \sum_{(i,j,t)\in F_{ql}\cap D}  (X_{ij}^t-\pi^{*}_{z_i^{*t}z_j^{*t}}) \right| > \frac{u\zeta^{2}}{2v_{n,T}} \right)
\leq&Q^2 \max_{1 \leq q,l \leq Q} \probaetoile\left( \left|  \sum_{(i,j,t)\in F_{ql}\cap D}  (X_{ij}^t-\pi^{*}_{z_i^{*t}z_j^{*t}})
      \right| > \frac{u\zeta^{2}}{2v_{n,T}} \right)\\
\leq& 2 Q^2 \exp\left(-\frac{2 u^2 \zeta^4}{4 v_{n,T}^2 Q^4}\frac{1}{{| D|}}\right).
\end{align*}
This leads to
\begin{align*}
\probaetoile\left(\left\{  \sum_{1 \leq q,l \leq Q} \frac{\left|\breve{\pi}_{ql}-\pi^{*}_{ql}\right| }{\pi^{*}_{ql}(1-\pi^{*}_{ql})} \left|  \sum_{(i,j,t)\in G_{ql}}  (X_{ij}^t-\pi^{*}_{z_i^{*t}z_j^{*t}}) \right| > u/2 \right\} \cap \Omega_{n,T} \right) 
& \leq \sum_{D\subset D_{n,T}(z^{1:T}) } 2 Q^2 \exp\left(-\frac{2 u^2 \zeta^4}{4 v_{n,T}^2 Q^4}\frac{1}{{|D|}}\right) \\
& \leq \sum_{k=1}^{2nr}\sum_{D\subset D_{n,T}(z^{1:T}) ; |D|=k} 2 Q^2 \exp\left(-\frac{2 u^2 \zeta^4}{4 v_{n,T}^2 Q^4}\frac{1}{{k}}\right) \\
&\leq 2 Q^2 \sum_{k=1}^{2nr} (2nr)^k \exp\left(-\frac{2 u^2 \zeta^4}{4 v_{n,T}^2 Q^4}\frac{1}{2nr}\right)\\
&\leq 2 Q^2 \exp\left(-\frac{ u^2 \zeta^4}{4 v_{n,T}^2 Q^4nr}\right)  (2nr)^{2nr+1}.
\end{align*}
For the second term of \eqref{eq:decomp_U2}, we get from a union bound and from Lemma~\ref{major} (that gives an upper
bound for $|D^* \cup \breve{D}|$) that 
\begin{align*}
&\probaetoile \left( \left\{ \sum_{1 \leq q,l \leq Q} \left| \frac{(\breve{\pi}_{ql}-\pi^{*}_{ql}) }{\pi^{*}_{ql}(1-\pi^{*}_{ql})} \right| \left| \sum_{1\leq q',l' \leq Q}  (\pi^{*}_{q'l'}-\pi^{*}_{ql}) |G_{qlq'l'}| \right|  > u/2 \right\} \cap \Omega_{n,T} \right)\\
\leq& \probaetoile \left( \sum_{1 \leq q,l \leq Q}  \left| \sum_{1\leq q',l' \leq Q}  (\pi^{*}_{q'l'}-\pi^{*}_{ql}) |G_{qlq'l'}| \right|  > \frac{u\zeta^{2}}{2v_{n,T}} \right)\\
\leq& Q^2 \max_{1 \leq q,l \leq Q} \probaetoile \left( \left| \sum_{1 \leq q',l' \leq Q}  (\pi^{*}_{q'l'}-\pi^{*}_{ql})|G_{qlq'l'}| \right|  > \frac{u\zeta^{2}}{2v_{n,T}Q^2} \right)
\leq Q^2 \probaetoile \left( 2nr  > \frac{u\zeta^{2}}{2v_{n,T}Q^2} \right),
\end{align*}
because $|\pi^{*}_{q'l'}-\pi^{*}_{ql}| \leq 1$, implying that 
\[
 \left| \sum_{q',l'}
  (\pi^{*}_{q'l'}-\pi^{*}_{ql})  |G_{qlq'l'}|\right| \leq \sum_{q',l'} |G_{qlq'l'}| = |G_{ql}|= |F_{ql} \cap (D^* \cup \breve
D)|\leq |D_{n,T}(z^{1:T})|\leq  2nr.
\]
Finally, we have the following upper bound for the second term of \eqref{eq:decomp_intermediate_bis}  
\begin{equation*}
\probaetoile \left( \Omega_{n,T} \cap \left\{|U_2| >r \log(nT)\right\} \right) \leq 
2 Q^2 \exp\left(-\frac{ r\zeta^4(\log(nT))^2 }{4 Q^4 v_{n,T}^2 n}\right)  (2nr)^{2nr+1}
+ Q^2 \probaetoile \left(v_{n,T}   > \frac{\zeta^{2} \log(nT)}{4Q^2n} \right). 
\end{equation*}

\paragraph*{Third term in the right-hand side of \eqref{eq:decomp_intermediate_bis}.}
We want to bound (in probability) the last term $U_3$. Distinguishing between the cases where $X_{ij}^t=0$ and $X_{ij}^t=1$, we have 
\begin{align*}
U_3 \coloneqq&\sum_{(i,j,t)\in D^*\cup \breve{D}} \log\left[1+\frac{(\breve{\pi}_{z^{*t}_i z^{*t}_j}-\pi^{*}_{z^{*t}_i z^{*t}_j})(X_{ij}^t-\pi^{*}_{z^{*t}_i z^{*t}_j})}{\pi^{*}_{z^{*t}_i z^{*t}_j}(1-\pi^{*}_{z^{*t}_i z^{*t}_j})} \right]\\
=&\sum_{(i,j,t)\in D^*\cup \breve{D}} \left( (1-X_{ij}^t) \log\left[1-\frac{(\breve{\pi}_{z^{*t}_i z^{*t}_j}-\pi^{*}_{z^{*t}_i z^{*t}_j})}{(1-\pi^{*}_{z^{*t}_i z^{*t}_j})} \right] + X_{ij}^t \log\left[1+\frac{(\breve{\pi}_{z^{*t}_i z^{*t}_j}-\pi^{*}_{z^{*t}_i z^{*t}_j})}{\pi^{*}_{z^{*t}_i z^{*t}_j}} \right] \right)\\
=& \sum_{1\leq q,l \leq Q}\sum_{(i,j,t)\in D^*\cup \breve{D}} \left( (1-X_{ij}^t)
   \log\left[1-\frac{(\breve{\pi}_{ql}-\pi^{*}_{ql})}{(1-\pi^{*}_{ql})} \right] + X_{ij}^t
   \log\left[1+\frac{(\breve{\pi}_{ql}-\pi^{*}_{ql})}{\pi^{*}_{ql}} \right] \right) \ind_{z_i^{*t}=q, z_j^{*t}=l} .
\end{align*}
For any $(q,l)\in \Q^2$, we further introduce the sets
\begin{align*}
F^*_{ql}&=\cup_{1\leq q',l' \leq Q} F_{q'l'ql} =\{(i,j,t)\in I_{n,T} ; z_i^{*t}=q,z_j^{*t}=l\} \\
G^*_{ql}&=\cup_{1 \leq q',l' \leq Q} G_{q'l'ql} = (D^*\cup \breve D) \cap F^*_{ql}= \{(i,j,t)\in D^*\cup \breve{D}  ; z_i^{*t}=q,z_j^{*t}=l\} .
\end{align*}
Centering the $X_{ij}^t$ (under the distribution $\probaetoile$), we get
\begin{align*}
 U_3 =& \sum_{1 \leq q,l \leq Q}\sum_{(i,j,t)\in D^*\cup \breve{D}} \left( (\pi^{*}_{ql} -X_{ij}^t) \log\left[1-\frac{(\breve{\pi}_{ql}-\pi^{*}_{ql})}{(1-\pi^{*}_{ql})} \right] + (X_{ij}^t-\pi^{*}_{ql}) \log\left[1+\frac{(\breve{\pi}_{ql}-\pi^{*}_{ql})}{\pi^{*}_{ql}} \right] \right) \ind_{z_i^{*t}=q, z_j^{*t}=l}\\
&+\sum_{1 \leq q,l \leq Q}\sum_{(i,j,t)\in D^*\cup \breve{D}} \left( (1-\pi^{*}_{ql}) \log\left[1-\frac{(\breve{\pi}_{ql}-\pi^{*}_{ql})}{(1-\pi^{*}_{ql})} \right] + \pi^{*}_{ql} \log\left[1+\frac{(\breve{\pi}_{ql}-\pi^{*}_{ql})}{\pi^{*}_{ql}} \right] \right) \ind_{z_i^{*t}=q, z_j^{*t}=l}\\
=& \sum_{1 \leq q,l \leq Q} \left( \log\left[1+\frac{(\breve{\pi}_{ql}-\pi^{*}_{ql})}{\pi^{*}_{ql}} \right]- \log\left[1-\frac{(\breve{\pi}_{ql}-\pi^{*}_{ql})}{(1-\pi^{*}_{ql})} \right]\right) \sum_{(i,j,t)\in G^*_{ql}} (X_{ij}^t-\pi^{*}_{ql}) \\
&+\sum_{1 \leq q,l \leq Q} |G^*_{ql}| \left( (1-\pi^{*}_{ql}) \log\left[1-\frac{(\breve{\pi}_{ql}-\pi^{*}_{ql})}{(1-\pi^{*}_{ql})} \right]+ \pi^{*}_{ql} \log\left[1+\frac{(\breve{\pi}_{ql}-\pi^{*}_{ql})}{\pi^{*}_{ql}} \right] \right) .
\end{align*}
Then, on the event $\Omega_{n,T}$
and for $n$ and $T$ large enough such that
$|(\breve{\pi}_{ql}-\pi^{*}_{ql})/(1-\pi^{*}_{ql})|\leq 1/2$ and $|(\breve{\pi}_{ql}-\pi^{*}_{ql})/\pi^{*}_{ql}|\leq 1/2$ for every $q$ and $l$, using the fact that $|\log(1+x)| \leq 2|x|$ for $x \in [-1/2,1/2]$, we have
\begin{align*}
|U_3| \leq& 4 \frac{v_{n,T}}{\zeta} \sum_{1 \leq q,l \leq Q} \left| \sum_{(i,j,t)\in G^*_{ql}} (X_{ij}^t-\pi^{*}_{ql}) \right| 
+ 4 \frac{v_{n,T}}{\zeta} \sum_{1 \leq q,l \leq Q} |G^*_{ql}|.
\end{align*}
Then, for every $u>0$,
\begin{align} \label{eq:decomp_U3}
\probaetoile \left( \Omega_{n,T} \cap \left\{ \left| U_3 \right| > u \right\} \right)
\leq& \probaetoile \left(  \sum_{1 \leq q,l \leq Q} \left| \sum_{(i,j,t)\in G^*_{ql}} (X_{ij}^t-\pi^{*}_{ql})
      \right| > \frac{u\zeta}{8 v_{n,T}} \right)
+ \probaetoile \left( v_{n,T} \sum_{1 \leq q,l \leq Q} |G^*_{ql} |> \frac{u \zeta}{8} \right).
\end{align}
For the first term of~\eqref{eq:decomp_U3}, using Hoeffding's inequality as before,
\begin{align*}
\probaetoile \left(  \sum_{1 \leq q,l \leq Q} \left| \sum_{(i,j,t)\in G^*_{ql}} (X_{ij}^t-\pi^{*}_{ql})
\right| > \frac{u\zeta}{8 v_{n,T}} \right)
\leq& \sum_{k =1}^{2nr} \sum_{D \subset D_{n,T} (z^{1:T}); |D|=k} \probaetoile \left( \sum_{1 \leq q,l \leq Q}
      \left| \sum_{(i,j,t)\in D\cap F^*_{ql}} (X_{ij}^t-\pi^{*}_{ql}) \right| > \frac{u\zeta}{8 v_{n,T}} \right)\\
 \leq& 2 Q^2 (2nr)^{2nr+1} \exp \left(-\frac{u^2 \zeta^2}{8^2 Q^4 v_{n,T}^2 n r} \right).
\end{align*}
For the second term of~\eqref{eq:decomp_U3}, we use 
\[
\probaetoile \left( v_{n,T} \sum_{1 \leq q,l \leq Q} |G^*_{ql} | > \frac{u \zeta}{8} \right) \leq \probaetoile \left(  v_{n,T}
  > \frac {u \zeta}{16nr} \right).
\]
Finally, we have the following upper bound for the third term of~\eqref{eq:decomp_intermediate_bis} 
\[
\probaetoile \left( \Omega_{n,T} \cap \left\{|U_3| >r \log(nT)\right\} \right) 
\leq 2 Q^2 (2nr)^{2nr+1} \exp \left(-\frac{r(\log(nT))^2 \zeta^2}{8^2 Q^4 v_{n,T}^2 n} \right)+  \probaetoile
    \left(  v_{n,T}  > \frac{\log (nT)\zeta}{16n} \right). 
\]

 \paragraph*{Combining the 3 bounds on the right-hand-side of~\eqref{eq:decomp_intermediate_bis}.}
\begin{align*}
&\probaetoile \left( \left\{ U_1+U_2-U_3 > -\log(1/(\epsilon y_{n,T})) - 3r\log (nT) \right\} \cap \Omega_{n,T}\right)\\
\leq&
\exp \left[\left[ \log(1/(\epsilon y_{n,T}))+ 5 r \log(nT)\right] \frac{ 2 K^{*}}{ C_{\zeta}} \right] \exp \left[ - nr \frac{ (\delta-\eta)^2 K^{*2}}{4 C_{\zeta}} \right] +
 2 Q^2 (2nr)^{2nr+1}\exp\left[-\frac{r  \zeta^4(\log(nT))^2}{4 Q^4 v_{n,T}^2 n}\right] \\ 
 &+ Q^2 \probaetoile \left( v_{n,T} > \frac{\zeta^{2}\log(nT)}{4  Q^2n} \right) +
 2 Q^2  (2nr)^{2nr+1}\exp \left[-\frac{r(\log(nT))^2 \zeta^2}{8^2 Q^4 v_{n,T}^2 n} \right]
 + \probaetoile \left( v_{n,T}  > \frac{\log(nT) \zeta}{16 n} \right).
\end{align*}
Now we choose the sequence $v_{n,T}$ such that $v_{n,T}=o(\sqrt{\log(nT)}/n)$ which is sufficient to imply  that
the quantities $\probaetoile \left( v_{n,T} > \zeta^{2}\log(nT)/(4 Q^2 n) \right)$ and $\probaetoile \left( v_{n,T}  >
  \log(nT) \zeta/(16n) \right)$ vanish as $n$ and $T$ increase. For large enough values of $n$ and $T$ and with $C_1$,
$C_2, C_3, C_4$ and $\kappa$ positive constants only depending on $Q, \zeta$ and $K^*$, we then have
\begin{align} 
&\probaetoile \left( \left\{ U_1+U_2-U_3 > -\log(1/(\epsilon y_{n,T})) - 3r\log (nT) \right\} \cap \Omega_{n,T}\right) \nonumber\\
\leq&
\exp \left[\left[ \log(1/(\epsilon y_{n,T}))+ 5 r \log(nT)\right] \frac{ 2 K^{*}}{ C_{\zeta}} \right] \exp \left[ - nr \frac{ (\delta-\eta)^2 K^{*2}}{4 C_{\zeta}} \right] +
 2 Q^2 (2nr)^{2nr+1}\exp\left[-\frac{r  \zeta^4(\log(nT))^2}{4 Q^4 v_{n,T}^2 n}\right] \nonumber \\ 
 & +
 2 Q^2  (2nr)^{2nr+1}\exp \left[-\frac{r(\log(nT))^2 \zeta^2}{8^2 Q^4 v_{n,T}^2 n} \right] \nonumber\\
\leq &\exp \Bigg[ - (\delta - \eta )^2 C_1 nr + C_2 \log(nT)r + C_4 \log(1/(\epsilon y_{n,T})) \Bigg] +
\kappa  \exp\Bigg[ 3 nr \log(nT)  - C_3 \frac{(\log(nT))^2r}{n v_{n,T}^2}  \Bigg] .
\label{eq:vanish}
\end{align}
Let us introduce
\begin{align*}
  u_{nT}&=\exp \left[ - (\delta - \eta )^2 C_1 n + C_2 \log(nT) + C_4 \log(1/(\epsilon y_{n,T}))  \right]  \\
w_{nT}&= \exp\left[ - C_3 \frac{(\log(nT))^2}{n v_{n,T}^2}+ 3 n \log(nT) \right].     
\end{align*}
Now we go back to~\eqref{eq:decomp_r}. 
Noticing that the number of configurations $z^{1:T}$ such that $\|z^{1:T}
   -z^{*1:T}\|_0=r$ is equal to $\dbinom{nT}{r} (Q-1)^r$, we have
\begin{align*}
\probaetoile \left( \left\{ \frac{\proba_{\breve{\theta}}(Z^{1:T}\neq z^{*1:T}\given X^{1:T})}{\proba_{\breve{\theta}}(Z^{1:T}=z^{*1:T}\given X^{1:T})} > \epsilon y_{n,T} \right\} \cap \Omega_{n,T}\right)
\leq& \sum_{r=1}^{nT} \dbinom{nT}{r} (Q-1)^r u_{nT}^r + \sum_{r=1}^{nT} \dbinom{nT}{r} (Q-1)^r
\kappa w_{nT}^r\\
\leq& [1+Q u_{nT}]^{nT} - 1 + \kappa \left( [1+Q w_{nT}]^{nT} -1 \right) .\\
\end{align*}
Finally, notice that as long as $\log T=o(n)$ and $\log(1/y_{n,T})=o(n)$ (resp. as long as $v_{n,T}=o(\sqrt{\log(nT)}/n)$), we have $nT u_{nT}$ (resp. $nT w_{nT}$) converges to 0.
Then we obtain for some universal positive constant $C$ and large enough $n$ and $T$
\[
\probaetoile \left( \left\{ \frac{\proba_{\breve{\theta}}(Z^{1:T}\neq z^{*1:T}\given X^{1:T})}{\proba_{\breve{\theta}}(Z^{1:T}=z^{*1:T}\given X^{1:T})} > \epsilon y_{n,T} \right\} \cap \Omega_{n,T}\right)
\leq C nT (u_{nT}+w_{nT}). 
\]
This leads directly to inequality~\eqref{eq:prop:ratio}.
\qed

\subsection{Proof of Theorem~\ref{prop:Gamma_consistency_bis}}
We fix some $\sigma \in \mathfrak{S}_Q$ and study the convergence in $\proba_{\theta^*}-$probability of $\hat{\gamma}_{\sigma(q)\sigma(l)}$ to $\gamma^*_{ql}$ with $\hat{\Gamma}$ as defined by the fixed point equation~\eqref{eq:criticalpoint}, i.e.
\[
\hat{\gamma}_{\sigma(q)\sigma(l)}=\frac{\sum_{t=1}^{T-1} \sum_{i=1}^n
	\proba_{\hat{\theta}_{\sigma}}\left(Z^t_i=q, Z^{t+1}_i=l \given X^{1:T}\right)} {\sum_{t=1}^{T-1} \sum_{i=1}^n
	\proba_{\hat {\theta}_{\sigma}}\left(Z^t_i=q \given X^{1:T} \right)}.
\]
First, let us denote 

\begin{align*}
A_{q,l}&= \frac{1}{n(T-1)}\sum_{t=1}^{T-1} \sum_{i=1}^n \proba_{\hat{\theta}_{\sigma}}\left(Z^t_i=q, Z^{t+1}_i=l
\given X^{1:T}\right),\\
B_{q}&= \frac{1}{n(T-1)}\sum_{t=1}^{T-1} \sum_{i=1}^n \proba_{\hat{\theta}_{\sigma}}\left(Z^t_i=q \given X^{1:T}\right).
\end{align*}
Then we can write the quantity at stake as
\[
\hat{\gamma}_{\sigma(q)\sigma(l)}-\gamma^*_{ql}=\frac{A_{q,l}}{B_{q}}-\gamma^*_{ql} = \frac{A_{q,l}-\alpha^*_q \gamma^*_{ql}}{B_{q}} + \alpha^*_q \gamma^*_{ql} \left(\frac{1}{B_{q}}-\frac{1}{\alpha^*_q} \right) 
\]
to obtain the following upper bound on the probability of interest
\begin{align} \label{eq:decomp_prob_gamma}
\proba_{\theta^*}\left(\left| \hat{\gamma}_{\sigma(q)\sigma(l)}-\gamma^*_{ql} \right| >\epsilon r_{n,T} \frac{\sqrt{\log n}}{\sqrt{nT}} \right) \leq \proba_{\theta^*}\left(\left| \frac{A_{q,l}-\alpha^*_q \gamma^*_{ql}}{B_{q}} \right| > \frac{\epsilon}{2} r_{n,T} \frac{\sqrt{\log n}}{\sqrt{nT}} \right) + \proba_{\theta^*}\left(\alpha^*_q \gamma^*_{ql} \left| \frac{1}{B_{q}}-\frac{1}{\alpha^*_q} \right| > \frac{\epsilon}{2} r_{n,T} \frac{\sqrt{\log n}}{\sqrt{nT}} \right).
\end{align} 
\paragraph*{First term of the right-hand side of~\eqref{eq:decomp_prob_gamma}.} 
For the first term in \eqref{eq:decomp_prob_gamma}, for any $0<\lambda<\delta$ (implying $\lambda<\alpha^*_q$ for any $q \in \Q$),
\begin{align} \label{eq:decomp_prob_gamma_2}
\proba_{\theta^*}\left(\left| \frac{A_{q,l}-\alpha^*_q \gamma^*_{ql}}{B_{q}} \right| > \frac{\epsilon}{2} r_{n,T} \frac{\sqrt{\log n}}{\sqrt{nT}} \right) =& \proba_{\theta^*}\left(\left| \frac{A_{q,l}-\alpha^*_q \gamma^*_{ql}}{B_{q}} \right| > \frac{\epsilon}{2} r_{n,T} \frac{\sqrt{\log n}}{\sqrt{nT}} \given[\Bigg] B_{q} \geq \alpha^*_q - \lambda \right) \proba_{\theta^*}\left( B_{q} \geq \alpha^*_q - \lambda\right) \nonumber \\ &+\proba_{\theta^*}\left(\left| \frac{A_{q,l}-\alpha^*_q \gamma^*_{ql}}{B_{q}} \right| > \frac{\epsilon}{2} r_{n,T} \frac{\sqrt{\log n}}{\sqrt{nT}} \given[\Bigg] B_{q} < \alpha^*_q - \lambda \right) \proba_{\theta^*}\left(B_{q} < \alpha^*_q - \lambda \right) \nonumber \\
\leq& \proba_{\theta^*}\left(\left| A_{q,l}-\alpha^*_q \gamma^*_{ql} \right| > \frac{\epsilon}{2} r_{n,T} \frac{\sqrt{\log n}}{\sqrt{nT}} (\alpha^*_q - \lambda) \right) + \proba_{\theta^*}\left(B_{q} < \alpha^*_q - \lambda \right).
\end{align}
First, we upper bound the probability $\proba_{\theta^*}\left(\left| A_{q,l}-\alpha^*_q \gamma^*_{ql} \right| > \epsilon r_{n,T} \frac{\sqrt{\log n}}{\sqrt{nT}} \right)$ for any $\epsilon>0$, using the following lemma.

\begin{lemme} \label{lem:cv_gamma_intermediaire}
	If $\log(T)=o(n)$, for any $\epsilon>0$, for any sequence $\{r_{n,T}\}_{n,T \geq 1}$ increasing to infinity such that $r_{n,T}=o\left(\sqrt{nT/\log n}\right)$ and any $\eta \in (0,\delta)$, we have for any $\sigma \in \mathfrak{S}_Q$
	\begin{align*}
	&\proba_{\theta^*}\left( \left|\frac{1}{n(T-1)}\sum_{t=1}^{T-1} \sum_{i=1}^n \proba_{\hat{\theta}_{\sigma}}\left(Z^t_i=q, Z^{t+1}_i=l
	\given X^{1:T}\right) - \alpha^*_q \gamma^*_{ql} \right| > \epsilon r_{n,T} \frac{\sqrt{\log n}}{\sqrt{nT}} \right) \leq \proba_{\theta^*} \left( \| \hat{\pi}_{\sigma}- \pi^* \|_{\infty} > v_{n,T} \right) +o(1)
	\end{align*}
	with $v_{n,T}$ a sequence decreasing to $0$ such that $v_{n,T}=o\left(\sqrt{\log(nT)} /n\right)$.
\end{lemme}

Then, for the second term of \eqref{eq:decomp_prob_gamma_2}, notice that $B_{q}= \sum_{l=1}^Q A_{q,l}$ and $\sum_{l=1}^Q \gamma^*_{ql}=1$. We then have, if $\log(T)=o(n)$ and $v_{n,T}=o\left(\sqrt{\log(nT)}/n\right)$, using Lemma~\ref{lem:cv_gamma_intermediaire} again,
\begin{align*}
\proba_{\theta^*}\left(B_{q} < \alpha^*_q - \lambda \right) =& \proba_{\theta^*}\left(B_{q} -\alpha^*_q< - \lambda \right) = \proba_{\theta^*}\left(\sum_{l=1}^Q (A_{q,l} -\alpha^*_q \gamma^*_{ql})< - \lambda \right) \leq \sum_{l=1}^Q \proba_{\theta^*}\left(A_{q,l} -\alpha^*_q \gamma^*_{ql}< - \lambda/Q \right)\\
\leq& \sum_{l=1}^Q \proba_{\theta^*}\left( \left|A_{q,l} -\alpha^*_q \gamma^*_{ql}\right| > \lambda/Q \right)
\leq Q \proba_{\theta^*} \left( \| \hat{\pi}_{\sigma}- \pi^* \|_{\infty} > v_{n,T} \right) + o(1).
\end{align*}
Finally, for the first term of \eqref{eq:decomp_prob_gamma}, if $y_{n,T}$ is such that $1/y_{n,T}=o\left(\sqrt{nT/\log(n)}\right)$, if $v_{n,T}=o\left(\sqrt{\log(nT)}/n\right)$ and as long as $\log(T)=o(n)$, we obtain
\begin{align} \label{eq:gamma_maj_1}
&\proba_{\theta^*}\left(\left| \frac{A_{q,l}-\alpha^*_q \gamma^*_{ql}}{B_{q}} \right| > \frac{\epsilon}{2} r_{n,T} \frac{\sqrt{\log n}}{\sqrt{nT}} \right) 
\leq (Q+1) \proba_{\theta^*} \left( \| \hat{\pi}_{\sigma}- \pi^* \|_{\infty} > v_{n,T} \right) + o(1).
\end{align}

\paragraph*{Second term of the right-hand side of~\eqref{eq:decomp_prob_gamma}.}
For the second term of \eqref{eq:decomp_prob_gamma}, we split it on two complementary events as before. For any $0<\lambda<\delta$, we have
\begin{align*} 
\proba_{\theta^*}\left(\alpha^*_q \gamma^*_{ql} \left| \frac{1}{B_{q}}-\frac{1}{\alpha^*_q} \right| >\frac{\epsilon}{2} r_{n,T} \frac{\sqrt{\log n}}{\sqrt{nT}} \right) 
=& \proba_{\theta^*}\left(\alpha^*_q \gamma^*_{ql} \left| \frac{1}{B_{q}}-\frac{1}{\alpha^*_q} \right| >\frac{\epsilon}{2} r_{n,T} \frac{\sqrt{\log n}}{\sqrt{nT}} \given[\Bigg] B_{q} \geq \alpha^*_q - \lambda \right) \proba_{\theta^*}\left( B_{q} \geq \alpha^*_q - \lambda\right) \nonumber \\ &+\proba_{\theta^*}\left(\alpha^*_q \gamma^*_{ql} \left| \frac{1}{B_{q}}-\frac{1}{\alpha^*_q} \right| >\frac{\epsilon}{2} r_{n,T} \frac{\sqrt{\log n}}{\sqrt{nT}} \given[\Bigg] B_{q} < \alpha^*_q - \lambda \right) \proba_{\theta^*}\left(B_{q} < \alpha^*_q - \lambda \right) \nonumber \\
\leq& \proba_{\theta^*}\left(\alpha^*_q \gamma^*_{ql} \left| \frac{1}{B_{q}}-\frac{1}{\alpha^*_q} \right| >\frac{\epsilon}{2} r_{n,T} \frac{\sqrt{\log n}}{\sqrt{nT}} \given[\Bigg] B_{q} \geq \alpha^*_q - \lambda \right) \proba_{\theta^*}\left( B_{q} \geq \alpha^*_q - \lambda\right) \\ &+ \proba_{\theta^*}\left(B_{q} < \alpha^*_q - \lambda \right). \numberthis \label{eq:cv_gamma_term2}
\end{align*}
We already gave an upper bound on the second term in the right-hand side of~\eqref{eq:cv_gamma_term2}. Let us give one for the first term.
Notice that as $\alpha^*_q\geq \delta$ and if $B_q\geq \alpha^*_q - \lambda \geq \delta-\lambda>0$, we have by the mean value theorem
\[
\left|\frac{1}{B_{q}}-\frac{1}{\alpha^*_q}\right| \leq \frac{1}{(\delta-\lambda)^2} \left|B_{q}-\alpha^*_q\right| .
\]
We can then write for the first term in the right-hand side of~\eqref{eq:cv_gamma_term2}, as long as $\log(T)=o(n)$, for $\{y_{n,T}\}_{n,T \geq 1}$ such that $1/y_{n,T}=o\left(\sqrt{nT/\log n}\right)$ and with $v_{n,T}$ such that $v_{n,T}=o\left(\sqrt{\log(nT)}/n\right)$, still using Lemma~\ref{lem:cv_gamma_intermediaire}
\begin{align*}
&\proba_{\theta^*}\left(\alpha^*_q \gamma^*_{ql} \left| \frac{1}{B_{q}}-\frac{1}{\alpha^*_q} \right| >\frac{\epsilon}{2} r_{n,T} \frac{\sqrt{\log n}}{\sqrt{nT}} \given[\Bigg] B_{q} \geq \alpha^*_q - \lambda \right)
\proba_{\theta^*}\left( B_{q} \geq \alpha^*_q - \lambda\right) \nonumber \\\leq& \proba_{\theta^*}\left(\left| B_{q}-\alpha^*_q \right| > \frac{(\delta-\lambda)^2 \epsilon}{2 \alpha^*_q \gamma^*_{ql} } r_{n,T} \frac{\sqrt{\log n}}{\sqrt{nT}} \right)  = \proba_{\theta^*}\left(\left| \sum_{l=1}^Q (A_{q,l}-\alpha^*_q \gamma^*_{ql}) \right| > \frac{(\delta-\lambda)^2 \epsilon}{2 \alpha^*_q \gamma^*_{ql} }  r_{n,T} \frac{\sqrt{\log n}}{\sqrt{nT}} \right) \nonumber \\
\leq& \sum_{l=1}^Q \proba_{\theta^*}\left(\left| A_{q,l}-\alpha^*_q \gamma^*_{ql} \right| > \frac{(\delta-\lambda)^2 \epsilon}{2 \alpha^*_q \gamma^*_{ql} Q}  r_{n,T} \frac{\sqrt{\log n}}{\sqrt{nT}} \right) 
\leq Q \proba_{\theta^*} \left( \| \hat{\pi}_{\sigma}- \pi^* \|_{\infty} > v_{n,T} \right) + o(1).
\end{align*}
We finally obtain for the second term of the right-hand side of~\eqref{eq:decomp_prob_gamma} 
\begin{align} \label{eq:gamma_maj_2}
\proba_{\theta^*}\left(\alpha^*_q \gamma^*_{ql} \left| \frac{1}{B_{q}}-\frac{1}{\alpha^*_q} \right| >\frac{\epsilon}{2} r_{n,T} \frac{\sqrt{\log n}}{\sqrt{nT}} \right) 
\leq 2Q \proba_{\theta^*} \left( \| \hat{\pi}_{\sigma}- \pi^* \|_{\infty} > v_{n,T} \right) + o(1).
\end{align} 
We conclude the proof by summing the upper bounds obtained in \eqref{eq:gamma_maj_1} and \eqref{eq:gamma_maj_2} 
\begin{align*}
\proba_{\theta^*}\left(\left| \hat{\gamma}_{\sigma(q)\sigma}(l)-\gamma^*_{ql} \right| > \epsilon r_{n,T} \frac{\sqrt{\log n}}{\sqrt{nT}} \right) 
\leq & (3Q+1)
\proba_{\theta^*} \left( \| \hat{\pi}_{\sigma}- \pi^* \|_{\infty} > v_{n,T} \right) + o(1)
\end{align*}
and by noticing that $\proba_{\theta^*}(\| \hat{\Gamma}_{\sigma}- \Gamma^* \|_{\infty} > \epsilon r_{n,T} \sqrt{\log n}/\sqrt{nT}) \leq \sum_{1 \leq q,l \leq Q} \proba_{\theta^*}(| \hat{\gamma}_{\sigma(q)\sigma(l)}- \gamma^*_{ql} | > \epsilon r_{n,T} \sqrt{\log n}/\sqrt{nT}) $.
\qed

\subsection{Proof of Corollary~\ref{cor:Gamma_consistency}}
Denoting by $\sigma_{n,T}$ the permutation minimizing the distance between $\hat{\pi}$ (permuted) and $\pi^*$ for every $(n,T) \in \n \times \T$, i.e. $\sigma_{n,T}=\argmin_{\sigma \in \mathfrak{S}_Q} \|\hat{\pi}_{\sigma} - \pi^*\|_{\infty}$, we apply Theorem~\ref{prop:Gamma_consistency_bis} to $\hat{\theta}_{\sigma_{n,T}}$ in order to get
\begin{align*}
\proba_{\theta^*}\left( \min_{\sigma \in \mathfrak{S}_Q} \|  \hat{\Gamma}_{\sigma} - \Gamma^* \|_{\infty} > \epsilon r_{n,T} \frac{\sqrt{\log n}}{\sqrt{nT}} \right) \leq& \proba_{\theta^*}\left( \|  \hat{\Gamma}_{\sigma_{n,T}} - \Gamma^* \|_{\infty} > \epsilon r_{n,T} \frac{\sqrt{\log n}}{\sqrt{nT}}  \right)\\ \leq& Q^2(3Q+1)
\proba_{\theta^*} \left( \min_{\sigma \in \mathfrak{S}_Q } \| \hat{\pi}_{\sigma} - \pi^* \|_{\infty} > v_{n,T} \right) + o(1) \xrightarrow[n,T \to \infty]{}0,
\end{align*}
\qed
\subsection{Proof of Theorem~\ref{prop:vitesse_cv_M_variationnel}}
We use the following lemma, that states that the quantity we optimize in the VEM algorithm and the log-likelihood are asymptotically equivalent.
\begin{lemme} \label{lem:variationel_asymp_eq}
	We have the following inequality $\proba_{\theta^*}$-a.s.
	\[
	\sup_{\theta \in \Theta} \left| \frac{2}{n (n-1) T} \mathcal{J}(\hat{\chi}(\theta),\theta) - \frac{2}{n (n-1) T} \ell(\theta) \right| \leq \frac{2 \log(1/\delta)}{n-1}. 
	\]
\end{lemme}

	We have that for any $\epsilon>0$, for $n$ and $T$ large enough,
	\[
	\proba_{\theta^*} \left( \sup_{\theta \in \Theta} \left| \frac{2}{n(n-1)T} \mathcal{J}(\hat{\chi}(\theta),\theta) - \frac{2}{n(n-1)T} \ell(\theta) \right| > \frac{\epsilon r_{n,T}}{\sqrt{n}} \right) \leq \proba_{\theta^*} \left( \frac{2 \log(1/\delta)}{n-1}  > \frac{\epsilon r_{n,T}}{\sqrt{n}} \right) =0
	\]
	We then conclude by combining this result with Theorem~\ref{prop:vitesse_cv_M}.
	\qed

\subsection{Proof of Corollary~\ref{cor:vitesse_pitilde_bis}}
	This is a direct consequence of Theorem~\ref{prop:vitesse_cv_M_variationnel} and Lemma~\ref{lem:cv_M_to_cv_pi_bis} applied with the functions $F_{n,T}=\frac{2}{n(n-1)T} \mathcal{J}(\hat{\chi}(\cdot),\cdot)$. 
\qed

\subsection{Proof of Theorem~\ref{prop:consistency_gamma_tilde}}
	This proof is quite similar to that of Theorem~\ref{prop:Gamma_consistency_bis}. We fix some $\sigma \in \mathfrak{S}_Q$ and study the convergence in $\proba_{\theta^*}-$probability of $\tilde{\gamma}_{\sigma(q)\sigma(l)}$ to $\gamma^*_{ql}$ with $\tilde{\Gamma}$ as defined by the fixed point equation~\eqref{eq:fixed_point_variationnel}, i.e.
	\[
	\tilde{\gamma}_{\sigma(q)\sigma(l)} = \frac{\sum_{i=1}^n \sum_{t=1}^{T-1} \hat{\eta}^t_{iql}(\tilde{\theta}_{\sigma})}{\sum_{i=1}^n \sum_{t=1}^{T-1} \hat{\tau}^{t}_{iq}(\tilde{\theta}_{\sigma})}.
	\]
	First, let us denote 
	
	\begin{align*}
	A_{q,l}&= \frac{1}{n(T-1)} \sum_{i=1}^n \sum_{t=1}^{T-1} \hat{\eta}^t_{iql}(\tilde{\theta}_{\sigma}) = \frac{1}{n(T-1)} \sum_{i=1}^n \sum_{t=1}^{T-1} \mathbb{Q}_{\hat{\chi}(\tilde \theta_{\sigma})}(Z_i^{t}=q,Z_i^{t+1}=l),\\
	B_{q}&= \frac{1}{n(T-1)}\sum_{i=1}^n\sum_{t=1}^{T-1}\hat{\tau}^{t}_{iq}(\tilde{\theta}_{\sigma}) = \frac{1}{n(T-1)}\sum_{i=1}^n\sum_{t=1}^{T-1} \mathbb{Q}_{\hat{\chi}(\tilde \theta_{\sigma})}(Z_i^{t}=q).
	\end{align*}
	Then we can write the quantity at stake as
	\[
	\tilde{\gamma}_{\sigma(q)\sigma(l)}-\gamma^*_{ql}=\frac{A_{q,l}}{B_{q}}-\gamma^*_{ql} = \frac{A_{q,l}-\alpha^*_q \gamma^*_{ql}}{B_{q}} + \alpha^*_q \gamma^*_{ql} \left(\frac{1}{B_{q}}-\frac{1}{\alpha^*_q} \right). 
	\]
	We follow the line of the proof of Theorem~\ref{prop:Gamma_consistency_bis}, using Lemma~\ref{lem:cv_gamma_intermediate_variationnel} below instead of Lemma~\ref{lem:cv_gamma_intermediaire} in order to obtain the result.
	\begin{lemme}\label{lem:cv_gamma_intermediate_variationnel}
		For any $\epsilon>0$, for any sequence $\{r_{n,T}\}_{n,T \geq 1}$ increasing to infinity such that $r_{n,T}=o\left(\sqrt{nT/\log n}\right)$ and any $\eta \in (0,\delta)$, we have for any $\sigma \in \mathfrak{S}_Q$
		\begin{align*}
		&\proba_{\theta^*}\left( \left|\frac{1}{n(T-1)} \sum_{i=1}^n \sum_{t=1}^{T-1} \mathbb{Q}_{\hat{\chi}(\tilde \theta_{\sigma})}(Z_i^{t}=q,Z_i^{t+1}=l) - \alpha^*_q \gamma^*_{ql} \right| > \epsilon r_{n,T} \frac{\sqrt{\log n}}{\sqrt{nT}} \right) 
		\leq 2 \proba_{\theta^*} \left( \| \tilde{\pi}_{\sigma}- \pi^* \|_{\infty} > v_{n,T} \right) + o(1)
		\end{align*}
		with $v_{n,T}$ a sequence decreasing to $0$ such that $v_{n,T}=o(\sqrt{\log(nT)} /n)$.
	\end{lemme}
	\qed

\section*{Acknowledgement}
Work partly supported by the grant ANR-18-CE02-0010 of the French National Research Agency ANR (project EcoNet).

\bibliographystyle{abbrvnat}
\bibliography{biblidyn}

\newpage
\appendix

\section{Proofs of main results for the finite time case} \label{appendix:proofs_fixed_T}
	\subsection{Proof of Corollary~\ref{cor:vitesse_cv_pi_fixedT}}
	When the number of time steps is fixed and the connection probabilities vary over time, the conditional log-likelihood is
	\begin{align*}
	\ell^T_{c} (\theta ; Z^{1:T})&= \sum_{t=1}^T \sum_{1\le i<j \le n} \Xtij \log \pi^t_{Z^{t}_i Z^{t}_j} + (1-\Xtij) \log (1-\pi^t_{Z^{t}_i Z^{t}_j})
	\end{align*}
	and the likelihood $\ell^T(\theta)$ is defined as in~\eqref{eq:log-likelihood} with $\ell^T_{c}(\cdot)$ instead of $\ell_{c}(\cdot)$. The maximum likelihood estimator is then
	\[
	\hat{\theta}=(\hat{\Gamma},\hat{\pi}^{1:T})= \argmax_{\theta \in \Theta^T} \ell^T(\theta).
	\]
	As before, we denote the normalized log-likelihood $M_{n,T}(\Gamma,\pi^{1:T}) = 2/(n(n-1)T) \ell^T(\theta)$.
	We introduce the following limiting quantity
	\begin{align*}
	\mathbb{M}^T(\pi^{1:T}) &= \frac{1}{T} \sum_{t=1}^T \mathbb{M}(\pi^t) = \frac{1}{T} \sum_{t=1}^T \sup_{A \in \mathcal{A}} \mathbb{M}(\pi^t,A) .
	\end{align*}
	We follow the lines of the proof of Theorem~\ref{prop:vitesse_cv_M} in order to prove that we have for any sequence $y_n \to +\infty$, for all $\epsilon>0$
	\begin{align} \label{eq:cv_M_fixed_T}
	\proba_{\theta^*}\left(\sup_{(\Gamma,\pi^{1:T})\in \Theta^T}\left| M_{n,T}(\Gamma,\pi^{1:T}) - \mathbb{M}^T(\pi^{1:T}) \right| >
	\frac{\epsilon y_n}{\sqrt{n}}\right) \mathop{\longrightarrow}_{n\to +\infty} 0 . 
	\end{align}
	Choosing $y_n=r_n^2$, we then use Lemma~\ref{lem:cv_M_to_cv_pi_bis} to conclude that, as $r_n^2/\sqrt{n}=o(1)$ by assumption, for any $\epsilon>0$,
	\[
	\proba_{\theta^*}\left(\min_{\sigma^{1},\ldots,\sigma^{T} \in \mathfrak{S}_Q}\|\hat{\pi}_{\sigma^{1:T}}^{1:T}-\pi^{*1:T} \|_{\infty}> \epsilon r_n/n^{1/4} \right) \xrightarrow[n \rightarrow \infty]{} 0.
	\]
	In particular, for every $t \in\T$, $\hat{\pi}^{t}$ converges in $\proba_{\theta^*}$-probability to $\pi^{*t}$ up to label switching. 
	Then, let us prove that on the event $\{ \min_{\sigma^1, \ldots,\sigma^T \in \mathfrak{S}_Q}  \| \hat{\pi}^{1:T} - \pi_{\sigma^{1:T}}^{*1:T} \|_{\infty} \leq \epsilon r_n n^{-1/4} \}$ (whose probability converges to $1$), for $n$ large enough, the permutation $\sigma^t$ 	minimizing the distance between $\pi^{*t}$ and $\hat{\pi}_{\sigma^t}^t$ is the same for every $t \in \T$. 
	We consider $n$ large enough such that $\epsilon r_n n^{-1/4}< \min_{1\leq q \neq l \leq Q} |\pi^{*}_{qq}-\pi^*_{ll}|/4$. 
	Denoting by $\sigma_m^1,\ldots,\sigma_m^T$ the permutations (depending on $n$) minimizing $\| \hat{\pi}^{1:T} - \pi_{\sigma^{1:T}}^{*1:T} \|_{\infty}$, 
	we have that, for any $1\leq t \neq t' \leq T$, if some $q,l \in \Q$ are such that $\sigma_m^t(q)=\sigma_m^{t'}(l)$, then
	\[
	\hat{\pi}^{t}_{\sigma_m^t(q)\sigma_m^t(q)}=\hat{\pi}^{t}_{\sigma_m^{t'}(l)\sigma_m^{t'}(l)}=\hat{\pi}^{t'}_{\sigma_m^{t'}(l)\sigma_m^{t'}(l)}
	\]
	and on the event we consider
	\begin{align*}
	| \pi^{*t}_{qq} -\pi^{*t}_{ll} | &= | \pi^{*t}_{qq} -\pi^{*t'}_{ll} | = | \pi^{*t}_{qq} - \hat{\pi}^{t}_{\sigma_m^t(q)\sigma_m^t(q)} + \hat{\pi}^{t'}_{\sigma_m^{t'}(l)\sigma_m^{t'}(l)}-\pi^{*t'}_{ll} |
	\leq |  \pi^{*t}_{qq} - \hat{\pi}^{t}_{\sigma_m^t(q)\sigma_m^t(q)}| + |\hat{\pi}^{t'}_{\sigma_m^{t'}(l)\sigma_m^{t'}(l)}-\pi^{*t'}_{ll} | \\
	&\leq  2 \epsilon r_n n^{-1/4} < \min_{1\leq q\neq l \leq Q} |\pi^*_{qq}-\pi^*_{ll}|/2,
	\end{align*}
	implying that $q=l$. This means that on this event, the permutation $\sigma_m^t$ minimizing the distance between $\pi^{*t}$ and $\hat{\pi}_{\sigma^t}^t$ is the same for every $t \in \T$.
	We can conclude that
	\begin{align*}
	\proba_{\theta^*}\left(\min_{\sigma \in \mathfrak{S}_Q}\|\hat{\pi}_{\sigma}^{1:T}-\pi^{*1:T} \|_{\infty}> \epsilon r_n/n^{1/4} \right)
	&=1-\proba_{\theta^*}\left(\min_{\sigma \in \mathfrak{S}_Q}\|\hat{\pi}_{\sigma}^{1:T}-\pi^{*1:T} \|_{\infty}\leq \epsilon r_n/n^{1/4} \right)
	\xrightarrow[n \rightarrow \infty]{} 0.
	\end{align*}
	
	\qed
	
\subsection{Proof of Theorem~\ref{prop:ratio_alt}}
First, let us introduce some notations, as in the proof of Theorem~\ref{prop:ratio}. For any fixed configuration $z^{*1:T} \in \Omega_{\eta}$, we define for any configuration $z^{1:T}$ and any parameter $\theta$
\begin{equation*}
D_{n,T}(z^{1:T}, \pi^{1:T}) \coloneqq \left\{(i,j,t) \in I_{n,T}; \pi^t_{z^{t}_i z^{t}_j} \neq \pi^t_{z^{*t}_i z^{*t}_j} \right\}
\end{equation*}
and for any $1\leq t \leq T$
\begin{equation*}
D^t_{n,T}(z^{t}, \pi^t) \coloneqq \left\{(i,j) \in \n^2; i<j \textrm{ and } \pi^t_{z^{t}_i z^{t}_j} \neq \pi^t_{z^{*t}_i z^{*t}_j} \right\},
\end{equation*}
and as before, we abbreviate to $D^*$ (resp. $\breve D$), the set $D_{n,T}(z^{1:T}, \pi^{*1:T}) $
(resp. $D_{n,T}(z^{1:T},\breve \pi^{1:T})$). We also introduce for any $q,l,q',l' \in \Q$ the quantities $F_{qlq'l'}$, $F_{ql}$, $G_{qlq'l'}$ and $G_{ql}$ as before, accordingly to this definition of $D_{n,T}(z^{1:T}, \pi^{1:T})$. Finally, we introduce for any $t \in \T$ and $q,l,q',l' \in \Q$ the quantities
\begin{align*}
F^t_{qlq'l'}&=  F^t_{qlq'l'} (z^{t},z^{*t})\coloneqq \{(i,j) \in \n^2; i<j  \textrm{ and } z_i^t=q,z_j^t=l,z_i^{*t}=q',z_j^{*t}=l'\}\\
F_{ql}^t&=  F^t_{ql} (z^{t})\coloneqq \cup_{1\leq q',l' \leq Q} F^t_{qlq'l'} =\{(i,j)\in\n^2 ; i<j  \textrm{ and }  z_i^t=q,z_j^t=l\}\\
G^t_{qlq'l'} &= G^t_{qlq'l'} (z^{t}, z^{*t},\pi^{*t},\breve{\pi}^t)\coloneqq (D^{*t}\cup \breve{D}^t) \cap F^t_{qlq'l'}\\
&=\{(i,j)\in\n^2 ;  i<j  \textrm{ and } z_i^t=q,z_j^t=l,z_i^{*t}=q',z_j^{*t}=l' \text{ and } (\pi^{*t}_{z_i^t z_j^t} \neq \pi^{*t}_{z_i^{*t} z_j^{*t}} \text{ or } \breve \pi^t_{z_i^t z_j^t} \neq \breve \pi^t_{z_i^{*t} z_j^{*t}} )\} \\
G^t_{ql} &= G^t_{ql} (z^{t}, z^{*t},\pi^{*t},\breve{\pi}^t)\coloneqq (D^{*t}\cup \breve{D}^t) \cap F^t_{ql} \\
&=\{(i,j)\in\n^2 ;  i<j  \textrm{ and } z_i^t=q,z_j^t=l \text{ and } (\pi^{*t}_{z_i^t z_j^t} \neq \pi^{*t}_{z_i^{*t} z_j^{*t}} \text{ or } \breve{\pi}^t_{z_i^t z_j^t} \neq \breve{\pi}^t_{z_i^{*t} z_j^{*t}} )\} .
\end{align*} 
Note that we can get an equivalent of Lemma~\ref{minor} with a similar proof that gives that for any configuration $z^{*1:T}$ in $\Omega_{\eta}$, for any configuration $z^{1:T}$ and any $\theta \in \Theta^{T}$,
\[
\left| D_{n,T}(z^{1:T},\pi^{1:T}) \right| \geq \frac{\gamma^2}{4}n r.
\]
In the same way, we have an equivalent of Lemma~\ref{major} (with a similar proof) that gives that for any $z^{t}$ and $z^{*t}$ two configurations at time $t$ such that
$ \|z^{t}-z^{*t}\|_0=r(t)$ and any parameter $\pi^t=(\pi^t_{ql})_{1\leq q,l\leq Q}$, we have
\begin{equation} \label{eq:equivalent_lemma2}
D^t_{n,T}(z^{t}, \pi^t)  \subset D^t_{n,T}(z^{t}) \coloneqq 
\left\{ (i,j) \in \n^2\times \T; (z_i^t ,z_j^t) \neq (z_i^{*t}, z_j^{*t}) \right\} 
\text{ and }  \left| D^t_{n,T}(z^{t}) \right| \leq 2nr(t).
\end{equation}
Going back to the proof of Theorem~\ref{prop:ratio_alt}, we follow the line of that of Theorem~\ref{prop:ratio}, with a few changes. We get the same decomposition as in equation~\eqref{eq:decomp_r}, replacing $\pi$ by $\pi^1, \ldots, \pi^T$ in the definitions of $U_1$, $U_2$ and $U_3$, and replacing the event $\Omega_{n,T}$ by $\Omega_n=\{\| \hat{\pi}^{1:T}- \pi^{*1:T}\|_{\infty}\leq v_n\}$. For $U_1$, the proof does not change. 
For $U_2$, we write (instead of \eqref{eq:prop3_secondterm})
\begin{align*}
|U_2| \leq& \left| \sum_{(i,j,t)\in D^*\cup \breve{D}} \sum_{1\leq q,l \leq Q} \frac{ \breve{\pi}^t_{ql}-\pi^{*t}_{ql} }{\pi^{*t}_{ql}(1-\pi^{*t}_{ql})} (X_{ij}^t-\pi^{*t}_{ql}) \ind_{z_i^t=q,z_j^t=l} \right|\\
\leq& \sum_{t=1}^T \sum_{1\leq q,l \leq Q} \left| \frac{ \breve{\pi}^t_{ql}-\pi^{*t}_{ql} }{\pi^{*t}_{ql}(1-\pi^{*t}_{ql})} \sum_{(i,j)\in G^t_{ql}} (X_{ij}^t-\pi^{*t}_{ql}) \right|
\leq \sum_{t=1}^T \sum_{1\leq q,l \leq Q} \frac{\left| \breve{\pi}^t_{ql}-\pi^{*t}_{ql}\right| }{\pi^{*t}_{ql}(1-\pi^{*t}_{ql})} \left| \sum_{q',l'} \sum_{(i,j)\in G^t_{qlq'l'}} (X_{ij}^t-\pi^{*t}_{ql}) \right|\\
\leq& \sum_{t=1}^T \sum_{1\leq q,l \leq Q} \frac{\left| \breve{\pi}^t_{ql}-\pi^{*t}_{ql}\right| }{\pi^{*t}_{ql}(1-\pi^{*t}_{ql})} \left| \sum_{q',l'} \sum_{(i,j)\in G^t_{qlq'l'}} (X_{ij}^t-\pi^{*t}_{q'l'}) \right| + \sum_{t=1}^T \sum_{1\leq q,l \leq Q} \frac{\left| \breve{\pi}^t_{ql}-\pi^{*t}_{ql}\right| }{\pi^{*t}_{ql}(1-\pi^{*t}_{ql})} \left| \sum_{q',l'}  (\pi^{*t}_{q'l'}-\pi^{*t}_{ql}) |G^t_{qlq'l'}| \right|.
\end{align*}
For every $u>0$, we thus have
\begin{align} \label{eq:decomp_U2_bis}
\probaetoile \left( \{ |U_2| >u \} \cap \Omega_n \right) \leq& \sum_{t=1}^T \probaetoile \left( \left\{\sum_{1\leq q,l \leq Q} \frac{\left| \breve{\pi}^t_{ql}-\pi^{*t}_{ql}\right| }{\pi^{*t}_{ql}(1-\pi^{*t}_{ql})} \left| \sum_{1\leq q',l' \leq Q} \sum_{(i,j) \in G^t_{qlq'l'}} (X_{ij}^t-\pi^{*t}_{q'l'}) \right| > \frac{u}{2T} \right\} \cap \Omega_n \right) \nonumber\\
&+ \sum_{t=1}^T \probaetoile \left( \left\{\sum_{1\leq q,l \leq Q} \frac{\left| \breve{\pi}^t_{ql}-\pi^{*t}_{ql}\right| }{\pi^{*t}_{ql}(1-\pi^{*t}_{ql})} \left| \sum_{1\leq q',l' \leq Q}  (\pi^{*t}_{q'l'}-\pi^{*t}_{ql}) |G^t_{qlq'l'}| \right| > \frac{u}{2T} \right\} \cap \Omega_n \right).
\end{align}
We start by dealing with the first term of \eqref{eq:decomp_U2_bis}. Notice that on the event $\Omega_{n}$, we have $\left| \breve{\pi}^t_{ql}-\pi^{*t}_{ql}\right|/(\pi^{*t}_{ql}(1-\pi^{*t}_{ql}))\leq v_{n}/\zeta^2$ for every $q,l \in \Q$. As the set $G^t_{ql}$ is random (because $\breve{D}^t$ is random), we write for every $t \in \T$, using~\eqref{eq:equivalent_lemma2},
\begin{align*}
& \probaetoile\left(\left\{  \sum_{1\leq q,l \leq Q}  \frac{\left|\breve{\pi}^t_{ql}-\pi^{*t}_{ql}\right| }{\pi^{*t}_{ql}(1-\pi^{*t}_{ql})}\left| \sum_{1 \leq q',l' \leq Q} \sum_{(i,j) \in G^t_{ql}}  (X_{ij}^t-\pi^{*t}_{q'l'}) \right| > \frac{u}{2T} \right\} \cap \Omega_{n} \right) \\
\leq &  \probaetoile\left( \sum_{1\leq q,l \leq Q} \left| \sum_{1\leq q',l' \leq Q} \sum_{(i,j) \in G^t_{ql}} (X_{ij}^t-\pi^{*t}_{q'l'}) \right| > \frac{u \zeta^2}{2Tv_{n}} \right)
\leq \sum_{D \subset D^t_{n,T}(z^t)}\probaetoile\left( \sum_{1\leq q,l \leq Q} \left| \sum_{1\leq q',l' \leq Q} \sum_{(i,j) \in F^t_{ql}\cap D} (X_{ij}^t-\pi^{*t}_{q'l'}) \right| > \frac{u \zeta^2}{2Tv_{n}} \right)
\end{align*}
where now $D$ is a deterministic set. By a union bound and Hoeffding's inequality, we have for any $D \subset D^t_{n,T}(z^t)$

\begin{align*}
\probaetoile\left( \sum_{1\leq q,l \leq Q} \left| \sum_{1\leq q',l' \leq Q} \sum_{(i,j) \in F^t_{ql}\cap D}  (X_{ij}^t-\pi^{*t}_{q'l'}) \right| > \frac{u \zeta^2}{2Tv_{n}} \right)
\leq&Q^2 \max_{1\leq q,l \leq Q} \probaetoile\left( \left| \sum_{1\leq q',l' \leq Q} \sum_{(i,j) \in F^t_{ql}\cap D}  (X_{ij}^t-\pi^{*t}_{q'l'}) \right| > \frac{u \zeta^2}{2Tv_{n} Q^2} \right)\\
\leq& 2 Q^2 \exp\left(-\frac{2 u^2 \zeta^4}{4 T^2 v_{n}^2 Q^4}\frac{1}{|D|}\right).
\end{align*}
This leads to, for the first term of \eqref{eq:decomp_U2_bis},
\begin{align*}
&\sum_{t=1}^T \probaetoile\left(\left\{  \sum_{1\leq q,l \leq Q} \left| \frac{(\breve{\pi}^t_{ql}-\pi^{*t}_{ql}) }{\pi^{*t}_{ql}(1-\pi^{*t}_{ql})} \right| \left| \sum_{1\leq q',l' \leq Q} \sum_{(i,j) \in G^t_{ql}}  (X_{ij}^t-\pi^{*t}_{q'l'}) \right| > \frac{u}{2T} \right\} \cap \Omega_{n} \right)\\ 
\leq& \sum_{t=1}^T \sum_{D \subset D^t_{n,T}(z^t)} 2 Q^2 \exp\left(-\frac{2 u^2 \zeta^4}{4 T^2 v_{n}^2 Q^4}\frac{1}{|D|}\right)
\leq \sum_{t=1}^T \sum_{k=1}^{2nr(t)}\sum_{D \subset D^t_{n,T}(z^t);|D|=k} 2 Q^2 \exp\left(-\frac{2 u^2 \zeta^4}{4 T^2 v_{n}^2 Q^4}\frac{1}{k}\right)\\
\leq& 2 Q^2 \sum_{t=1}^T \exp\left(-\frac{u^2 \zeta^4}{4 T^2 v_{n}^2 Q^4 nr(t)}\right) (2nr(t))^{2nr(t)+1}
\leq 2 Q^2 T \exp\left(-\frac{u^2 \zeta^4}{4 T^2 v_{n}^2 Q^4 nr}\right) (2nr)^{2nr+1}.
\end{align*}
For the second term of \eqref{eq:decomp_U2_bis}, we get from a union bound and from~\eqref{eq:equivalent_lemma2} that
\begin{align*}
&\sum_{t=1}^T \probaetoile \left( \left\{ \sum_{1\leq q,l \leq Q} \left| \frac{(\breve{\pi}^t_{ql}-\pi^{*t}_{ql}) }{\pi^{*t}_{ql}(1-\pi^{*t}_{ql})} \right| \left| \sum_{1\leq q',l' \leq Q}  (\pi^{*t}_{q'l'}-\pi^{*t}_{ql}) |G^t_{qlq'l'}| \right|  > \frac{u}{2T}\right\} \cap \Omega_{n} \right)\\
\leq& Q^2 \sum_{t=1}^T \max_{1\leq q,l \leq Q} \probaetoile \left( \left| \sum_{1\leq q',l'\leq Q}  (\pi^{*t}_{q'l'}-\pi^{*t}_{ql}) |G^t_{qlq'l'}| \right|  > \frac{u\zeta^{2}}{2Tv_{n}Q^2} \right)
\leq Q^2 T \probaetoile \left( 2nr  > \frac{u\zeta^{2}}{2v_{n}TQ^2} \right).
\end{align*}
Finally, we have the following upper bound for $U_2$ 
\begin{align*}
\probaetoile \left( \Omega_{n} \cap \left\{|U_2| >r \log(nT)\right\} \right) 
\leq& 2 Q^2 T \exp\left(-\frac{r \zeta^4 (\log(nT))^2 }{4 Q^4 T^2 v_{n}^2 n}\right) (2nr)^{2nr+1}
+ Q^2 T \probaetoile \left( v_{n} > \frac{\zeta^{2}\log(nT)}{4 Q^2 T n} \right).
\end{align*}	
For the third term $U_3$, denoting $G^{*t}_{ql}=\cup_{1\leq q',l'\leq Q} G^t_{ql}= \{(i,j) \in D^{*t}\cup \breve{D}^t; z_i^{*t}=q, z_j^{*t}=l \}$, we have
\begin{align*}
U_3 =& \sum_{1\leq q,l \leq Q}\sum_{(i,j,t)\in D^*\cup \breve{D}} \left( (\pi^{*t}_{ql} -X_{ij}^t) \log\left[1-\frac{(\breve{\pi}^t_{ql}-\pi^{*t}_{ql})}{(1-\pi^{*t}_{ql})} \right] + (X_{ij}^t-\pi^{*t}_{ql}) \log\left[1+\frac{(\breve{\pi}^t_{ql}-\pi^{*t}_{ql})}{\pi^{*t}_{ql}} \right] \right) \ind_{z_i^{*t}=q, z_j^{*t}=l}\\
&+\sum_{1\leq q,l \leq Q}\sum_{(i,j,t)\in D^*\cup \breve{D}} \left( (1-\pi^{*}_{ql}) \log\left[1-\frac{(\breve{\pi}^t_{ql}-\pi^{*t}_{ql})}{(1-\pi^{*t}_{ql})} \right] + \pi^{*t}_{ql} \log\left[1+\frac{(\breve{\pi}^t_{ql}-\pi^{*t}_{ql})}{\pi^{*t}_{ql}} \right] \right) \ind_{z_i^{*t}=q, z_j^{*t}=l}\\
=& \sum_{t=1}^T \sum_{1\leq q,l \leq Q} \left( \log\left[1+\frac{(\breve{\pi}^t_{ql}-\pi^{*t}_{ql})}{\pi^{*t}_{ql}} \right]- \log\left[1-\frac{(\breve{\pi}^t_{ql}-\pi^{*t}_{ql})}{(1-\pi^{*t}_{ql})} \right]\right) \sum_{(i,j) \in G^{*t}_{ql}} (X_{ij}^t-\pi^{*t}_{ql}) \\
&+ \sum_{t=1}^T \sum_{1\leq q,l \leq Q} |G^{*t}_{ql}| \left( (1-\pi^{*t}_{ql}) \log\left[1+\frac{(\breve{\pi}^t_{ql}-\pi^{*t}_{ql})}{\pi^{*t}_{ql}} \right]+ \pi^{*t}_{ql} \log\left[1-\frac{(\breve{\pi}^t_{ql}-\pi^{*t}_{ql})}{(1-\pi^{*t}_{ql})} \right]\right).
\end{align*}
Then, we have on the event $\Omega_{n}$ and for $n$ large enough such that $|(\breve{\pi}^t_{ql}-\pi^{*t}_{ql})/\pi^{*t}_{ql}| \leq 1/2$ and $|(\breve{\pi}^t_{ql}-\pi^{*t}_{ql})/(1-\pi^{*t}_{ql})| \leq 1/2$ for every $q$ and $l$, using the fact that $|\log(1+x)| \leq 2|x|$ for $x \in [-1/2,1/2]$,
\begin{align*}
|U_3| \leq& \sum_{t=1}^T 4 \frac{v_{n}}{\zeta} \sum_{1\leq q,l \leq Q} \left| \sum_{(i,j) \in G^{*t}_{ql}} (X_{ij}^t-\pi^{*t}_{ql}) \right| 
+ \sum_{t=1}^T 4 \frac{v_{n}}{\zeta} \sum_{1\leq q,l \leq Q} |G^{*t}_{ql}|.
\end{align*}
Then, for every $u>0$,
\begin{align} \label{eq:decomp_U3_bis}
\probaetoile \left( \Omega_{n} \cap \left\{ \left| U_3 \right| > u \right\} \right)
\leq& \sum_{t=1}^T \probaetoile \left( \sum_{1\leq q,l \leq Q} \left| \sum_{(i,j) \in Q^{*t}_{ql}} (X_{ij}^t-\pi^{*t}_{ql}) \right| > \frac{u \zeta}{8 v_{n}T} \right)
+ \sum_{t=1}^T \probaetoile \left(v_{n}\sum_{1\leq q,l \leq Q} |G^{*t}_{ql}| > \frac{u \zeta}{8T} \right).
\end{align}
For the first term of \eqref{eq:decomp_U3_bis}, using Hoeffding's inequality as before,
\begin{align*}
\sum_{t=1}^T \probaetoile \left( \sum_{q,l} \left| \sum_{(i,j) \in G^{*t}_{ql}} (X_{ij}^t-\pi^{*t}_{ql}) \right| > u\zeta/(8 v_{n}T) \right) \leq& \sum_{t=1}^T \sum_{k =1}^{2nr(t)} \sum_{D \subset D^t_{n,T}(z^t);|D|=k} \probaetoile \left( \sum_{q,l} \left| \sum_{(i,j) \in F^{*t}_{ql} \cap D} (X_{ij}^t-\pi^{*t}_{ql}) \right| > u\zeta/(8 v_{n}T) \right)\\
\leq&  2 Q^2 T \exp \left(-\frac{u^2 \zeta^2}{8^2 T^2 Q^4 v_{n}^2 nr} \right) (2nr)^{2nr+1},
\end{align*}
and for the second term of \eqref{eq:decomp_U3_bis},
\begin{align*}
\sum_{t=1}^T \probaetoile \left( v_{n} \sum_{q,l} |G^{*t}_{ql}| > \frac{u \zeta}{8T} \right) \leq T \probaetoile \left( v_{n}  > \frac{u \zeta}{16 T nr} \right).
\end{align*}
Finally, we have the following upper bound for $U_3$ 
\begin{align*}
\probaetoile \left( \Omega_{n} \cap \left\{|U_3| > r \log(nT)\right\} \right) \leq& 2 Q^2 T \exp \left(-\frac{r \zeta^2 (\log(nT))^2}{8^2 T^2 Q^4 v_{n}^2 n} \right) (2nr)^{2nr+1} + T \probaetoile \left( v_{n}  > \frac{\zeta \log(nT)}{16 T n} \right).
\end{align*}
Now we choose the sequence $v_n$ such that $v_{n}=o(\sqrt{\log n}/n)$ which is sufficient to imply that the quantities $\probaetoile \left( v_{n} > \zeta^{2}\log(nT)/(4 Q^2 T n) \right)$ and $\probaetoile \left( v_{n}  > \zeta \log(nT)/(16 T n) \right)$ vanish as $n$ increases and we gather the three upper bounds. 
For large enough values of $n$ and with $C_1$, $C_2$, $C_3$, $C_4$ and $\kappa$ positive constants only depending on $Q$, $\zeta$, $K^*$ and $T$, we then have 
\begin{align*}
&\probaetoile \left( \left\{ U_1+U_2-U_3 > - \log(1/(\epsilon y_{n})) - 3r\log (nT) \right\} \cap \Omega_{n}\right)\\
\leq& \exp \left[ (\log(1/(\epsilon y_{n}))+ 5 r \log(nT)) \frac{ 2 K^{*}}{ C_{\zeta}} \right] \exp \left[ - nr \frac{ (\delta-\eta)^2 K^{*2}}{4 C_{\zeta}} \right] +
2 Q^2 T \exp\left(-\frac{r \zeta^4 (\log(nT))^2 }{4 Q^4 T^2 v_{n}^2 n}\right) (2nr)^{2nr+1}\\
&+
2 Q^2 T \exp \left(-\frac{r \zeta^2 (\log(nT))^2}{8^2 T^2 Q^4 v_{n}^2 n} \right) (2nr)^{2nr+1} \\
\leq& \exp \left[ -(\delta-\eta)^2 C_1 nr + C_2 \log(nT)r + C_4 \log(1/(\epsilon y_{n})) \right] +
\kappa \exp\left[ 5 nr \log(nT) - C_3 \frac{(\log(nT))^2 r}{n v_{n}^2}\right] 
\end{align*}
Then, introducing
\begin{align*}
u_{nT}&=\exp \left[ -(\delta-\eta)^2 C_1 n + C_2 \log(nT) + C_4 \log(1/(\epsilon y_{n})) \right]\\
w_{nT}&= \exp\left[ - C_3 \frac{(\log(nT))^2}{n v_{n}^2} + 5 n \log(nT)\right],
\end{align*}
we conclude as in the proof of Theorem~\ref{prop:ratio}, noticing that $nT u_{nT}$ (resp. $nT w_{nT}$) converges to 0 as $n$ increases as long as $\log(1/y_{n})=o(n)$  
(resp. as long as $v_{n}=o(\sqrt{\log(n)}/n)$).
\qed
\subsection{Proof of Corollary~\ref{cor:vitesse_cv_pi_fixedT_variationnel}}
As in the proof of Theorem~\ref{prop:vitesse_cv_M_variationnel}, using the convergence in Equation~\eqref{eq:cv_M_fixed_T} and Lemma~\ref{lem:variationel_asymp_eq}, we obtain for any $\epsilon>0$
\begin{align*}
\proba_{\theta^*}\left(\sup_{\theta \in \Theta}\left| \frac{2}{n(n-1)T} \mathcal{J}(\hat{\chi}(\theta),\theta) - \mathbb{M}^T(\pi^{1:T}) \right| >
\frac{\epsilon r_n^2}{\sqrt{n}}\right) \mathop{\longrightarrow}_{n \to +\infty} 0.
\end{align*}
We then conclude by using Lemma~\ref{lem:cv_M_to_cv_pi_bis} applied with $F_{n,T}=\frac{2}{n(n-1)T} \mathcal{J}(\hat{\chi}(\cdot),\cdot)$.
\qed
\section{Proofs of technical lemmas}
\label{appendix:proofs_tech}
\subsection{Proof of Lemma~\ref{estgamma}}
As in the proof of Lemma E.2 from \cite{celisse2012consistency}, we use the method of Lagrange multipliers to find
  the fixed-point equation of the critical point. Recall that $\theta=(\Gamma,\pi)$ and let us denote the likelihood
$L(\Gamma, \pi) \coloneqq \exp \ell(\theta)=\proba_{\theta}(X^{1:T})$
and the conditional likelihood $L_c(z^{1:T},\pi)=\proba_{\theta}(X^{1:T} \given Z^{1:T}=z^{1:T})$. 
Recall the definition of $N_{ql}(z^{1:T})$ in~\eqref{eq:def_nql} and that 
\[
  \proba_{\theta}(Z^{1:T}= z^{1:T})= \prod_{1\leq q,l \leq Q}  \gamma_{ql}^{N_{ql}(z^{1:T})} \prod_{i=1}^n \alpha^1_{z_i^1} .
\]
We compute the derivative of the Lagrangian  with respect to each parameter $\gamma_{ql}$.

\begin{align*}
\frac \partial {\partial \gamma_{ql}} \left[ L(\Gamma, \pi) + \sum_{m=1}^Q \lambda_m \left( \sum_{k=1}^Q \gamma_{mk} -1 \right) \right]
  &= \frac\partial {\partial \gamma_{ql}}\left( \sum_{z^{1:T}} L_c(z^{1:T},\pi) \proba_{\theta}(Z^{1:T}= z^{1:T})   \right) + \lambda_q \\
  &= \sum_{z^{1:T}} L_c(z^{1:T},\pi) \frac{N_{ql}(z^{1:T})}{\gamma_{ql}} \proba_{\theta}(Z^{1:T}= z^{1:T}) + \lambda_q   \\
  &= \frac 1 {\gamma_{ql}}\left(  \sum_{t=1}^{T-1} \sum_{i=1}^n\sum_{z^{1:T}} L_c(z^{1:T},\pi)
    \proba_{\theta}(Z^{1:T}= z^{1:T}) \ind_{z_i^{t}=q,z_i^{t+1}=l}  + \lambda_q \gamma_{ql}\right)\\
 &= \frac 1 {\gamma_{ql}}\left(  \sum_{t=1}^{T-1} \sum_{i=1}^n \probat(X^{1:T}, Z_i^{t}=q,Z_i^{t+1}=l)  + \lambda_q \gamma_{ql}\right). 
\end{align*}
At the critical point $\breve \theta=(\breve \gamma, \breve \pi)$, we obtain that for each $(q,l)\in \Q^2$ we have
\[
\breve \gamma_{ql} \propto \sum_{t=1}^{T-1} \sum_{i=1}^n \proba_{\breve \theta}(X^{1:T}, Z_i^{t}=q,Z_i^{t+1}=l)
  \]
where $\propto$ means 'proportional to'. The constraint $\sum_l \gamma_{ql}=1$ gives the normalizing term and we obtain 
\begin{align*}
\breve{\gamma}_{ql} = \frac{\sum_{t=1}^{T-1} \sum_{i=1}^n \proba_{\breve{\theta}}(X^{1:T}, Z_i^{t}=q, Z_i^{t+1}=l)}{\sum_{t=1}^{T-1} \sum_{i=1}^n \proba_{\breve{\theta}}(X^{1:T}, Z_i^{t}=q)} = \frac{\sum_{t=1}^{T-1} \sum_{i=1}^n \proba_{\breve{\theta}}(Z_i^{t}=q, Z_i^{t+1}=l \given X^{1:T})}{\sum_{t=1}^{T-1} \sum_{i=1}^n \proba_{\breve{\theta}}(Z_i^{t}=q \given X^{1:T})} .
\end{align*}
\qed

\subsection{Proof of Lemma~\ref{lem:fixed_point_variational}}
We can write the quantity to optimize
\begin{align} \label{eq:variational}
\mathcal{J}(\chi,\theta)=&\esp_{\mathbb{Q}_{\chi}} \left[\log \proba_{\theta}(X^{1:T},Z^{1:T})\right] + \mathcal{H}(\mathbb{Q}_{\chi})\nonumber \\
=& \esp_{\mathbb{Q}_{\chi}} \left[\log \proba_{\theta}(X^{1:T}\given Z^{1:T})\right] + \esp_{\mathbb{Q}_{\chi}} \left[\log \proba_{\theta}(Z^{1:T})\right] - \esp_{\mathbb{Q}_{\chi}} \left[\log \mathbb{Q}_{\chi}(Z^{1:T})\right]\nonumber \\
=&\esp_{\mathbb{Q}_{\chi}} \left[\sum_{t=1}^T \sum_{i< j} X^t_{ij} \log \pi_{Z_i^t Z_j^t} + (1-X^t_{ij}) \log(1-\pi_{Z_i^t Z_j^t}) \right]\nonumber \\ &+ \esp_{\mathbb{Q}_{\chi}} \left[\sum_{i=1}^n \log \alpha_{Z^1_i} + \sum_{i=1}^n \sum_{t=1}^{T-1} \log \gamma_{Z^{t}_i Z^{t+1}_i} \right] - \esp_{\mathbb{Q}_{\chi}} \left[\sum_{i=1}^n \log \mathbb{Q}_{\chi}(Z^1_i) + \sum_{i=1}^n \sum_{t=1}^{T-1} \log \mathbb{Q}_{\chi}(Z^{t+1}_i \given Z^{t}_i)\right]\nonumber \\
=&\sum_{t=1}^T \sum_{i< j} \sum_{q,l} \tau^t_{iq} \tau^t_{jl} \left[ X^t_{ij} \log \pi_{ql} + (1-X^t_{ij}) \log(1-\pi_{ql})\right]\nonumber \\ &+ \sum_{i=1}^n \sum_{q=1}^Q \tau^1_{iq} \log \alpha_{q} + \sum_{i=1}^n \sum_{q,l} \sum_{t=1}^{T-1} \eta^{t}_{iql} \log \gamma_{ql} - \sum_{i=1}^n \sum_{q=1}^Q \tau_{iq}^1 \log \tau^1_{iq} - \sum_{i=1}^n \sum_{t=1}^{T-1} \sum_{q,l} \eta^t_{iql} \log \frac{\eta^t_{iql}}{\tau^t_{iq}}.
\end{align}
Using this expression, we can obtain directly the expected fixed-point equation for the variational estimator of the transition probability from $q$ to $l$.
\qed

\subsection{Proof of Lemma~\ref{lemmet1}}
  We rely on the notation introduced in the proof of Theorem~\ref{prop:vitesse_cv_M}. 
For any $t \in \T$, using classical dependency rules in directed acyclic graphs and the expression~\eqref{eq:zthat} of $\hat{z}^t$, we  write  
\begin{align*}
\log \probat(X^{t} \given X^{1:t-1}) 
&= \log \sum_{z^t}\probat(X^{t} \given Z^{t}=z^{t}) \probat(Z^{t}=z^{t}\given X^{1:t-1}) \\
&\leq \log \left[\probat(X^{t} \given {Z}^{t}=\hat{z}^{t}) \sum_{z^t} \probat(Z^{t}=z^{t} \given X^{1:t-1}) \right] = \log \probat(X^{t} \given {Z}^{t}=\hat{z}^{t})
		\end{align*}
and thus
	\begin{equation}
		\log \probat(X^{t} \given X^{1:t-1}) - \log \probat(X^{t} \given {Z}^{t}=\hat{z}^{t}) \leq 0 \label{ineg1}.
		\end{equation}
Using Bayes' rule, we have  
\[
\log \probat(X^{t} \given X^{1:t-1}) = \log \probat(X^{t}, Z^{t} \given X^{1:t-1}) - \log \probat(Z^{t} \given X^{1:t}).
\]
Taking the expectation of this quantity with respect to any distribution $\mathbb{Q}$ on $Z^t$, we obtain 
\begin{align*}
   \log \probat(X^{t} \given X^{1:t-1})
		&= \esp_{\mathbb{Q}} \left[ \log \probat(X^{t}, Z^{t} \given X^{1:t-1}) \right] + \mathrm{KL}\left(\mathbb{Q} ;\probat(Z^{t} \given X^{1:t}) \right) + \mathcal{H}(\mathbb{Q})\\
		&\geq \esp_{\mathbb{Q}} \left[ \log \probat(X^{t}, Z^{t} \given X^{1:t-1}) \right] + \mathcal{H}(\mathbb{Q})\\
		&\geq  \esp_{\mathbb{Q}} \left[ \log \probat(X^{t} \given Z^{t}) \right] + \esp_{\mathbb{Q}} \left[ \log \probat(Z^{t} \given X^{1:t-1}) \right] +\mathcal{H}(\mathbb{Q}),
		\end{align*}
where $\mathrm{KL}\left(\mathbb{Q} ;\probat(Z^{t} \given X^{1:t}) \right) = \esp_{\mathbb{Q}} \left[ \log \mathbb{Q}(Z^t) - \log \probat(Z^{t} \given X^{1:t}) \right]$ is a  Kullback-Leibler divergence (thus non negative) and  $\mathcal{H}(\mathbb{Q}) = - \esp_{\mathbb{Q}} \left[ \log \mathbb{Q}(Z^t) \right]$ is the entropy of $\mathbb{Q}$.

Taking now $\mathbb{Q}$ as the Dirac distribution located on $\hat{z}^t$ , we have $\mathcal{H}(\mathbb{Q})=0$ and 
		\begin{equation}
		\log \probat(X^{t} \given X^{1:t-1}) \geq \log \probat(X^{t} \given {Z}^{t}=\hat{z}^{t}) + \log \probat({Z}^{t}=\hat{z}^{t} \given X^{1:t-1}) \label{ineg2}.
		\end{equation}
Now, combining Inequalities~\eqref{ineg1} and \eqref{ineg2}, we obtain
\[
\log \probat({Z}^{t}=\hat{z}^{t} \given X^{1:t-1}) \leq \log \probat(X^{t} \given X^{1:t-1}) - \log \probat(X^{t}
   \given {Z}^{t}=\hat{z}^{t}) \leq 0,
\]
giving the expected result.
\qed

\subsection{Proof of Lemma~\ref{lem:app_talagrand}}
To prove this lemma, we first establish a control of the expectation of the random variable appearing in the statement. 

 \begin{lemme} \label{lemmebound}
We have the following inequality for $z^{*1:T}$ and $z^{1:T}$ any configurations and any $\theta \in \Theta$ 
\[
  \espetoile \left[ \sup_{(z^{1:T},\pi) \in \Q^{nT}\times [\zeta,1-\zeta]^{Q^2}} \left|\frac{2}{n(n-1)T} \sum_{t=1}^T  \sum_{i<j}
      (\Xtij-\pi^{*}_{z^{*t}_i z^{*t}_j}) \log\left(\frac{\pi_{z^t_i z^t_j}}{1-\pi_{z^t_i z^t_j}}\right) \right|
    \given[\Big] Z^{1:T}=z^{*1:T} \right] \leq \sqrt{\frac{2}{n(n-1)T}}\Lambda 
\]
with $\Lambda=2 \log[(1-\zeta)/\zeta]$.
\end{lemme} 
We now turn to the proof of Lemma~\ref{lem:app_talagrand}. 
Let us first recall Talagrand's inequality~\cite[see for e.g.][page 170, Equation (5.50)]{massart2007concentration}. 

\begin{theorem*}[Talagrand's inequality]
Let $\{Y_{ij}^t\}_{1 \leq i<j \leq n, 1\leq t \leq T}$ denote independent and centered random variables. Define 
\[
\forall g\coloneqq \{g_{ij}^t\}_{1\leq i<j \leq n, 1\leq t \leq T} \in \mathcal{G}, \quad S_{n,T}(g)=\sum_{1\le i<j\le n } \sum_{t=1}^T Y^t_{ij} g^t_{ij} ,
\]
where $\mathcal{G}\subset\mathbb{R}^{n(n-1)T/2}$. 
Let us further assume that there exist $b>0$ and $\sigma^2>0$ such that $|Y^t_{ij} g^t_{ij}|\leq b$ for every
$(i,j,t)\in \n^2\times \T$ and any $g\in \mathcal{G}$ and $\sup_{g\in \mathcal{G}} \sum_{i<j}\sum_t \var(Y^t_{ij}g^t_{ij}) \leq \sigma^2$.
Then, for every $\beta > 0$ and $x>0$, for any finite set  $\{g_1,\ldots,g_{2^{n(n-1)T/2}}\}$ of elements of
  $\mathcal{G}$, we have 
\begin{equation}
\proba\left( \max_{g \in \{g_1,\ldots,g_{2^{n(n-1)T/2}}\}}  S_{n,T}(g) \geq \esp\left[\max_{g \in \{g_1,\ldots,g_{2^{n(n-1)T/2}}\}} S_{n,T}(g)\right] (1+\beta) + \sqrt{2 \sigma^2 x} + b(\beta^{-1}+3^{-1})x \right) \leq e^{-x}.
\end{equation}
\end{theorem*}

First, notice that $\argmin_{\varpi \in [\zeta,1-\zeta]} \log(\varpi/(1-\varpi))=\zeta$ and $\argmax_{\varpi \in
  [\zeta,1-\zeta]} \log(\varpi/(1-\varpi))=1-\zeta$ so that we have
\begin{multline*}
\probaetoile\left( \sup_{(z^{1:T},\pi) \in \Q^{nT}\times [\zeta,1-\zeta]^{Q^2}} \frac{2}{n(n-1)T} \left| \sum_{t=1}^T  \sum_{i<j} (\Xtij-\pi^{*}_{z^{*t}_i z^{*t}_j}) \log\left(\frac{\pi_{z^t_i z^t_j}}{1-\pi_{z^t_i z^t_j}}\right) \right| > \epsilon \right) \\ \leq \probaetoile\left( \max_{\varpi \in \{\zeta,1-\zeta\}^{n(n-1)T/2}} \frac{2}{n(n-1)T} \left| \sum_{t=1}^T  \sum_{i<j} (\Xtij-\pi^{*}_{z^{*t}_i z^{*t}_j}) \log\left(\frac{\varpi_{i,j}^t}{1-\varpi_{i,j}^t}\right) \right| > \epsilon \right)
\end{multline*}
with $\varpi\coloneqq\{\varpi_{i,j}^t\}_{1\leq i<j \leq n, 1 \leq t \leq T}$. 
The set $\{\zeta,1-\zeta\}^{n(n-1)T/2}$ is finite, of size ${2^{n(n-1)T/2}}$.
Let us now apply Talagrand's inequality to our setup. 
Note that for every $(i,j,t) \in \n^2\times \T$, for any $\pi \in [\zeta,1-\zeta]^{Q^2}$, we have
\[
\left|(\Xtij-\pi^{*}_{z^{*t}_i z^{*t}_j}) \log\left(\frac{\pi_{z^t_i z^t_j}}{1-\pi_{z^t_i z^t_j}}\right) \right| \leq \log[(1-\zeta)/\zeta]=\frac{\Lambda}{2} 
\] 
almost surely thanks to Assumption~\ref{itm:hyp_pi}, and with $\Lambda$ as defined in Lemma~\ref{lemmebound}. 
Combining this result with Lemma~\ref{lemmebound} and writing $\Omega=(1+\beta) \Lambda \sqrt{n(n-1)T/2}  + \sqrt{n(n-1) T
  (\Lambda/2)^2 x_{n,T}} + (1/\beta+1/3) (\Lambda/2) x_{n,T} $, we have for any $\epsilon>0$, for any $\beta>0$, applying Talagrand's inequality with $b=\Lambda/2$ and $\sigma^2=n(n-1)T/2(\Lambda/2)^2$,
	\begin{align*}
	& \probaetoile \left( \sup_{(z^{1:T},\pi) \in \Q^{nT}\times [\zeta,1-\zeta]^{Q^2}} \frac{2}{n(n-1)T} \left| \sum_{t=1}^T  \sum_{i<j} (\Xtij-\pi^{*}_{z^{*t}_i z^{*t}_j}) \log\left(\frac{\pi_{z^t_i z^t_j}}{1-\pi_{z^t_i z^t_j}}\right) \right| > \epsilon \right)  \\
	\leq & \probaetoile \left( \max_{\varpi \in \{\zeta,1-\zeta\}^{n(n-1)T/2}} \frac{2}{n(n-1)T} \left| \sum_{t=1}^T  \sum_{i<j} (\Xtij-\pi^{*}_{z^{*t}_i z^{*t}_j}) \log\left(\frac{\varpi_{i,j}^t}{1-\varpi_{i,j}^t}\right) \right| > \epsilon \right) \\
	\leq & \probaetoile\left( \epsilon < \max_{\varpi \in \{\zeta,1-\zeta\}^{n(n-1)T/2}} \frac{2}{n(n-1)T} \left| \sum_{t=1}^T  \sum_{i<j} (\Xtij-\pi^{*}_{z^{*t}_i z^{*t}_j}) \log\left(\frac{\varpi_{i,j}^t}{1-\varpi_{i,j}^t}\right) \right| \leq \frac{2}{n(n-1)T} \Omega \right) \\
	&+ \probaetoile\left(\max_{\varpi \in \{\zeta,1-\zeta\}^{n(n-1)T/2}} \frac{2}{n(n-1)T} \left| \sum_{t=1}^T  \sum_{i<j} (\Xtij-\pi^{*}_{z^{*t}_i z^{*t}_j}) \log\left(\frac{\varpi_{i,j}^t}{1-\varpi_{i,j}^t}\right) \right| > \frac{2}{n(n-1)T} \Omega \right)\\
	\leq&  \probaetoile \left(\frac{2}{n(n-1)T} \Omega > \epsilon\right) + 2 e^{-x_{n,T}}
	\leq \ind_{\epsilon< 2\Omega/(n(n-1)T)} + 2 e^{-x_{n,T}}.\\
	\end{align*}
\qed

\subsection{Proof of Lemma~\ref{prop:whp}}
For any $\eta \in (0,\delta)$, Hoeffding's inequality~\cite[see for example Theorem 2.8 from][]{boucheron2013concentration} gives that
  \begin{align*}
          \probat \left( \forall t \in \T, \forall q \in \Q, \frac{N_q(Z^{t})}{n} \geq \alpha_q-\eta \right)
	&= 1-\probat \left( \exists t \in \T, \exists q \in \Q; \frac{1}{n} \sum_{i=1}^n \ind_{Z^t_i=q} < \alpha_q-\eta  \right)\\
	&\geq 1-\sum_{q=1}^Q \sum_{t=1}^T \exp\left(-2 \eta^2 n\right)
	\geq 1- QT\exp\left(-2 \eta^2 n\right),
	\end{align*}
  which concludes the proof.
  \qed

\subsection{Proof of Lemma~\ref{lem:cv_matriceA}}
First notice that $\argmax_{A \in \mathcal{A}} \mathbb{M}(\pi,A)$ may not be unique, it is in fact a closed subset of $\mathcal{A}$. However, we choose a fixed element $\bar{A}_{\pi}$ in this subset in the following.
Letting $\epsilon>0$ and $\eta \in (0,\delta)$ and using Lemma~\ref{prop:whp}, we can split the probability as
\begin{align*}
\proba_{\theta^*} \left(\frac{1}{T} \sum_{t=1}^T \sup_{\pi \in [\zeta,1-\zeta]^{Q^2}} \left|\mathbb{M}(\pi,\bar{A}_{\pi}^t)  -
\mathbb{M}(\pi,\bar{A}_{\pi}) \right| > \frac{\epsilon r_n}{6 \sqrt{n}} \right) \leq& \proba_{\theta^*} \left(\left\{\frac{1}{T} \sum_{t=1}^T \sup_{\pi \in [\zeta,1-\zeta]^{Q^2}} \left|\mathbb{M}(\pi,\bar{A}_{\pi}^t)  -
\mathbb{M}(\pi,\bar{A}_{\pi}) \right| > \frac{\epsilon r_n}{6 \sqrt{n}} \right\} \cap \Omega_{\eta}(\theta^*) \right) \\ &+ QT \exp\left(-2 \eta^2 n \right),
\end{align*}
recalling that
\[
\Omega_{\eta}(\theta)\coloneqq \left\{z^{1:T} \in \Q^{nT} ; \forall t \in \T, \forall q \in \Q,
\frac{N_q(z^{t})}{n} \geq \alpha_q-\eta \right\}.
\]
We thus want to bound the quantity $\proba_{\theta^*} \left(T^{-1} \sum_{t=1}^T \sup_{\pi \in [\zeta,1-\zeta]^{Q^2}}
  \left|\mathbb{M}(\pi,\bar{A}_{\pi}^t)  -	\mathbb{M}(\pi,\bar{A}_{\pi}) \right| > \epsilon r_n/(6 \sqrt{n}) \right)$ on the event $\left\{Z^{1:T}
  \in \Omega_{\eta}(\theta^*)\right\}$, which means bounding  
\[
\proba_{\theta^*} \left(\frac{1}{T} \sum_{t=1}^T \sup_{\pi \in [\zeta,1-\zeta]^{Q^2}} \left|\mathbb{M}(\pi,\bar{A}_{\pi}^t)  -
	\mathbb{M}(\pi,\bar{A}_{\pi}) \right| > \frac{\epsilon r_n}{6 \sqrt{n}} \given[\Bigg] Z^{1:T} \in \Omega_{\eta}(\theta^*) \right).
\]
Let us denote for any matrix $P$ of size $m \times n$ the norm $\|P\|_{\infty}=\max_{(i,j) \in \llbracket 1,m \rrbracket \times \llbracket 1,n \rrbracket} |P_{ij}|$. Then note that, for any matrix $\breve{A}$ with coefficients in $[0,1]$, for any $\pi \in [\zeta,1-\zeta]^{Q^2}$, using Assumption~\ref{itm:hyp_gamma} and \ref{itm:hyp_pi}, 
\begin{align*}
\left( \mathbb{M}(\pi,\bar{A}_{\pi}) - \mathbb{M}(\pi,\breve{A})\right)  &\leq \sum_{q,l} \alpha^*_q \alpha^*_l \sum_{q',l'} |\bar{a}_{qq'}\bar{a}_{ll'}-\breve{a}_{qq'}\breve{a}_{ll'}| \sup_{\pi \in [\zeta,1-\zeta]^{Q^2}}  | \pi^*_{ql} \log \pi_{q'l'} + (1-\pi^*_{ql}) \log (1-\pi_{q'l'})|\\
&\leq 2 (1-\delta)^2 (1-\zeta) \log(1/\zeta) \sum_{q,l}\sum_{q',l'} |\bar{a}_{qq'}\bar{a}_{ll'}-\breve{a}_{qq'}\breve{a}_{ll'}|\\
&\leq 2 (1-\delta)^2 (1-\zeta) \log(1/\zeta) Q^4 2 \| \breve{A}-\bar{A}_{\pi} \|_{\infty} \coloneqq c \| \breve{A}-\bar{A}_{\pi} \|_{\infty}
\end{align*}
with $c=4 (1-\delta)^2 (1-\zeta) \log(1/\zeta) Q^4$. On the event $\Omega_{\eta}(\theta^*)$ we then have
	\begin{align*}
	&\proba_{\theta^*} \left( \frac{1}{T} \sum_{t=1}^T \sup_{\pi \in [\zeta,1-\zeta]^{Q^2}}  \left|\mathbb{M}(\pi,\bar{A}_{\pi}^t)  -
	\mathbb{M}(\pi,\bar{A}_{\pi}) \right| > \frac{\epsilon r_n}{6 \sqrt{n}} \right) \\ =&  1 - \proba_{\theta^*} \left(\frac{1}{T} \sum_{t=1}^T \sup_{\pi \in [\zeta,1-\zeta]^{Q^2}}  \left|\mathbb{M}(\pi,\bar{A}_{\pi}^t)  -
	\mathbb{M}(\pi,\bar{A}_{\pi}) \right| \leq \frac{\epsilon r_n}{6 \sqrt{n}} \right) \\
	\leq& 1- \proba_{\theta^*} \left(\forall t \in \T, \sup_{\pi \in [\zeta,1-\zeta]^{Q^2}}  \left(\mathbb{M}(\pi,\bar{A}_{\pi})  - \mathbb{M}(\pi,\bar{A}_{\pi}^t)\right) \leq \frac{\epsilon r_n}{6 \sqrt{n}} \right)\\
	\leq& 1- \proba_{\theta^*} \left(\forall t \in \T, \forall \pi \in [\zeta,1-\zeta]^{Q^2},  \left(\mathbb{M}(\pi,\bar{A}_{\pi})  - \mathbb{M}(\pi,\bar{A}_{\pi}^t)\right) \leq \frac{\epsilon r_n}{6 \sqrt{n}} \right)\\
	\leq& 1-\proba_{\theta^*} \left(\forall t \in \T, \forall \pi \in [\zeta,1-\zeta]^{Q^2}, \exists \breve{A} \in \mathcal{A}^t(Z^{1:T});  \left(\mathbb{M}(\pi,\bar{A}_{\pi})  - \mathbb{M}(\pi,\breve{A})\right) \leq \frac{\epsilon r_n}{6 \sqrt{n}} \right)\\
	\leq& 1-\proba_{\theta^*} \left(\forall t \in \T, \forall \pi \in [\zeta,1-\zeta]^{Q^2}, \exists \breve{A} \in \mathcal{A}^t(Z^{1:T}); \| \breve{A}-\bar{A}_{\pi} \|_{\infty}< \frac{\epsilon r_n}{6c \sqrt{n}}\right).
	\end{align*}
	We then show that for any $\epsilon>0$, for every $t \in \T$ and every $\pi \in [\zeta,1-\zeta]^{Q^2}$, for any $n$ such that $n > 6c\sqrt{n} / [\epsilon r_n  (\delta - \eta)]$, there exists some $\breve{A}\in \mathcal{A}^t(Z^{1:T})$ such that $\| \breve{A}-\bar{A}_{\pi} \|_{\infty}< \epsilon r_n /(6c \sqrt{n})$, i.e. such that for every $q,l$, $| \breve{a}_{ql} - \bar{a}_{ql} | < \epsilon r_n /(6c \sqrt{n})$. For every $1\leq q \leq Q$, we can construct $\breve{A}_{q\cdot}=(\breve{a}_{q1},\ldots,\breve{a}_{qQ})$ as follows. 
	On the event $\Omega_{\eta}(\theta^*)$, for every $q \in \Q$, for any $n$ such that $n > 6c\sqrt{n} / [\epsilon r_n  (\delta - \eta)]$, we have $N_q(Z^t) \epsilon r_n/(6c\sqrt{n})> 1$ for every $t \in \T$.
		We then construct $(\breve{n}_{ql})_{1\leq l\leq Q}$ as follows and take $\breve{a}_{ql}=\breve{n}_{ql}/N_q(Z^{1:T})$ for every $l \in \Q$.
		\begin{itemize}
			\item for $l=1$ choose $\breve{n}_{q1}$ as the closest integer to $N_q(Z^t) \bar{a}_{q1}$. 
			It is in the interval $(N_q(Z^t) \bar{a}_{q1} - 1, N_q(Z^t) \bar{a}_{q1} + 1)$ so we have $|\bar{a}_{q1} - \breve{n}_{q1}/N_q(Z^t)| < 1/N_q(Z^{t})< \epsilon r_n/(6c\sqrt{n})$. Moreover, note that $0 \leq \breve{n}_{q1} \leq N_q(Z^t)$ because $0 \leq N_q(Z^t) \bar{a}_{q1} \leq N_q(Z^t)$.
			\item Repeat for $l=2,\ldots,Q$ \begin{itemize}
				\item if $\sum_{l'=1}^{l-1} (N_q(Z^t) \bar{a}_{ql'}-\breve{n}_{ql'}) \geq 0$ choose $\breve{n}_{ql}$ as the closest bigger (or equal) integer to $N_q(Z^t) \bar{a}_{ql}$.
				\item if $\sum_{l'=1}^{l-1} (N_q(Z^t) \bar{a}_{ql'}-\breve{n}_{ql'}) < 0$ choose $\breve{n}_{ql}$ as the closest smaller (or equal) integer to  $N_q(Z^t) \bar{a}_{ql}$.
			\end{itemize}
				As before, $\breve{n}_{ql}$ is in the interval $(N_q(Z^t) \bar{a}_{ql} - 1, N_q(Z^t) \bar{a}_{ql} + 1)$ so we have $|\bar{a}_{ql} - \breve{n}_{ql}/N_q(Z^t)| < 1/N_q(Z^{1:T}) < \epsilon r_n/(6c \sqrt{n})$. Moreover $0 \leq \breve{n}_{ql} \leq N_q(Z^t)$ because $0 \leq N_q(Z^t) \bar{a}_{ql} \leq N_q(Z^t)$. We also have (by induction) 
				\[
				\left|\sum_{l'=1}^{l} (N_q(Z^t) \bar{a}_{ql'}-\breve{n}_{ql'})\right| = \left|\left(\sum_{l'=1}^{l-1} N_q(Z^t) \bar{a}_{ql'}-\breve{n}_{ql'}\right) + N_q(Z^t) \bar{a}_{ql}-\breve{n}_{ql} \right| < 1.
				\]
		\end{itemize}
	In the end, we have $|\sum_{l=1}^{Q} (N_q(Z^t) \bar{a}_{ql}-\breve{n}_{ql})| < 1$ i.e. $|N_q(Z^t) - \sum_{l=1}^{Q}\breve{n}_{ql}| < 1$, meaning that $\sum_{l=1}^{Q}\breve{n}_{ql}=N_q(Z^t)$, both $N_q(Z^t)$ and $\sum_{l=1}^{Q}\breve{n}_{ql}$ being integers. Then, if $n > 6c\sqrt{n} / [\epsilon r_n (\delta-\eta)]$, there exists $\breve{A} \in \mathcal{A}^t(Z^{1:T})$ such that $\| \breve{A}-\bar{A}_{\pi} \|_{\infty}<\epsilon r_n/(6c\sqrt{n})$. This leads to
	\begin{align*}
	\proba_{\theta^*} \left( \frac{1}{T} \sum_{t=1}^T \left|\mathbb{M}(\pi,\bar{A}_{\pi}^t) -
	\mathbb{M}(\pi,\bar{A}_{\pi}) \right| > \frac{\epsilon r_n}{6\sqrt{n}} \right)
	&\leq QT \exp(-2\eta^2n) + 1- \ind_{n > 6c \sqrt{n} / [\epsilon r_n (\delta-\eta)]}
	\end{align*}
	which concludes the proof.
\qed

\subsection{Proof of Lemma~\ref{lem:thirdterm}}
We can upper bound the expectation as follows
\begin{align*}
\esp_{\theta^*}\left[\left| \frac{N_{q}(Z^1) N_{l}(Z^1)}{n(n-1)} - \alpha^*_q \alpha^*_l \right| \right] =& \esp_{\theta^*}\left[\left| \left(\frac{N_{q}(Z^1)}{n} - \alpha^*_q\right) \frac{N_{l}(Z^1)}{n-1} + \alpha^*_q \left(\frac{N_{l}(Z^1)}{n-1} - \alpha^*_l\right) \right| \right]\\
\leq& \esp_{\theta^*}\left[\left| \frac{N_{q}(Z^1)}{n} - \alpha^*_q\right| \frac{N_{l}(Z^1)}{n-1} \right] +  \alpha^*_q \esp_{\theta^*}\left[ \left|\frac{N_{l}(Z^1)}{n-1} - \alpha^*_l\right| \right]\\
\leq& \sqrt{\esp_{\theta^*}\left[\left( \frac{N_{q}(Z^1)}{n} - \alpha^*_q\right)^2\right] \esp_{\theta^*}\left[ \frac{N_{l}(Z^1)^2}{(n-1)^2} \right]} + \alpha^*_q \sqrt{\esp_{\theta^*}\left[ \left(\frac{N_{l}(Z^1)}{n-1} - \alpha^*_l\right)^2 \right]}.
\end{align*}
   We have for any $q \in \Q$
   \begin{align*}
   \espetoile\left[N_{q}(Z^1)^2\right]=\sum_{i,j} \espetoile\left[\ind_{Z^1_i=q}\ind_{Z^1_j=q}\right]=\sum_{i} \alpha^*_q + \sum_{i\neq j} \alpha^{*2}_q= n \alpha^*_q + n(n-1) \alpha^{*2}_q.
   \end{align*}
   This implies that
   \begin{align*}
    \esp_{\theta^*}\left[ \left(\frac{N_{q}(Z^1)}{n} - \alpha^*_q\right)^2 \right]= \esp_{\theta^*}\left[ \frac{N_{q}(Z^1)^2}{n^2}\right] - \alpha^{*2}_q = \frac{1}{n} \alpha^*_q + \frac{n-1}{n} \alpha^{*2}_q - \alpha^{*2}_q = \frac{1}{n} \alpha^*_q (1-\alpha^*_q),
   \end{align*}
   and identically
   \begin{align*}
    \esp_{\theta^*}\left[ \left(\frac{N_{l}(Z^1)}{n-1} - \alpha^*_l\right)^2 \right]&= \esp_{\theta^*}\left[ \frac{N_{l}(Z^1)^2}{(n-1)^2}\right] + \alpha^{*2}_l - 2 \frac{n}{n-1}\alpha^{*2}_l 
    =  \frac{n}{(n-1)^2} \alpha^*_l - \frac{1}{n-1} \alpha^{*2}_l =\frac{1}{n-1} \alpha^*_l \left(\frac{n}{n-1}-\alpha^*_l\right).
   \end{align*}
  This leads to
   \begin{align}
   \esp_{\theta^*}\left[\left| \frac{N_{q}(Z^1) N_{l}(Z^1)}{n(n-1)} - \alpha^*_q \alpha^*_l \right| \right] 
   \leq& \sqrt{ \frac{1}{n} \alpha^*_q (1-\alpha^*_q)\left( \frac{n}{(n-1)^2}\alpha^*_q +\frac{n}{n-1}\alpha^{*2}_q \right)} + \alpha^*_q \sqrt{\frac{1}{n-1} \alpha^*_l \left(\frac{n}{n-1}-\alpha^*_l\right)} \nonumber\\
   \leq& \sqrt{ \frac{1}{(n-1)^2} +\frac{1}{n-1}} + \sqrt{\frac{n}{(n-1)^2}} \leq 2 \frac{\sqrt{n}}{n-1}\label{eq:maj_esp_thirdterm},
   \end{align}
   using the fact that $0 \leq \alpha_q^* \leq 1$ for every $q \in \Q$.
\qed

\subsection{Proof of Lemma~\ref{lem:cv_M_to_cv_pi_bis}}
We first consider the case when $T \rightarrow \infty$, and $\pi$ is constant over time. 
We use the following lemma.
	\begin{lemme}\label{lem:fct_M}
		For any $\theta \in \Theta$, we have for $\epsilon$ small enough ($0<\epsilon <\min_{1\leq q\neq q' \leq Q} \max_{1\leq l\leq Q} |\pi^*_{ql}-\pi^*_{q'l}|/2$)
		\[
		\min_{\sigma \in \mathfrak{S}_Q}\|\pi_{\sigma}-\pi^*\|_{\infty} > \epsilon \implies \mathbb{M}(\pi^*)-\mathbb{M}(\pi) > \frac{2\delta^2}{Q^2} \epsilon^2.
		\]
	\end{lemme}
	This gives an upper bound on the probability of interest
	\[
	\proba_{\theta^*}\left( \min_{\sigma \in \mathfrak{S}_Q}\|\hat{\pi}_{\sigma}-\pi^*\|_{\infty}>\epsilon \sqrt{v_{n,T}} \right) \leq \proba_{\theta^*}\left(\mathbb{M}(\pi^*)-\mathbb{M}(\hat{\pi})> \frac{2\delta^2}{Q^2} \epsilon^2 v_{n,T} \right).
	\]
	By definition of $\hat{\theta}=(\hat{\Gamma},\hat{\pi})$, we write
	\[
	\mathbb{M}(\pi^*) = F_{n,T}(\hat{\Gamma},\pi^*) + \mathbb{M}(\pi^*) - F_{n,T}(\hat{\Gamma},\pi^*) \leq F_{n,T}(\hat{\Gamma},\hat{\pi}) + \mathbb{M}(\pi^*) - F_{n,T}(\hat{\Gamma},\pi^*),
	\]
	implying that
	\[
	\mathbb{M}(\pi^*) - \mathbb{M}(\hat{\pi}) \leq 	\left[ F_{n,T}(\hat{\Gamma},\hat{\pi}) - \mathbb{M}(\hat{\pi}) \right] + \left[ \mathbb{M}(\pi^*) - F_{n,T}(\hat{\Gamma},\pi^*) \right].
	\]
	We then obtain the following upper bound, that converges to $0$ as $n$ and $T$ increase by assumption,
	\begin{align*}
	\proba_{\theta^*}\left(\min_{\sigma \in \mathfrak{S}_Q} \|\hat{\pi}_{\sigma}-\pi^*\|_{\infty}> \epsilon \sqrt{v_{n,T}} \right) \leq& \proba_{\theta^*}\left(F_{n,T}(\hat{\Gamma},\hat{\pi}) - \mathbb{M}(\hat{\pi})> \frac{\delta^2}{Q^2} \epsilon^2 v_{n,T} \right)+ \proba_{\theta^*}\left(  \mathbb{M}(\pi^*) - F_{n,T}(\hat{\Gamma},\pi^*) > \frac{\delta^2}{Q^2} \epsilon^2 v_{n,T}\right). 
	\end{align*}
	When the number of time steps $T$ is fixed and $\pi$ is allowed to vary over time, the proof is almost the same. Indeed, $\min_{\sigma^{1},\ldots,\sigma^T \in \mathfrak{S}_Q} \|\hat{\pi}_{\sigma^{1:T}}^{1:T}-\pi^{*1:T}\|_{\infty}> \epsilon \sqrt{v_n}$ means that there exists $t \in \T$ such that $\min_{\sigma^t \in \mathfrak{S}_Q} \|\hat{\pi}_{\sigma^t}^t-\pi^{*t}\|_{\infty}> \epsilon \sqrt{v_n}$ and we can apply Lemma~\ref{lem:fct_M} to this $\hat{\pi}^t$ to obtain that $\mathbb{M}(\pi^{*t})-\mathbb{M}(\hat{\pi}^t) > 2 \epsilon^2 \delta^2 v_n/Q^2$.
	This implies that $\mathbb{M}^T(\pi^{*1:T})-\mathbb{M}^T(\hat{\pi}^{1:T}) > 2 \epsilon^2 \delta^2  v_n /(TQ^2)$, which allows to conclude in the same way as before.
\qed

\subsection{Proof of Lemma~\ref{lem:ratio_3terms}}
We have 
\begin{align*}
\log \frac{\proba_{\breve{\theta}}(Z^{1:T}=z^{1:T} \given X^{1:T})}{\proba_{\breve{\theta}}(Z^{1:T}=z^{*1:T} \given X^{1:T})}&=\log \frac{\proba_{\breve{\theta}}(X^{1:T} \given Z^{1:T}=z^{1:T})}{\proba_{\breve{\theta}}(X^{1:T}\given Z^{1:T}=z^{*1:T})}+ \log \frac{\proba_{\breve{\theta}}(Z^{1:T}=z^{1:T})}{\proba_{\breve{\theta}}(Z^{1:T}=z^{*1:T})}\\
&= \sum_{t=1}^T \sum_{1\leq i< j \leq n} \left( X_{ij}^t \log\frac{\breve{\pi}_{z^t_i z^t_j}}{\breve{\pi}_{z^{*t}_i z^{*t}_j}} +
(1-X_{ij}^t) \log\frac{1-\breve{\pi}_{z^t_i z^t_j}}{1-\breve{\pi}_{z^{*t}_i z^{*t}_j}} \right) + \sum_{i=1}^n \log
\frac{\breve{\alpha}_{z^1_i}}{\breve{\alpha}_{z^{*1}_i}} + \sum_{t=1}^{T-1} \sum_{i=1}^n \log \frac{\breve{\gamma}_{z^{t}_i
		z^{t+1}_i}}{\breve{\gamma}_{z^{*t}_i z^{*t+1}_i}} .
\end{align*}
We decompose this sum as
\begin{align}
\label{decompotrans}
\log \frac{\proba_{\breve{\theta}}(Z^{1:T}=z^{1:T}\given X^{1:T})}{\proba_{\breve{\theta}}(Z^{1:T}=z^{*1:T}\given X^{1:T})}=&\sum_{t=1}^T \sum_{1\leq i<j \leq n} \left( X_{ij}^t \log\frac{\pi^{*}_{z^t_i z^t_j}}{\pi^{*}_{z^{*t}_i z^{*t}_j}} + (1-X_{ij}^t) \log\frac{1-\pi^{*}_{z^t_i z^t_j}}{1-\pi^{*}_{z^{*t}_i z^{*t}_j}} \right) + \sum_{i=1}^n \log \frac{\breve{\alpha}_{z^1_i}}{\breve{\alpha}_{z^{*1}_i}} + \sum_{t=1}^{T-1} \sum_{i=1}^n \log \frac{\breve{\gamma}_{z^{t}_i z^{t+1}_i}}{\breve{\gamma}_{z^{*t}_i z^{*t+1}_i}} \nonumber\\
&+\sum_{t=1}^T \sum_{1 \leq i<j \leq n} \left( X_{ij}^t \log\frac{\breve{\pi}_{z^t_i z^t_j}}{\pi^{*}_{z^{t}_i z^{t}_j}}  \frac{\pi^{*}_{z^{*t}_i z^{*t}_j}}{\breve{\pi}_{z^{*t}_i z^{*t}_j}} + (1-X_{ij}^t) \log\frac{1-\breve{\pi}_{z^t_i z^t_j}}{1-\pi^{*}_{z^t_i z^t_j}} \frac{1-\pi^{*}_{z^{*t}_i z^{*t}_j}}{1-\breve{\pi}_{z^{*t}_i z^{*t}_j}} \right).
\end{align}
In the first sum of the right-hand side of~\eqref{decompotrans}, the terms are different from zero only for triplets
$(i,j,t)$ in $D^*$. Similarly in the last sum, the terms are different from zero for triplets $(i,j,t)$ in $D^* \cup
\breve{D}$. As a consequence, we obtain
\begin{align*}
\log \frac{\proba_{\breve{\theta}}(Z^{1:T}=z^{1:T}\given X^{1:T})}{\proba_{\breve{\theta}}(Z^{1:T}=z^{*1:T}\given X^{1:T})}=&\sum_{(i,j,t)\in D^*} \left( X_{ij}^t \log\frac{\pi^{*}_{z^t_i z^t_j}}{\pi^{*}_{z^{*t}_i z^{*t}_j}} + (1-X_{ij}^t) \log\frac{1-\pi^{*}_{z^t_i z^t_j}}{1-\pi^{*}_{z^{*t}_i z^{*t}_j}} \right) + \sum_{i=1}^n \log \frac{\breve{\alpha}_{z^1_i}}{\breve{\alpha}_{z^{*1}_i}} + \sum_{t=1}^{T-1} \sum_{i=1}^n \log \frac{\breve{\gamma}_{z^{t}_i z^{t+1}_i}}{\breve{\gamma}_{z^{*t}_i z^{*t+1}_i}} \nonumber \\
&+ \sum_{(i,j,t)\in D^*\cup \breve{D}} \left( X_{ij}^t \log\frac{\breve{\pi}_{z^t_i z^t_j}}{\pi^{*}_{z^{t}_i z^{t}_j}}  \frac{\pi^{*}_{z^{*t}_i z^{*t}_j}}{\breve{\pi}_{z^{*t}_i z^{*t}_j}} + (1-X_{ij}^t) \log\frac{1-\breve{\pi}_{z^t_i z^t_j}}{1-\pi^{*}_{z^t_i z^t_j}} \frac{1-\pi^{*}_{z^{*t}_i z^{*t}_j}}{1-\breve{\pi}_{z^{*t}_i z^{*}_j}} \right).
\end{align*}
We now write the last sum in the right-hand side as 
\begin{align*}
&\sum_{(i,j,t)\in D^*\cup \breve{D}} \left( X_{ij}^t \log\frac{\breve{\pi}_{z^t_i z^t_j}}{\pi^{*}_{z^{t}_i z^{t}_j}}  \frac{\pi^{*}_{z^{*t}_i z^{*t}_j}}{\breve{\pi}_{z^{*t}_i z^{*t}_j}} + (1-X_{ij}^t) \log\frac{1-\breve{\pi}_{z^t_i z^t_j}}{1-\pi^{*}_{z^t_i z^t_j}} \frac{1-\pi^{*}_{z^{*t}_i z^{*t}_j}}{1-\breve{\pi}_{z^{*t}_i z^{*t}_j}} \right)\\
=&\sum_{(i,j,t)\in D^*\cup \breve{D}} \left\{ X_{ij}^t \left[ \log\left( 1+\frac{\breve{\pi}_{z^t_i
		z^t_j}-\pi^{*}_{z^t_i z^t_j}}{\pi^{*}_{z^{t}_i z^{t}_j}}\right) +\log \frac{\pi^{*}_{z^{*t}_i
		z^{*t}_j}}{\breve{\pi}_{z^{*t}_i z^{*t}_j}}\right] + (1-X_{ij}^t) \left[\log \left(1-\frac{\breve{\pi}_{z_i^t           z_j^t}-\pi^{*}_{z_i^t z_j^t}}{1-\pi^{*}_{z_i^t z_j^t}}\right) +\log \frac{1-\pi^{*}_{z^{*t}_i
		z^{*t}_j}}{1-\breve{\pi}_{z^{*t}_i z^{*t}_j}}\right] \right\}.
\end{align*}
Distinguishing between the cases where $X_{ij}^t=1$ and $X_{ij}^t=0$, we obtain
\begin{align*}
&\sum_{(i,j,t)\in D^*\cup \breve{D}} \left( X_{ij}^t \log\frac{\breve{\pi}_{z^t_i z^t_j}}{\pi^{*}_{z^{t}_i z^{t}_j}}  \frac{\pi^{*}_{z^{*t}_i z^{*t}_j}}{\breve{\pi}_{z^{*t}_i z^{*t}_j}} + (1-X_{ij}^t) \log\frac{1-\breve{\pi}_{z^t_i z^t_j}}{1-\pi^{*}_{z^t_i z^t_j}} \frac{1-\pi^{*}_{z^{*t}_i z^{*t}_j}}{1-\breve{\pi}_{z^{*t}_i z^{*t}_j}} \right)\\
=&\sum_{(i,j,t)\in D^*\cup \breve{D}}\log\left[1+\frac{(\breve{\pi}_{z^t_i z^t_j}-\pi^{*}_{z^t_i z^t_j})(X_{ij}^t-\pi^{*}_{z^t_i z^t_j})}{\pi^{*}_{z^{t}_i z^{t}_j}(1-\pi^{*}_{z^t_i z^t_j})}  \right]
-\sum_{(i,j,t)\in D^*\cup \breve{D}} \log\left[1+\frac{(\breve{\pi}_{z^{*t}_i z^{*t}_j}-\pi^{*}_{z^{*t}_i z^{*t}_j})(X_{ij}^t-\pi^{*}_{z^{*t}_i z^{*t}_j})}{\pi^{*}_{z^{*t}_i z^{*t}_j}(1-\pi^{*}_{z^{*t}_i z^{*t}_j})} \right].
\end{align*}
In the end, we decompose
\begin{align*}
\log \frac{\proba_{\breve{\theta}}(Z^{1:T}=z^{1:T}\given X^{1:T})}{\proba_{\breve{\theta}}(Z^{1:T}=z^{*1:T}\given X^{1:T})}=&\sum_{(i,j,t)\in D^*} \left( X_{ij}^t \log\frac{\pi^{*}_{z^t_i z^t_j}}{\pi^{*}_{z^{*t}_i z^{*t}_j}} + (1-X_{ij}^t) \log\frac{1-\pi^{*}_{z^t_i z^t_j}}{1-\pi^{*}_{z^{*t}_i z^{*t}_j}} \right) + \sum_{i=1}^n \log \frac{\breve{\alpha}_{z^1_i}}{\breve{\alpha}_{z^{*1}_i}} + \sum_{t=1}^{T-1} \sum_{i=1}^n \log \frac{\breve{\gamma}_{z^{t}_i z^{t+1}_i}}{\breve{\gamma}_{z^{*t}_i z^{*t+1}_i}} \nonumber\\
&+ \sum_{(i,j,t)\in D^*\cup \breve{D}}\log\left[1+\frac{(\breve{\pi}_{z^t_i z^t_j}-\pi^{*}_{z^t_i z^t_j})(X_{ij}^t-\pi^{*}_{z^t_i z^t_j})}{\pi^{*}_{z^{t}_i z^{t}_j}(1-\pi^{*}_{z^t_i z^t_j})}  \right] \nonumber\\
&-\sum_{(i,j,t)\in D^*\cup \breve{D}} \log\left[1+\frac{(\breve{\pi}_{z^{*t}_i z^{*t}_j}-\pi^{*}_{z^{*t}_i z^{*t}_j})(X_{ij}^t-\pi^{*}_{z^{*t}_i z^{*t}_j})}{\pi^{*}_{z^{*t}_i z^{*t}_j}(1-\pi^{*}_{z^{*t}_i z^{*t}_j})} \right],
\end{align*}
which gives the result.

\subsection{Proof of Lemma~\ref{minor}}
We first notice that
\[
\left| D_{n,T}(z^{1:T},\pi)\right| = \frac{1}{2}\left| \left\{ (i,j,t)\in \n^2\times \T; \pi_{z_i^t z_j^t} \neq \pi_{z_i^{*t} z_j^{*t}} \right\} \right| = \frac{1}{2}
\sum_{t=1}^T 	\left| \left\{ (i,j)\in \n^2; \pi_{z_i^t z_j^t} \neq \pi_{z_i^{*t} z_j^{*t}} \right\} \right|.
\]
For every $t \in \T$, we can apply Proposition B.4. from
\cite{celisse2012consistency}, as their Assumption (A4) is required to hold only for $z^{*t}$ (see proof) and is valid {on $\Omega_{\eta}(\theta)$} with the constant $\delta-\eta$. We obtain
\[
\left| \left\{ (i,j)\in \n^2; \pi_{z_i^t z_j^t} \neq \pi_{z_i^{*t} z_j^{*t}} \right\} \right| \geq \frac{(\delta-\eta)^2}{2}n r(t).
\]
We conclude by noticing that  $\sum_{t=1}^T r(t)=r$.

\subsection{Proof of Lemma~\ref{major}}
The inclusion of the sets is straightforward. Now we have
\begin{align*}
\left| \left\{ (i,j,t)\in \n^2\times \T; \pi_{z_i^t z_j^t} \neq \pi_{z_i^{*t} z_j^{*t}} \right\} \right| \leq &
                                                                                                                \left|  \left\{ (i,j,t) \in \n^2\times \T; (z_i^t ,z_j^t) \neq (z_i^{*t}, z_j^{*t}) \right\}  \right|\\
  \leq &\left|  \left\{ (i,j,t) \in \n^2\times \T; z_i^t  \neq z_i^{*t}\right\}  \right| + \left|  \left\{ (i,j,t) \in \n^2\times \T; z_j^t  \neq z_j^{*t}\right\}  \right|\\
\leq& 2\sum_{t=1}^T nr(t) 
\leq 2nr.
\end{align*}

\subsection{Proof of Lemma~\ref{lem:cv_gamma_intermediaire}}
First, let us decompose the quantity at stake as follows
\begin{align} \label{eq:cv_to_alpha_gamma}
&\proba_{\theta^*}\left( \left|\frac{1}{n(T-1)}\sum_{t=1}^{T-1} \sum_{i=1}^n \proba_{\hat{\theta}_{\sigma}}\left(Z^t_i=q, Z^{t+1}_i=l
\given X^{1:T}\right) - \alpha^*_q \gamma^*_{ql} \right| > \epsilon r_{n,T} \frac{\sqrt{\log n}}{\sqrt{nT}}  \right) \nonumber \\
\leq& \proba_{\theta^*}\left( \left|\frac{1}{n(T-1)}\sum_{t=1}^{T-1} \sum_{i=1}^n \proba_{\hat{\theta}_{\sigma}}\left(Z^t_i=q, Z^{t+1}_i=l \given X^{1:T}\right) - \frac{N_{ql}(Z^{1:T})}{n(T-1)} \right| > \frac{\epsilon}{2} r_{n,T} \frac{\sqrt{\log n}}{\sqrt{nT}}  \right) \nonumber \\&+ \proba_{\theta^*}\left( \left|\frac{N_{ql}(Z^{1:T})}{n(T-1)} - \alpha^*_q \gamma^*_{ql} \right| > \frac{\epsilon}{2} r_{n,T} \frac{\sqrt{\log n}}{\sqrt{nT}} \right),
\end{align}
and upper bound the two terms in the right-hand side of \eqref{eq:cv_to_alpha_gamma}. For the first one we will follow the proof of Theorem 3.9 from \cite{celisse2012consistency}. Let $z^{1:T}$ denote a fixed configuration. We work on the set $\{Z^{1:T}=z^{1:T}\}$ and write
\begin{align*}
V_1(z^{1:T})\coloneqq &\left|\frac{1}{n(T-1)}\sum_{t=1}^{T-1} \sum_{i=1}^n \proba_{\hat{\theta}_{\sigma}}\left(Z^t_i=q, Z^{t+1}_i=l \given X^{1:T}\right) - \frac{N_{ql}(z^{1:T})}{n(T-1)} \right| \\
\leq&\left|\frac{1}{n(T-1)}\sum_{t=1}^{T-1} \sum_{i=1}^n \proba_{\hat{\theta}_{\sigma}}\left(Z^t_i=q, Z^{t+1}_i=l \given X^{1:T}\right)
\ind_{(z^{t}_i, z^{t+1}_i) = (q,l)} - \frac{N_{ql}(z^{1:T})}{n(T-1)} \right|  \\
&+ \frac{1}{n(T-1)}\sum_{t=1}^{T-1} \sum_{i=1}^n \proba_{\hat{\theta}_{\sigma}}\left(Z^t_i=q, Z^{t+1}_i=l \given X^{1:T}\right) \ind_{(z^{t}_i, z^{t+1}_i) \neq (q,l)}\\
\leq& \frac{1}{n(T-1)}\sum_{t=1}^{T-1} \sum_{i=1}^n \left(1-\proba_{\hat{\theta}_{\sigma}}\left((Z^t_i, Z^{t+1}_i)=(z^{t}_i,      z^{t+1}_i) \given X^{1:T}\right)\right) \ind_{(z^{t}_i, z^{t+1}_i) = (q,l)} \\
& + \frac{1}{n(T-1)}\sum_{t=1}^{T-1} \sum_{i=1}^n \proba_{\hat{\theta}_{\sigma}}\left((Z^t_i, Z^{t+1}_i) \neq (z^{t}_i, z^{t+1}_i) \given X^{1:T}\right) \ind_{(z^{t}_i, z^{t+1}_i) \neq (q,l)}\\
\leq 2 &\proba_{\hat{\theta}_{\sigma}}\left(Z^{1:T} \neq z^{1:T} \given X^{1:T}\right).
\end{align*}
Then
\begin{align*} 
\proba_{\theta^*}\left(V_1(Z^{1:T}) >\frac{\epsilon}{2} r_{n,T} \frac{\sqrt{\log n}}{\sqrt{nT}} \right)
= & \esp_{\theta^*} \left[ \proba_{\theta^*}\left(V_1(Z^{1:T}) > \frac{\epsilon}{2} r_{n,T} \frac{\sqrt{\log n}}{\sqrt{nT}} \given[\Big] Z^{1:T} \right) \right] \\  
\leq & \sum_{z^{1:T}} \proba_{\theta^*}\left(\proba_{\hat{\theta}_{\sigma}}\left(Z^{1:T} \neq z^{1:T} \given X^{1:T}\right)  > \frac{\epsilon}{4} r_{n,T} \frac{\sqrt{\log n}}{\sqrt{nT}}  \given[\Big] Z^{1:T}=z^{1:T}\right) \proba_{\theta^*} \left(Z^{1:T}=z^{1:T}\right) \\
\leq & \sum_{z^{1:T}} \proba_{\theta^*}\left( \frac{\proba_{\hat{\theta}_{\sigma}}\left(Z^{1:T} \neq z^{1:T} \given X^{1:T}\right)}{\proba_{\hat{\theta}_{\sigma}}\left(Z^{1:T} = z^{1:T} \given X^{1:T}\right)} > \frac{\epsilon}{4} r_{n,T} \frac{\sqrt{\log n}}{\sqrt{nT}} \given[\Big] Z^{1:T}=z^{1:T}\right) \proba_{\theta^*} \left(Z^{1:T}=z^{1:T}\right) \\
\leq &  
QT\exp(-2\eta^2n) + 
\proba_{\theta^*} \left( \| \hat{\pi}_{\sigma}- \pi^* \|_{\infty} > v_{n,T} \right) \\
&
+C nT \exp \left[ - (\delta - \eta )^2 C_1 n + C_2 \log(nT) + C_4 \log\left(\frac{4 \sqrt{nT}}{\epsilon r_{n,T} \sqrt{\log n}}\right) \right]\\ &+ CnT \exp\left[ - C_3    \frac{(\log(nT))^2}{n
	v_{n,T}^2}+ 3 n \log(nT)\right] , \label{eq:cv_gamma_intermediaire} \numberthis
\end{align*} 
where the last inequality comes from Theorem~\ref{prop:ratio} where the bound is uniform with respect to $z^{1:T}$.

Now, for the second term of \eqref{eq:cv_to_alpha_gamma}, we use the following lemma. 
\begin{lemme} \label{lem:mixing} 
	There exist $c_1,c_2>0$ such that for any $\epsilon>0$, for any sequence $\{r_{n,T}\}_{n,T \geq 1}$, we have, as long as $\epsilon r_{n,T} \sqrt{\log n}/(2 \alpha^*_q \gamma^*_{ql} \sqrt{nT}) < 1$,
	\begin{align} \label{eq:mixing}
	\proba_{\theta^*}\left( \left|\frac{N_{ql}(Z^{1:T})}{n(T-1)} - \alpha^*_q \gamma^*_{ql} \right| > \frac{\epsilon}{2} r_{n,T} \frac{\sqrt{\log n}}{\sqrt{nT}} \right) \leq c_1 \exp\left( - c_2\epsilon^2 r_{n,T}^2 \right).  
	\end{align}
\end{lemme}
We then combine the two upper bounds obtained in~\eqref{eq:cv_gamma_intermediaire} and~\eqref{eq:mixing} in order to conclude, the assumption  $\epsilon r_{n,T} \sqrt{\log n}/(2 \alpha^*_q \gamma^*_{ql} \sqrt{nT}) < 1$ being satisfied for $n$ and $T$ large enough because $r_{n,T}=o(\sqrt{nT/\log n})$. 
We obtain the expected result, using the fact that $\log(T)=o(n)$, that $r_{n,T}$ increases to infinity and that $v_{n,T}=o\left(\sqrt{\log(nT)}/n\right)$,
\begin{align*}
&\proba_{\theta^*}\left( \left|\frac{1}{n(T-1)}\sum_{t=1}^{T-1} \sum_{i=1}^n \proba_{\hat{\theta}_{\sigma}}\left(Z^t_i=q, Z^{t+1}_i=l
\given X^{1:T}\right) - \alpha^*_q \gamma^*_{ql} \right| > \epsilon y_{n,T} \right) 
\leq \proba_{\theta^*} \left( \| \hat{\pi}_{\sigma}- \pi^* \|_{\infty} > v_{n,T} \right) + o(1).
\end{align*}
\qed

\subsection{Proof of Lemma~\ref{lem:variationel_asymp_eq}}
	We have the following inequalities by definition of $\hat{z}^{1:T}$, $\mathcal{J}(\chi,\theta)$ and $\hat{\chi}(\theta)$ and because the Kullback-Leibler divergence is non-negative
	\begin{align} \label{eq:ineg_variationnel}
	\mathcal{J}(\hat{z}^{1:T},\theta) \leq \mathcal{J}(\hat{\chi}(\theta),\theta) \leq \ell(\theta) \leq \ell_c(\theta,\hat{z}^{1:T}),
	\end{align}
	with $\mathcal{J}(\hat{z}^{1:T},\theta)=\ell(\theta) - KL(\delta_{\hat{z}^{1:T}}, \probat(\cdot | X^{1:T}))$.
	We write this Kullback-Leibler divergence (from $\probat(\cdot | X^{1:T})$ to $\mathbb{Q}_{\chi}=\delta_{\hat{z}^{1:T}}$, with $\chi=(\tau,\eta)$ such that $\tau_{iq}^t=\hat{z}_{iq}^t$ and $\eta^t_{iql}=\hat{z}_{iq}^{t} \hat{z}_{il}^{t+1}$) as follows
	\begin{align*}
	KL(\delta_{\hat{z}^{1:T}}, \probat(\cdot | X^{1:T})) =& - \log \probat(\hat{z}^{1:T} | X^{1:T}).
	\end{align*} 
	We then obtain
	\begin{align*}
	\mathcal{J}(\hat{z}^{1:T},\theta)=& \log \probat(X^{1:T}) + \log \probat(\hat{z}^{1:T} | X^{1:T}) = \probat(X^{1:T}|\hat{z}^{1:T}) + \log \probat(\hat{z}^{1:T}) \\
	=& \ell_c (\theta;\hat{z}^{1:T}) + \sum_{i=1}^n \log \alpha_{\hat{z}_i^1} + \sum_{i=1}^n \sum_{t=2}^T \log \gamma_{\hat{z}_i^{t-1} \hat{z}_i^{t}}.
	\end{align*}
	Combined with \eqref{eq:ineg_variationnel}, this leads to the following inequality for any parameter $\theta \in \Theta$
	\[
	\left| \mathcal{J}(\hat{\chi}(\theta),\theta) - \ell(\theta) \right| \leq \left| \mathcal{J}(\hat{z}^{1:T},\theta) - \ell_c(\theta,\hat{z}^{1:T}) \right| \leq - \sum_{i=1}^n \log \alpha_{\hat{z}_i^1} - \sum_{i=1}^n \sum_{t=2}^T \log \gamma_{\hat{z}_i^{t-1} \hat{z}_i^{t}} \leq nT \log (1/\delta).
	\]
	We can conclude that 
	\[
	\sup_{\theta \in \Theta} \left| \frac{2}{n (n-1) T} \mathcal{J}(\hat{\chi}(\theta),\theta) - \frac{2}{n (n-1) T} \ell(\theta) \right| \leq \frac{2 \log (1/\delta)}{n-1}.
	\]
	\qed

\subsection{Proof of Lemma~\ref{lem:cv_gamma_intermediate_variationnel}}
		This proof is quite similar to that of Lemma~\ref{lem:cv_gamma_intermediaire}.
		For any $\epsilon >0$, let us write 
		\begin{align*}
		&\proba_{\theta^*}\left( \left|\frac{1}{n(T-1)} \sum_{i=1}^n \sum_{t=1}^{T-1} \mathbb{Q}_{\hat{\chi}(\tilde{\theta}_{\sigma})}(Z_i^{t}=q,Z_i^{t+1}=l) - \alpha^*_q \gamma^*_{ql} \right| > \epsilon r_{n,T} \frac{\sqrt{\log n}}{\sqrt{nT}} \right) \\
		\leq& \proba_{\theta^*}\left( \left|\frac{1}{n(T-1)}\sum_{i=1}^n \sum_{t=1}^{T-1} \mathbb{Q}_{\hat{\chi}(\tilde{\theta}_{\sigma})}(Z_i^{t}=q,Z_i^{t+1}=l) - \frac{N_{ql}(Z^{1:T})}{n(T-1)} \right| > \frac{\epsilon}{2} r_{n,T} \frac{\sqrt{\log n}}{\sqrt{nT}} \right) \\&+ \proba_{\theta^*}\left( \left|\frac{N_{ql}(Z^{1:T})}{n(T-1)} - \alpha^*_q \gamma^*_{ql} \right| > \frac{\epsilon}{2} r_{n,T} \frac{\sqrt{\log n}}{\sqrt{nT}} \right)
		\end{align*}
		and upper bound the two probabilities in the right-hand side of this inequality.
		We already proved in Lemma~\ref{lem:cv_gamma_intermediaire} that the second term converges to $0$ thanks to the assumptions on the sequence $\{r_{n,T}\}_{n,T \geq 1}$.
		For the first term, let $z^{1:T}$ denote a fixed configuration. 
		Let us work on the set $\{Z^{1:T}=z^{1:T}\}$ and use the same method as in the proof of Lemma~\ref{lem:cv_gamma_intermediaire},
		\begin{multline*}
		\frac{1}{n(T-1)} \sum_{i=1}^n \sum_{t=1}^{T-1} \mathbb{Q}_{\hat{\chi}(\tilde{\theta}_{\sigma})}(Z_i^{t}=q,Z_i^{t+1}=l)=\\ \frac{1}{n(T-1)} \sum_{i=1}^n \sum_{t=1}^{T-1} \mathbb{Q}_{\hat{\chi}(\tilde{\theta}_{\sigma})}(Z_i^{t}=q,Z_i^{t+1}=l) \ind_{z^{t}_i=q,z^{t+1}_i=l} +\frac{1}{n(T-1)} \sum_{i=1}^n \sum_{t=1}^{T-1} \mathbb{Q}_{\hat{\chi}(\tilde{\theta}_{\sigma})}(Z_i^{t}=q,Z_i^{t+1}=l) \ind_{(z^{t}_i,z^{t+1}_i) \neq (q,l)},
		\end{multline*}
		leading to
		\begin{align*}
		\left| \frac{1}{n(T-1)} \sum_{i=1}^n \sum_{t=1}^{T-1} \mathbb{Q}_{\hat{\chi}(\tilde{\theta}_{\sigma})}(Z_i^{t}=q,Z_i^{t+1}=l) - \frac{N_{ql}(z^{1:T})}{n(T-1)} \right| 
		\leq& 2 \mathbb{Q}_{\hat{\chi}(\tilde{\theta}_{\sigma})}(Z^{1:T}\neq z^{1:T}).
		\end{align*}
Then we obtain
\begin{align*}
	&\proba_{\theta^*}\left( \left|\frac{1}{n(T-1)}\sum_{i=1}^n \sum_{t=1}^{T-1} \mathbb{Q}_{\hat{\chi}(\tilde{\theta}_{\sigma})}(Z_i^{t}=q,Z_i^{t+1}=l) - \frac{N_{ql}(Z^{1:T})}{n(T-1)} \right| > \frac{\epsilon}{2} r_{n,T} \frac{\sqrt{\log n}}{\sqrt{nT}} \right)\\
	\leq& \sum_{z^{1:T}} \proba_{\theta^*}\left( \mathbb{Q}_{\hat{\chi}(\tilde{\theta}_{\sigma})}(Z^{1:T}\neq z^{1:T}) > \frac{\epsilon}{4} r_{n,T} \frac{\sqrt{\log n}}{\sqrt{nT}} \given[\Big] Z^{1:T}=z^{1:T} \right) \proba_{\theta^*}\left( Z^{1:T}=z^{1:T} \right).
\end{align*}
For each $z^{1:T}$, we use the following lemma.
\begin{lemme}\label{lem:variational_distrib}
	Denoting $\tilde{\proba}_{\sigma}(\cdot)=\proba_{\tilde{\theta}_{\sigma}}(Z^{1:T}=\cdot \given X^{1:T})$, we have the following inequality for any configuration $z^{1:T}$
	\[
	\left|\mathbb{Q}_{\hat{\chi}(\tilde \theta_{\sigma})}(z^{1:T})-\tilde{\proba}_{\sigma}(z^{1:T}) \right| \leq \sqrt{ -\frac{1}{2} \log\left(\tilde{\proba}_{\sigma}(z^{1:T})\right)}.
	\]
	\end{lemme}
This gives us
\begin{align*}
&\proba_{\theta^*}\left( \mathbb{Q}_{\hat{\chi}(\tilde{\theta}_{\sigma})}(Z^{1:T}\neq z^{1:T}) > \frac{\epsilon}{4} r_{n,T} \frac{\sqrt{\log n}}{\sqrt{nT}} \given[\Big] Z^{1:T}=z^{1:T} \right) \\ 
\leq& \proba_{\theta^*}\left( \left| \mathbb{Q}_{\hat{\chi}(\tilde{\theta}_{\sigma})}(Z^{1:T}\neq z^{1:T})- \tilde{\proba}_{\sigma}(Z^{1:T}\neq z^{1:T}) \right| > \frac{\epsilon}{8} r_{n,T} \frac{\sqrt{\log n}}{\sqrt{nT}} \given[\Big] Z^{1:T}=z^{1:T} \right) \\ &+ \proba_{\theta^*}\left( \tilde{\proba}_{\sigma}(Z^{1:T}\neq z^{1:T})  > \frac{\epsilon}{8} r_{n,T} \frac{\sqrt{\log n}}{\sqrt{nT}} \given[\Big] Z^{1:T}=z^{1:T} \right)\\
\leq& \proba_{\theta^*}\left( \sqrt{ -\frac{1}{2} \log\left(\tilde{\proba}_{\sigma}(z^{1:T})\right)} > \frac{\epsilon}{8} r_{n,T} \frac{\sqrt{\log n}}{\sqrt{nT}} \given[\Big] Z^{1:T}=z^{1:T} \right) + \proba_{\theta^*}\left( \tilde{\proba}_{\sigma}(Z^{1:T}\neq z^{1:T})  > \frac{\epsilon}{8} r_{n,T} \frac{\sqrt{\log n}}{\sqrt{nT}} \given[\Big] Z^{1:T}=z^{1:T} \right)\\
\leq& \proba_{\theta^*}\left( \tilde{\proba}_{\sigma}(Z^{1:T}\neq z^{1:T}) > 1- \exp\left(- \frac{\epsilon^2 r_{n,T}^2 \log n}{32 n T}\right) \given[\Big] Z^{1:T}=z^{1:T} \right) + \proba_{\theta^*}\left( \tilde{\proba}_{\sigma}(Z^{1:T}\neq z^{1:T})  > \frac{\epsilon}{8} r_{n,T} \frac{\sqrt{\log n}}{\sqrt{nT}} \given[\Big] Z^{1:T}=z^{1:T} \right). \numberthis \label{eq:lem:cv_gamma_intermediate_variationnel}
\end{align*}
Noticing that the assumptions on $\{r_{n,T}\}_{n,T \geq 1}$ imply that 
\[
-\log\left[1- \exp\left(- \frac{\epsilon^2 r_{n,T}^2 \log n}{32 n T}\right)\right]=o(n) \quad\textrm{ and }\quad -\log\left[r_{n,T} \frac{\sqrt{\log n}}{\sqrt{nT}}\right]=o(n),
\]
we can conclude by applying the result of Theorem~\ref{prop:ratio} with the estimator $\tilde{\theta}_{\sigma}=(\tilde{\Gamma}_{\sigma},\tilde{\pi}_{\sigma})$ for both terms of the right-hand side of~\eqref{eq:lem:cv_gamma_intermediate_variationnel}.
\qed

\subsection{Proof of Lemma~\ref{lemmebound}}
	The proof follows the lines of the proof of Lemma C.3. from~\cite{celisse2012consistency}.
	Let $\esp^*_{\theta^*}[\cdot]$ denote the expectation given $Z^{1:T}=z^{*1:T}$,
	i.e. $\esp^*_{\theta^*}[\cdot]=\espetoile[\cdot\given Z^{1:T}=z^{*1:T}]$. 
	Introducing a ghost sample $\{\tilde X^{t}_{ij}\}_{i,j,t}$ that is independent of $\{X^{t}_{ij}\}_{i,j,t}$ and has the same
	distribution, we  write
	\begin{align*} 
	E \coloneqq &\espetoile^*\left[ \sup_{(z^{1:T},\pi)\in \Q^{nT}\times [\zeta,1-\zeta]^{Q^2}} \left|\frac{2}{n(n-1)T} \sum_{t=1}^T  \sum_{i<j} (\Xtij-\pi^{*}_{z^{*t}_i z^{*t}_j})
	\log\left(\frac{\pi_{z^t_i z^t_j}}{1-\pi_{z^t_i z^t_j}}\right) \right| \right] \\
	= & \espetoile^*\left\{ \sup_{(z^{1:T},\pi) \in \Q^{nT}\times [\zeta,1-\zeta]^{Q^2}} \left|\frac{2}{n(n-1)T} \espetoile^*\left[ \sum_{t=1}^T  \sum_{i<j} (\Xtij-\tilde X^t_{ij}) \log\left(\frac{\pi_{z^t_i z^t_j}}{1-\pi_{z^t_i z^t_j}} \right) \given[\Big] \{X_{ij}^t\}_{i,j,t} \right] \right| \right\}\\
	\leq & \espetoile^*\left\{ \espetoile^*\left[ \sup_{(z^{1:T},\pi) \in \Q^{nT}\times [\zeta,1-\zeta]^{Q^2}} \frac{2}{n(n-1)T}  \left|\sum_{t=1}^T  \sum_{i<j} (\Xtij-\tilde X^t_{ij}) \log\left(\frac{\pi_{z^t_i z^t_j}}{1-\pi_{z^t_i z^t_j}} \right) \right| \given[\Big] \{X_{ij}^t\}_{i,j,t} \right] \right\}\\
	\leq & \esp^*_{\theta^*, X,\tilde X}\left[ \sup_{(z^{1:T},\pi) \in \Q^{nT}\times [\zeta,1-\zeta]^{Q^2}} \frac{2}{n(n-1)T} \left|\sum_{t=1}^T  \sum_{i<j} (\Xtij-\tilde X^t_{ij}) \log\left(\frac{\pi_{z^t_i z^t_j}}{1-\pi_{z^t_i z^t_j}} \right)\right| \right],
	\end{align*}
	where $\esp^*_{\theta^*, X,\tilde X}[\cdot]$ denotes the expectation with respect to $\{X,\tilde X\}=\{X_{ij}^t, \tilde X_{ij}^{t}\}_{i,j,t}$ under the true parameter $\theta^*$ and given $Z^{1:T}=z^{*1:T}$.
	At this point, we notice that, if $\{\epsilon^{t}_{ij}\}_{i,j,t}\coloneqq\epsilon$ are $n^2T$ independent Rademacher variables, then the random variables
	\[
	\esp_{\epsilon}  \left|\sum_{t=1}^T  \sum_{i<j} \epsilon_{ij}^t (\Xtij-\tilde X^t_{ij})
	\log\left(\frac{\pi_{z^t_i z^t_j}}{1-\pi_{z^t_i z^t_j}} \right)\right| 
	\quad \text{ and } \quad \left|\sum_{t=1}^T
	\sum_{i<j} (\Xtij-\tilde X^t_{ij}) \log\left(\frac{\pi_{z^t_i z^t_j}}{1-\pi_{z^t_i z^t_j}} \right)\right| 
	\]
	follow the same distribution, which implies that 
	\begin{multline*}
	\esp^*_{\theta^*, X,\tilde X}\left[ \sup_{(z^{1:T},\pi) \in \Q^{nT}\times [\zeta,1-\zeta]^{Q^2}} \frac{2}{n(n-1)T} \esp_{\epsilon}  \left|\sum_{t=1}^T  \sum_{i<j} \epsilon_{ij}^t (\Xtij-\tilde
	X^t_{ij}) \log\left(\frac{\pi_{z^t_i z^t_j}}{1-\pi_{z^t_i z^t_j}} \right)\right| \right] 
	\\= \esp^*_{\theta^*, X,\tilde X} \left[
	\sup_{(z^{1:T},\pi) \in \Q^{nT}\times [\zeta,1-\zeta]^{Q^2}} \frac{2}{n(n-1)T} \left|\sum_{t=1}^T  \sum_{i<j} (\Xtij-\tilde X^t_{ij}) \log\left(\frac{\pi_{z^t_i
			z^t_j}}{1-\pi_{z^t_i z^t_j}} \right)\right| \right]. 
	\end{multline*}
	As a consequence, we have
	\begin{align*}
	E 
	\leq& \esp^*_{\theta^*, X,\tilde X}\left[ \sup_{(z^{1:T},\pi) \in \Q^{nT}\times [\zeta,1-\zeta]^{Q^2}} \frac{2}{n(n-1)T} \esp_{\epsilon}  \left|\sum_{t=1}^T \sum_{i<j} \epsilon_{ij}^t (\Xtij-\tilde X^t_{ij}) \log\left(\frac{\pi_{z^t_i z^t_j}}{1-\pi_{z^t_i z^t_j}} \right)\right| \right] \\
	\leq& \espetoile^*\left[ \sup_{(z^{1:T},\pi) \in \Q^{nT}\times [\zeta,1-\zeta]^{Q^2}} \frac{2}{n(n-1)T} \esp_{\epsilon}  \left|\sum_{t=1}^T  \sum_{i<j} \epsilon_{ij}^t
	\Xtij \log\left(\frac{\pi_{z^t_i z^t_j}}{1-\pi_{z^t_i z^t_j}} \right)\right|  \right]
	\\ &+ \espetoile^*\left[ \sup_{(z^{1:T},\pi) \in \Q^{nT}\times [\zeta,1-\zeta]^{Q^2}} \frac{2}{n(n-1)T} \esp_{\epsilon} \left|\sum_{t=1}^T  \sum_{i<j} \epsilon_{ij}^t \tilde X^t_{ij} \log\left(\frac{\pi_{z^t_i z^t_j}}{1-\pi_{z^t_i z^t_j}} \right)\right| \right]\\
	\leq& 2 \espetoile^*\left[ \sup_{(z^{1:T},\pi) \in \Q^{nT}\times [\zeta,1-\zeta]^{Q^2}} \frac{2}{n(n-1)T} \esp_{\epsilon} \left|\sum_{t=1}^T  \sum_{i<j} \epsilon_{ij}^t \Xtij \log\left(\frac{\pi_{z^t_i z^t_j}}{1-\pi_{z^t_i z^t_j}} \right)\right| \right].
	\end{align*}
	Then using Jensen's inequality, Assumption~\ref{itm:hyp_pi} and the bound $\var_{\epsilon}(\epsilon_{ij}^t \Xtij) \le 1$, we get 
	\begin{align*}
	E
	\leq & 2 \espetoile^*\left[ \sup_{(z^{1:T},\pi) \in \Q^{nT}\times [\zeta,1-\zeta]^{Q^2}} \frac{2}{n(n-1)T} \sqrt{ \esp_{\epsilon}\left[ \left(\sum_{t=1}^T  \sum_{i<j} \epsilon_{ij}^t \Xtij \log\left(\frac{\pi_{z^t_i z^t_j}}{1-\pi_{z^t_i z^t_j}} \right)\right)^2 \right]} \right]\\
	\leq & 2 \espetoile^*\left[ \sup_{(z^{1:T},\pi) \in \Q^{nT}\times [\zeta,1-\zeta]^{Q^2}} \frac{2}{n(n-1)T} \sqrt{ \var_{\epsilon} \left[ \sum_{t=1}^T  \sum_{i<j} \epsilon_{ij}^t \Xtij \log\left(\frac{\pi_{z^t_i z^t_j}}{1-\pi_{z^t_i z^t_j}} \right) \right]} \right]\\
	\leq& 2 \espetoile^*\left[ \frac{2}{n(n-1)T} \sup_{\pi \in [\zeta,1-\zeta]}\log\left( \frac{\pi} {1-\pi} \right) \sqrt{\frac{n(n-1)T}{2}} \right]
	\leq \sqrt{\frac{2}{n(n-1)T}} \Lambda,
	\end{align*}
	where $\Lambda=2 \log[(1-\zeta)/\zeta]$, concluding the proof.
	\qed

\subsection{Proof of Lemma~\ref{lem:fct_M}}
We assume that $\min_{\sigma \in \mathfrak{S}_Q}\|\pi_{\sigma}-\pi^*\|_{\infty} >\epsilon$. Without loss of generality, assume that the permutation (or one of the permutations) minimizing this distance is the identity. Let us write, using the fact that $I_Q$ the identity matrix of size $Q$ maximizes in $A$ (over the set of $Q\times Q$ stochastic matrices) the quantity $\mathbb{M}(\pi^*,A)$ (see the proof of Theorem 3.6 in \cite{celisse2012consistency}) and denoting $(\bar{a}_{qq'})_{q,q' \in \Q}$ the coefficients of  $\bar{A}_{\pi}$ (thus depending on $\pi$),
\[
\mathbb{M}(\pi^*)-\mathbb{M}(\pi)= \sum_{q,l} \alpha_q^* \alpha_l^* \sum_{q',l'} \bar{a}_{qq'} \bar{a}_{ll'} \left[ \pi^*_{ql} \log \frac{\pi^*_{ql}}{\pi_{q'l'}} + (1-\pi^*_{ql}) \log \frac{1-\pi^*_{ql}}{1-\pi_{q'l'}} \right] = \sum_{q,l} \alpha_q^* \alpha_l^* \sum_{q',l'} \bar{a}_{qq'} \bar{a}_{ll'} K(\pi^*_{ql},\pi_{q'l'})
\]
denoting $K(p_1,p_2)=p_1 \log (p_1/p_2) + (1-p_1) \log [(1-p_1)/(1-p_2)]>0$ the Kullback-Leibler divergence from a Bernoulli distribution with parameter $p_2$ to a Bernoulli distribution with parameter $p_1$. For every $q$, there exists $q'\coloneqq f(q)$ such that $\bar{a}_{qq'} \geq 1/Q$ because $\bar{A}_{\pi}$ is a stochastic matrix.
Using Assumption~\ref{itm:hyp_gamma}, we obtain
\[
\mathbb{M}(\pi^*)-\mathbb{M}(\pi) \geq \frac{\delta^2}{Q^2} \sum_{q,l} K(\pi^*_{ql},\pi_{f(q) f(l)}) \geq \frac{\delta^2}{Q^2} \sum_{q,l} 2 (\pi^*_{ql}-\pi_{f(q) f(l)})^2
\]
thanks to a result on Kullback-Leibler divergence for Bernoulli distributions (see for instance \cite{phdbubeck}, Chapter 10, Section 2, Lemma 10.3).
We then want to show that there exist $q,l$ such that $|\pi^*_{ql}-\pi_{f(q) f(l)}|>\epsilon$. 
\begin{itemize}
	\item If $f$ is a permutation, the assumption $\min_{\sigma \in \mathfrak{S}_Q}\|\pi_{\sigma}-\pi^*\|_{\infty} >\epsilon$ gives the expected result. 
	\item If $f$ is not a permutation, it is not injective and there exist $q_1\neq q_2$ such that $f(q_1)=f(q_2)$. Thanks to Assumption~\ref{itm:hyp_id}, take $l_0 \in \Q$ such that $|\pi_{q_1l_0}-\pi_{q_2l_0}|=\max_{l \in \Q} |\pi_{q_1l}-\pi_{q_2l}|>0$. Then 
	\[
	|\pi^*_{q_1 l_0}- \pi_{f(q_1) f(l_0)}| + |\pi_{f(q_2) f(l_0)} - \pi^*_{q_2 l_0}| \geq |\pi^*_{q_1 l_0}- \pi_{f(q_1) f(l_0)} + \pi_{f(q_2) f(l_0)} - \pi^*_{q_2 l_0}| = |\pi^*_{q_1 l_0}-\pi^*_{q_2 l_0}|>0 
	\]
	leading to either $|\pi^*_{q_1 l_0}- \pi_{f(q_1) f(l_0)}|\geq |\pi^*_{q_1 l_0}-\pi^*_{q_2 l_0}| /2>\epsilon$ or  $|\pi^*_{q_2 l_0}- \pi_{f(q_2) f(l_0)}|\geq |\pi^*_{q_1 l_0}-\pi^*_{q_2 l_0}| /2>\epsilon$, using the fact that $\epsilon< \min_{1\leq q\neq q' \leq Q} \max_{1\leq l\leq Q} |\pi^*_{ql}-\pi^*_{q'l}|/2$.
\end{itemize}
So, as there exist $q,l$ such that $|\pi^*_{ql}- \pi_{f(q)f(l)}|> \epsilon$, we have
\[
\mathbb{M}(\pi^*)-\mathbb{M}(\pi) > \frac{2\delta^2}{Q^2} \epsilon^2.
\]
\qed
\subsection{Proof of Lemma \ref{lem:mixing}}
For any node $i\in \n$, the Markov chain $\{Z_i^t\}_{t\ge 1}$ is geometrically ergodic because its transition matrix $\Gamma$
satisfies Doeblin's condition thanks to Assumption~\ref{itm:hyp_gamma}. 
For any $z\in \Q$, let us denote $\delta_{z}$ the Dirac mass at $z$. 
There exists a positive constant $A$ and some $r \in
(0,1)$ such that $\forall q \in\Q$ and $\forall t \geq 1$, we have  
\[
\left\| \delta_q\Gamma^t-\alpha \right\|_{TV} \leq Ar^t, 
\]
where $\|\cdot \|_{TV}$ is the total variation norm. This leads to 
\[
\left\| \delta_q\Gamma^t-\alpha \right\|_{TV} = \frac{1}{2} \left\| \delta_q\Gamma^t-\alpha \right\|_{1} = \frac{1}{2} \sum_{l \in \Q} | \Gamma^t(q,l)- \alpha_{l} | \leq Ar^t.
\]
We now consider the Markov chain $\{Z^t=(Z_1^t,\dots, Z_n^t)\}_{t\ge 1}$ of the $n$ nodes evolving through time. Note
that it is irreducible and aperiodic. Moreover, its
transition matrix is given by $P_n= \Gamma^{\otimes n}$, the $n$-th Kronecker power of $\Gamma$ and its stationary
distribution is $\alpha^{\otimes n}$.  
For any $z=(z_1,\ldots,z_n) \in \Q^n$, let us denote $\mu_{n,z}= \otimes_{i=1}^n \delta_{z_i}$.  
For every $t \geq 1$, we can decompose
\begin{align*}
\left\| \mu_{n,z} P_n^t-\alpha^{\otimes n} \right\|_{TV} &=\left\| \left(\overset{n}{\underset{i=1}{\otimes}} \delta_{z_i}\right) (\Gamma^{\otimes n})^t-\alpha^{\otimes n}\right\|_{TV}
=\left\| \left(\overset{n}{\underset{i=1}{\otimes}} \delta_{z_i}\right) (\Gamma^t)^{\otimes n}-\alpha^{\otimes n}\right\|_{TV}
=\left\| \overset{n}{\underset{i=1}{\otimes}} \left(\delta_{z_i}  \Gamma^t \right)-\alpha^{\otimes n}\right\|_{TV}\\
&= \frac{1}{2} \left\| \overset{n}{\underset{i=1}{\otimes}}\left( \delta_{z_i}  \Gamma^t \right)-\alpha^{\otimes n}\right\|_{1}
= \frac{1}{2} \sum_{(z'_1,\ldots,z'_n)\in \Q^n} \left| \prod_{i=1}^n \Gamma^t(z_i,z_i')- \prod_{i=1}^n \alpha_{z'_i} \right|.
\end{align*}
We use 
\[
\prod_{i=1}^n \Gamma^t(z_i,z_i')- \prod_{i=1}^n \alpha_{z'_i}= \sum_{i=1}^n
\left\{ \left( \prod_{j=1}^{i-1} \alpha_{z'_j} \right) \left[ \Gamma^t(z_i,z_i')-\alpha_{z'_i}\right]  \prod_{k=i+1}^{n}
(\mu_{z_k}\Gamma^t)_{z_k'}\right\}. 
\]
So, reorganizing the terms, we write
\begin{align*}
\left\| \mu_{n,z} P_n^t-\alpha^{\otimes n}\right\|_{TV} 
&\leq \frac{1}{2} \sum_{(z'_1,\ldots,z'_n)\in \Q^n} \sum_{i=1}^n \left\{ \left(\prod_{j=1}^{i-1} \alpha_{z'_j} \right)\left|  \Gamma^t(z_i,z_i')-\alpha_{z'_i} \right| \prod_{k=i+1}^{n}  (\mu_{z_k}\Gamma^t)_{z_k'} \right\}\\
&\leq \frac{1}{2} \sum_{i=1}^n \sum_{z'_1} \alpha_{z'_1}   \ldots \sum_{z'_{i-1}} \alpha_{z'_{i-1}} \sum_{z'_i} \left| \Gamma^t(z_i,z_i')-\alpha_{z'_i} \right| \sum_{z'_{i+1}}\Gamma^t(z_{i+1},z_{i+1}')\ldots \sum_{z'_n} \Gamma^t(z_n,z_n')\\
&\leq \frac{1}{2} \sum_{i=1}^n \sum_{z'_i} \left|\Gamma^t(z_i,z_i')-\alpha_{z'_i} \right|
\leq n A r^t. 
\end{align*}
Let us recall the definition of an $\epsilon$-mixing time. For any Markov transition matrix $M$ over the set
$\mathcal{X}$ with stationary distribution $\alpha$, for any $\epsilon >0$, the $\epsilon$-mixing time of the Markov chain is defined as 
\[
\tau(\epsilon) = \min \{ t \ge 1; \max_{x \in \mathcal{X}} \|\delta_xM^t -\alpha\|_{TV} \le \epsilon\}.
\]
Denoting by $\tau_n(\epsilon)$ the $\epsilon$-mixing time of the Markov chain $\{Z^t\}_{t \ge 1}$, we thus obtain 
\[
\tau_n(\epsilon) \leq \frac{\log(nA/\epsilon)}{\log(1/r)}.
\]
Now, we introduce a new Markov chain $Y=\{Y^t\}_{t\geq1}$, that is defined by
\[
Y^t=(Z^t,Z^{t+1}) \quad \forall t \geq 1.
\]
Notice that it is irreducible and aperiodic, with stationary
distribution $\rho$ defined  for every state $(q_1^t,\ldots,q_n^t,q_1^{t+1},\ldots,q_n^{t+1})$ by 
\begin{align*}
\rho_{(q_1^t,\ldots,q_n^t,q_1^{t+1},\ldots,q_n^{t+1})}=\alpha_{q_1^t} \ldots \alpha_{q_n^t} \gamma_{q_1^{t}q_1^{t+1}} \ldots \gamma_{q_n^{t}q_n^{t+1}}.
\end{align*}
It is easily seen that for any $\epsilon >0$,  its $\epsilon$-mixing time $\tau_{Y,n}(\epsilon)$ equals
$\tau_n(\epsilon)+1$. 
We apply Theorem 3 from~\cite{chung_et_al:LIPIcs:2012:3437}, for any $\eta \leq 1/8$, considering the weight function $f(Y^t)=\sum_{i=1}^n$ for every $t\geq 1$ (of expectation $n \alpha^*_q \gamma^*_{ql}$ under the stationary distribution). Then $N_{ql}(Z^{1:T})=\sum_{t=1}^{T-1} f(Y^t)$, and denoting $\epsilon_{n,T}=\epsilon r_{n,T} \sqrt{\log n}/(2 \alpha^*_q \gamma^*_{ql} \sqrt{nT})$, we obtain that there exist
$c_1,c_2>0$ such that for any $\epsilon>0$, as long as $\epsilon_{n,T} \leq 1$
\begin{align*}
\proba_{\theta^*}\left( \left|\frac{N_{ql}(Z^{1:T})}{n(T-1)} - \alpha^*_q \gamma^*_{ql} \right| >\frac{\epsilon}{2} r_{n,T} \frac{\sqrt{\log n}}{\sqrt{nT}} \right) 
=& \proba_{\theta^*}\left( N_{ql}(Z^{1:T}) > (1+\epsilon_{n,T}) n \alpha^*_q \gamma^*_{ql} (T-1) \right) \\ &+ \proba_{\theta^*}\left( N_{ql}(Z^{1:T}) < (1-\epsilon_{n,T}) n \alpha^*_q \gamma^*_{ql} (T-1) \right)\\
\leq& c_1 \exp\left( - \frac{\epsilon_{n,T}^2 n \alpha^*_q \gamma^*_{ql} (T-1)}{72 \tau_{Y,n}(\eta)}\right) \leq c_1 \exp\left( - c_2\epsilon^2 r_{n,T}^2 \right).
\end{align*}
\qed

\subsection{Proof of Lemma~\ref{lem:variational_distrib}}
For any configuration $z^{1:T}$,
	\begin{align*}
	\left|\mathbb{Q}_{\hat{\chi}(\tilde{\theta}_{\sigma})}(z^{1:T})-\tilde{\proba}_{\sigma}(z^{1:T}) \right| \leq \left\| \mathbb{Q}_{\hat{\chi}(\tilde{\theta}_{\sigma})}- \tilde{\proba}_{\sigma} \right\|_{TV} \leq \sqrt{ \frac{1}{2} KL(\mathbb{Q}_{\hat{\chi}(\tilde{\theta}_{\sigma})},\tilde{\proba}_{\sigma})} \leq \sqrt{ \frac{1}{2} KL(\delta_{z^{1:T}},\tilde{\proba}_{\sigma})} \leq \sqrt{ -\frac{1}{2} \log\left(\tilde{\proba}_{\sigma}(z^{1:T})\right)},
	\end{align*}
	the third inequality being true because by definition $\mathbb{Q}_{\hat{\chi}(\tilde{\theta}_{\sigma})}$ minimizes $KL(\cdot,\tilde{\proba}_{\sigma})$ over the set of variational distributions. 
	\qed

\end{document}